\documentclass{amsart}
\usepackage{amssymb,stmaryrd}
\usepackage{amsfonts}
\usepackage{amstext}
\usepackage{algorithmic}
\usepackage{algorithm}
\usepackage{graphicx}
\usepackage{epstopdf}
\usepackage[all]{xy}
\usepackage{MnSymbol}

\numberwithin{equation}{section}
\newtheorem{theorem}{Theorem}[section]
\newtheorem{lemma}[theorem]{Lemma}
\newtheorem{proposition}[theorem]{Proposition}
\newtheorem{corollary}[theorem]{Corollary}

\theoremstyle{definition}

\newtheorem{example}[theorem]{Example}

\theoremstyle{remark}

\newcommand{\R}{{\mathbb{R}}}

\newcommand{\C}{{\mathbb{C}}}

\renewcommand{\(}{{(}}

\newcommand{\CC}{{\mathcal{C}}}
\newcommand{\CE}{{\mathcal{E}}}
\newcommand{\CD}{{\mathcal{D}}}

\newcommand{\CR}{{\mathcal{R}}}

\newcommand{\wedgeq}{{\wedge\kern-5pt\cdot}}

\newcommand{\tens}{\otimes}

\newcommand{\id}{{\rm id}}

\newcommand{\extd}{{\rm d}}
\newcommand{\del}{{\partial}}
\newcommand{\eps}{\epsilon}

\newcommand{\la}{{\triangleright}}

\begin{document}

\title{Semiquantisation Functor and Poisson-Riemannian Geometry, I}


\author{Edwin J.\ Beggs \&\ Shahn Majid}
\address{Dept of Mathematics, Swansea University\\ Singleton Parc, Swansea SA2 8PP\\
 +\\ 
School of Mathematical Sciences, Queen Mary, University of London\\
 Mile End Rd, London E1 4NS, UK}

\email{e.j.beggs@swansea.ac.uk, s.majid@qmul.ac.uk}


\begin{abstract} We study noncommutative bundles and Riemannian geometry at the semiclassical level of first order in a deformation parameter $\lambda$, using a functorial approach. The data for quantisation of the cotangent bundle is  known to be a Poisson structure and Poisson preconnection and we now show that this data defines to a functor $Q$ from the monoidal category of classical vector bundles equipped with connections  to the monodial category of bimodules equipped with bimodule connections over the quantised algebra. We adapt this functor to quantise the wedge product of the exterior algebra and in the Riemannian case, the metric and the Levi-Civita connection. Full metric compatibility requires vanishing of an obstruction in the classical data, expressed in terms of a generalised Ricci 2-form $\CR$, without which our quantum Levi-Civita connection is still the best possible. We apply the theory to the Schwarzschild black-hole and to Riemann surfaces as examples, as well as verifying our results on the 2D bicrossproduct model quantum spacetime. The quantized Schwarzschild black-hole in particular has features similar to those encountered in $q$-deformed models, notably  the necessity of nonassociativity of any rotationally invariant quantum differential calculus of classical dimensions. 
\end{abstract}
\maketitle 

\section{Introduction}

Noncommutative geometry aims to extend notions of geometry to situations where the `coordinate algebra' is noncommutative. Such algebras could arise on `quantisation' in the passage 
from a classical mechanical system to a quantum one or, it is now widely accepted, as a more accurate description of spacetime itself so as to include quantum corrections arising out
of quantum gravity. There are different approaches to the formulation of the right axioms in the noncommutative case and we mention notably the approach of Connes\cite{Con} coming out of cyclic cohomology,  ring-theoretic projective module approach due to Van den Bergh,  Stafford and others, e.g. \cite{vdb}, and a constructive approach coming out of quantum groups but not limited to them in which the different layers of geometry starting with the differential graded algebra are built up, typically guided by quantum symmetry and other considerations. Our work falls within this last approach and particularly a set of axioms of `noncommutative Riemannian geometry' using bimodule connections\cite{Mou,DV1,DV2,BegMa2,BegMa3,BegMa:twi,BegMa4,Ma:gra,Ma:alm}. We refer in particular to \cite{BegMa4} for a nontrivial 2D example containing a gravitational source and worked out in complete detail in this bimodule connection approach.

In this paper we take  a step towards the general problem of `quantisation' of all the rest of the geometrical structure beyond the algebra itself, within the bimodule connections approach. We consider only perturbative phenomena in the sense of order by order in a deformation parameter $\lambda$ and hence will miss `deep quantum' effects, but on the plus side this will allow us to construct concrete examples from familiar classical geometrical data and moreover, by explicitly building structures on the classical counterparts, we will have a ready-made identification between these and corresponding quantum objects which would otherwise be open to interpretation such as ordering ambiguities. This then provides a route to making experimental predictions. Thus, while the deformation problem is by no means adequate for the full theory of noncommutative geometry, it  nevertheless has practical value. 

In fact we are only going to solve here the problem of quantising Riemannian and other differential geometry to first order in $\lambda$. This will already be quite a significant task as we shall see,  but at this level we will arrive at a fairly complete and functorial picture (we would then envision to develop the same ideas order by order possibly in an $A_\infty$ algebra approach). Physically speaking, the minimal hypothesis is that noncommutative geometry represents an effective description of spacetime to include first planck-scale correction and in that case it may only be the first order in $\lambda$ that are immediately relevant (it is tempting to suppose an absolute significance to the noncommutative geometry but that is a further assumption).  Also, in distance units the value of $\lambda$  at around $10^{-35}$m is extremely small making these effects only just now beginning to be measurable in principle, in which case $O(\lambda^2)$ effects can be expected to be so much even smaller as to be beyond any possible relevance. This motivates a deeper analysis of the semiclassical level where we work to order $\lambda$. Mathematically speaking, this means that we are working at the Poisson level and in principle we could set our our main results such as the functor in Theorem~\ref{functor} in an entirely Poisson setting. However, we will develop the theory equivalently in a novel way that keeps better the connection with algebra, namely as exact noncommutative geometry but working over the ring of dual numbers $\C[\lambda]/(\lambda^2)$ where we set $\lambda^2=0$.  We call the construction of such noncommutative geometry more precisely {\em semiquantisation} rather than the more familiar term `semiclassicalisation'.

Starting with a classical manifold $M$ expressed algebraically as $C^\infty(M)$, it is well-known that a Poisson bivector $\omega^{ij}$ controls its associative deformations and moreover that the `quantisation problem' is then solved notably by Kontsevich\cite{Konts} at a formal level and in the symplectic case and in the presence of a symplectic connection,  more geometrically by Fedosov\cite{Fed}. These results are, however, only for the algebra and for actual noncommutative geometry we must `quantise' also differential structures, bundles, connections, and so forth. By differential structure we mean the algebra $\Omega(M)$ of differential forms and in \cite{Haw,BegMa1} it was shown that deformations of the 1-forms $\Omega^1(M)$ are controlled to this order by Poisson-preconnections $\hat\nabla$ on $\Omega^1(M)$ (these are only defined along hamiltonian vector fields). The curvature of the preconnection entails nonassociativity of the bimodule product (so a breakdown of $(a\extd b)c=a((\extd b)c)$ for elements $a,b,c$ of the noncommutative algebra) at $O(\lambda^2)$.  More precisely, for convenience, we will assume an actual connection $\nabla$ with torsion but in most cases we only make use of the combination $\omega^{ij}\nabla_j$ which could more generally be replaced by a preconnection. Our starting point,  implicit in \cite{BegMa1}, is that Poission-compatibility in terms of torsion amounts to the condition for $(\omega,\nabla)$,
\begin{equation}\label{poissoncompat} \omega^{ij}_{; m}=\omega^{ki}T^j_{km}+\omega^{jk}T^i_{km}
 \end{equation}
where $T$ is the torsion of $\nabla$ (see Lemma~\ref{compatT}). The startling conclusion in Sections 3,4 of the paper is that with this assumption we have canonically a functor $Q$ that quantises to lowest order the monoidal category of classical bundles over $M$ (Theorem~\ref{functor}) and also the structure of the exterior algebra $(\Omega(M),\extd)$, i.e. the wedge product (and tensor products more generally), see Theorem~\ref{dga}. 

Then in Section 5 we quantise a Riemannian structure consisting of metric $g$ and Levi-Civita connection $\widehat{\nabla}$, which will require a further condition on the classical geometry (i.e. not every classical Riemannian manifold will be quantisable even at our first order level). We tend to see that as a good rather than a bad thing, as ultimately constraints of this type could (when the theory is fully developed) explain such things as Einstein's equation. The thinking is that if classical geometry emerges out of quantum gravity via noncommuative geometry then the greater rigidity of the latter can and should imprint constraints on what can emerge at the classical level and hence explain them. Thus in the 2D toy model in \cite{BegMa4} the constraints of noncommutative geometry forces a curved metric on the chosen spacetime algebra and this metric describes either a strongly gravitational source at the origin in space or a toy model of a big-bang cosmology with fluid matter, depending on the interpretation and sign of a parameter. Mathematically speaking, the semiclassical analysis allows us to clarify  obstructions or `anomalies' to the quantisation process and one can take the view just stated that vanishing of the anomaly is a quantisability constraint or, if we wish to consider more general models, one can take the view that either the anomaly is an order $\lambda$ effect to live with or that it is a signal that the theory needs to be extended, for example extra dimensions, to absorb the anomaly. 

Our quantum wedge product $\wedge_1$ in Theorem~\ref{dga} consists of a functorial part $\wedge_Q$ plus an order $\lambda$ nonlinear correction that can be attributed to the non-linear nature of the problem in enforcing the Leibniz rule (the calculus is both used in the functor and now is being quantised by the functor). As a result the quantum metric $g_1=g_Q-\lambda\CR$ consists of a functorial part $g_Q$ plus an order $\lambda$  correction give by a certain 2-form $\CR$ (which we call the `generalised Ricci 2-form') viewed here as an antisymmetric tensor, see Proposition~\ref{nablaQCR}. The only requirement at this point beyond (\ref{poissoncompat}) is 
\begin{equation}\label{nablag} \nabla g=0\end{equation}
 which we assume throughout in order to be able to apply our functorial methods. This system (\ref{poissoncompat})-(\ref{nablag}) is already  a strong constraint which, by the time one adds symmetry requirements, explains the anomaly for differential calculus in quantum group models \cite{BegMa1,BegMa:coch} as forced curvature of $\nabla$, as well as the need for gravity in \cite{BegMa4} as forced curvature of the metric. For the Schwarzschild black-hole metric we will find a similar anomaly to the one for quantum group models, namely that $\nabla$ has to have curvature. Thus any rotationally invariant deformation of the black-hole will necessarily entail nonassociativity at $O(\lambda^2)$ if we assume classical dimensions (an alternative is the use of an  extra cotangent direction  in \cite{Ma:alm}). This is a tangible result of our semiclassical analysis in the present paper. 
   
Similarly, both the Poisson-compatible connection $\nabla$ and the classical Levi-Civita connection $\widehat\nabla$ get functorially quantised as $\nabla_Q$ and $\nabla_{QS}$ respectively, where we write $\widehat\nabla=\nabla+S$ and quantise each term functorially. As before, the functorial $\nabla_{QS}$ needs an order $\lambda$ correction, i.e. we construct the quantum Levi-Civita connection in the form $\nabla_1=\nabla_{QS}+\lambda K$ for some tensor $K$ which is uniquely determined by requiring $\nabla_1$ to be quantum torsion free and requiring merely the symmetric part of $\nabla_1g_1$ to vanish; see Theorem~\ref{qlevi} and our `quantum Koszul' formula (\ref{koszul}). The antisymmetric part of $\nabla_1 g_1$, however, is a new phenomenon which does not exist classically; an order $\lambda$ obstruction independent of $K$ and proportional to the left hand side of the equation
\begin{equation}\label{CRcompat} \widehat\nabla \CR +\omega^{ij}\,g_{rs}\,S^s_{jn}(R^r{}_{mki}+S^r_{km;i})\,\extd x^k\tens\extd x^m \wedge \extd x^n  =0,\end{equation}
where $R$ is the curvature of $\nabla$. This additional (\ref{CRcompat}) is therefore necessary and sufficient for the existence of a torsion free fully metric compatible $\nabla_1$, and when this does exist  it is given by our above unique $\nabla_1$. Otherwise the latter remains `best possible' in the sense of killing the part of $\nabla_1 g_1$ that can be controlled, with the antisymmetric part remaining as an anomaly. The Schwarzschild black-hole will again be similar to the quantum group case in that there will be this order $\lambda$ obstruction to full metric compatibility.

While the above holds formally over most fields, for physics we want to work over $\R$ at the classical level while at the quantum level over $\C$ but in a $*$-algebra setting where the classical reality is extended as Hermiticity (we will say `reality constraint' as an umbrella terms for the relevant constraint but one could also say loosely `unitarity'). Thus our quantum metric $g_1$ will be complex but subject to such a `reality' constraint  and similarly we would like our quantum connection $\nabla_1$ to be  suitably `real' in the sense of star-preserving\cite{BegMa3}. Again our functorial construction $\nabla_{QS}$ works (Theorem~\ref{starpresK}) in that there is always a unique order $\lambda$ adjustment $K$ to make it star-preserving, leading to a  canonical  $\nabla_1$ from this point of view. Fortunately, Corollary~\ref{jointcompat} says that when a quantum torsion free metric compatible connection exists it is necessarily star-preserving and coincides with the $\nabla_1$ given by this reality/unitarity requirement, so our two constructions coincide in this case. Even when a fully metric compatible one does not exist, our unique star-preserving $\nabla_1$ still tends to be the same as our unique `best possible' quantum Levi-Civita found before, at least in nice cases that we have looked at such as the Schwarzschild black-hole. 

An alternative to the above straight metric compatibility is\cite{BegMa3} to work  with the corresponding sequilinear `Hermitian' metric $(\star\tens\id)g_1$ and ask for $\nabla_1$ to be quantum torsion free and Hermitian-metric compatible. Here again there is no obstruction and we show, Proposition~\ref{hermetriccompat}, that we can always find suitable $K$. In general these $\nabla_1$ will not be unique but in some cases there could be a unique one. For example, in the geometry of quantum groups  there is a unique perturbative such connection  for $C_q(SU_2)$ with its 3D calculus  in \cite[Thm. 7.9]{BegMa3}.  Another weaker notion of metric compatibility is vanishing  cotorsion -- meaning $(\wedge_1\tens\id)\nabla_1 g_1=0$ and in nice cases such as the standard $q$-sphere\cite{Ma:sph} quantum  torsion free, cotorsion free turns out to be the same as quantum torsion free, Hermitian-metric compatible. We will see that this is also the case for the Schwarzschild black-hole.
 
We remark that one could view $\wedge_1,g_1,\nabla_1$ above as the functor $Q$ applied to preadjusted nonstandard classical maps $\wedge', g',\widehat\nabla'=\nabla+S'$ differing from the usual ones  by an order $\lambda$ correction, thus $g'=g-\lambda\CR$ etc. One could then say that to be quantised functorially the classical metric needs to acquire an antisymmetric order $\lambda$ correction (in our conventions $\lambda$ is imaginary so $g'$ remains Hermitian) and similarly for $\wedge'$ and $S'$. This interpretation could be a direct path to predictions, although it is not in our scope to pursue that point of view here.

Section 6 turns to examples in the simplest case where $S=0$, i.e. where the quantising connection and the Levi-Civita connection coincide. Their quantisation $\nabla_{Q}=\nabla_{QS}$ does not need any order $\lambda$ adjustment as it is automatically star-preserving and `best possible' in terms of metric compatibility. The constraint  (\ref{poissoncompat}) simplifies to $\nabla\omega=0$, (\ref{nablag}) is automatic,  while the condition (\ref{CRcompat})  for full metric-compatibility simplifies to $\nabla\CR=0$. The latter is automatically solved for example if $(\omega,g)$ is K\"ahler-Einstein. It is similarly solved for any surface of constant curvature and we give hyperbolic space and sphere in detail. This section thus provides the simplest class of solutions. The only downside is that the Levi-Civita connection typically has significant curvature in examples of interest and if we take this for our quantizing connection then the quantum differential calculus will be nonassociative at $O(\lambda^2)$. 

Finally, in Section~7 we given two examples where the quantising connection $\nabla$ is very different from the Levi-Civita one. The first is the 2D bicrossproduct model quantum spacetime with curved metric in \cite{BegMa4} but analysed now at the semiclassical level. Here all our conditions  (\ref{poissoncompat})-(\ref{CRcompat}) hold, $\nabla$  has zero curvature but a lot of torsion, and application of the general machinery indeed yields a unique torsion free metric-compatible and star-preserving quantum connection in agreement with one of the two quantum Levi-Civita connections in\cite{BegMa4} (the other is `deep quantum' with no $\lambda\to 0$ limit). This provides a nontrivial check on our analysis. We then turn to the Schwarzschild black hole metric with a rotationally invariant $\omega$. The latter is not covariantly constant so $\nabla$ cannot be the Levi-Civita one, and rather we find a 4-functional parameter moduli of rotationally invariant $\nabla$ with torsion. These cannot be adjusted to have zero curvature so there will be nonassociativity at $O(\lambda^2)$ as explained above. Also as promised, there is a nonzero obstruction to the antisymmetric part of the metric compatibility and turns out to be `topological' in the sense of independent even of the choice of $\nabla$ within the considered moduli. Nevertheless, we find the unique `best possible' quantum Levi-Civita which turns out also to be the unique star-preserving $\nabla_1$.

While our focus above has been on the quantum Levi-Civita connection, our functorial quantisation $\nabla_Q$ of any quantising connection $\nabla$ obeying (\ref{poissoncompat})--(\ref{nablag}) also has nice properties by itself (see Proposition~\ref{nablaQCR}). Here $\nabla_Q$ is quantum star-preserving and, if $\nabla\CR=0$, quantum metric compatible but will generally have quantum torsion. This quantum connection  has a more direct relevance to teleparallel gravity\cite{tele} where, when the manifold is parallelisable, one can take $\nabla$ to be the Weitzenb\"ock connection instead of the Levi-Civita one that we took in Section~6. This has torsion but zero curvature and working with it is equivalent to General Relativity but interpreted differently, with $S$ above now viewed as the contorsion tensor and our results similarly viewed as its quantisation $Q(S)+\lambda K$. The Weitzenb\"ock $\nabla$ with its zero curvature corresponds to an associative quantum differential calculus at $O(\lambda^2)$, but as the black-hole example showed one may need to have some small amount of curvature to have a compatible $\omega$, i.e. to be quantisable. This application, as well as the general theory of the quantum Ricci tensor,  quantum Laplacians and quantum  complex structures in the sense of \cite{ebsmNCcomplex} at order $\lambda$ are deferred to forthcoming work, as is the important case of  $\C P^n$ as an example of a K\"ahler-Einstein manifold, the classical Riemannian geometry of which is linked to Berry phase and higher uncertainty relations in quantum mechanics\cite{Brody}, among other applications.

Also, we note that the equation (\ref{poissoncompat}) has a striking similarity to weak-metric-compatibility \cite{Ma:recon}
\[ g^{ij}_{; m}=g^{ki}T^j_{km}+g^{jk}T^i_{km}.\]
which applies to metric-connection pairs arising from cleft central extensions of the classical exterior algebra by an extra closed 1-form $\theta'$ with $\theta'{}^2=0$ and $\theta'$ graded-commutative, whereas in our present paper we extend by $\lambda$ a central scalar with $\lambda^2=0$ as explained above. The $\theta'$ approach was used to associatively quantise the Schwarzschild black-hole in \cite{Ma:alm} in contrast to our approach now. It would seem that these two different ideas might be unified into a single  construction. This and the higher order theory are some other directions for further work.

\section{Preliminaries}
\subsection{Classical differential geometry}  \label{kuycvgy}

We assume that the reader is comfortable with classical differential geometry and recall its noncommutative algebraic generalisation in a bimodule approach. For classical geometry suffice it to say that we assume $M$ is a smooth manifold with further smooth structures notably the exterior algebra $(\Omega(M),\extd)$ but more generally we could start with any graded-commutative classical differential graded algebra with further structure (i.e. the graded-commutative case of the next section). One small generalisation: we allow complexifications. However, one could work with real values and a trivial $*$-operation or one could have the a complex version and, in the classical case, pick out the real part. We use the following categories based on vector bundles on $M$:

\begin{tabular}{ccc}Name & Objects & Morphisms \\$\CE_0$  & vector bundles over $M$ & bundle maps \\$\tilde{\CD_0}$ & $(E,\nabla)$ bundle and connection & bundle maps \\$\CD_0$  & $(E,\nabla)$ bundle and connection & bundle maps intertwining the connections\end{tabular}

The condition for the bundle map $\theta:E\to F$ to  intertwine the connections $(E,\nabla_E)$ and $(F,\nabla_F)$ is that
$\theta(\nabla_{Ei}e)=\nabla_{Fi}(\theta(e))$, where $e:M\to E$ is a section of the bundle $E$. To fit the viewpoint of noncommutative geometry, we will talk about sections of the bundles rather than the bundles themselves, and from that point of view we would write $e\in E$, and consider $E$ as a module over the algebra of smooth functions $C^\infty(M)$.

We use two tensor products: $E\tens F$ will denote the algebraic tensor product, and $E\tens_0 F$ will denote the tensor product over $C^\infty(M)$. Thus $E\tens_0 F$ obeys the relation $e.a\tens_0 f=e\tens_0 a.f$ for all $a\in C^\infty(M)$, and this corresponds to the usual tensor product of vector bundles, wheras $E\tens F$ is much larger. We can use the tensor product to rewrite a connection $(E,\nabla_E)$ as a map $\nabla_E:E\to\Omega^1 (M)\tens_0 E$ by using the formula $\nabla_E(e)=\extd x^i\tens_0\nabla_i(e)$, and the Leibniz rule becomes
\begin{eqnarray*}
\nabla_E(a.e) \,=\, \extd a\tens_0 e + a.\nabla_E(e)\ .
\end{eqnarray*}

As most readers will be more familiar with tensor calculus on manifolds  than with the commutative case of the algebraic version above, we use the former throughout for computations in the classical case. We adopt here standard conventions for curvature and torsion tensors as well Christoffel symbols for a linear connection. On forms and in a local coordinate system we have
\[ \nabla_j\extd x^i=-\Gamma^i_{jk}\extd x^k\]
while $T_\nabla=\wedge\nabla-\extd:\Omega^1(M) \to \Omega^2(M)$ has the torsion tensor
\begin{eqnarray} \label{vuyihgfxxy}
T_\nabla(\extd x^i)= -\Gamma^i_{jk}\extd x^j\wedge\extd x^k=\tfrac12\,T^i_{jk}\extd x^k\wedge\extd x^j,\quad T^i_{jk}=\Gamma^i_{jk}-\Gamma^i_{kj}\ .
\end{eqnarray}
Similarly, for the curvature tensor
\[ R_\nabla(\extd x^k)=\tfrac12\,\extd x^i\wedge\extd x^j\tens_0[\nabla_i,\nabla_j]\extd x^k=\tfrac12\,R^k{}_{mij}\extd x^j\wedge\extd  x^i \tens_0  \extd x^m.\]
The summation convention is understood unless specified otherwise. 

We also recall the interior product $\righthalfcup\,:\mathrm{Vec}M\tens \Omega^n(M)\to \Omega^{n-1}(M)$ defined by $v\, \righthalfcup\,\eta$ being the evaluation for $\eta\in\Omega^1(M)$, or in terms of indices $v^i\,\eta_i$, extended recursively to higher degrees by
\begin{eqnarray*}
v\, \righthalfcup\,(\xi\wedge\eta)=(v\, \righthalfcup\,\xi)\wedge\eta+(-1)^{|\xi|}\xi\wedge(v\, \righthalfcup\,\eta)\ .
\end{eqnarray*}

\subsection{Noncommutative bundles and connections}

Here we briefly summarise the elements of noncommutative differential geometry that we will be concerned with in our bimodule approach\cite{DV1,DV2,BegMa3,BegMa4}. The following picture can be generalised at various places, but for readability we will not refer to this further. 
 The associative algebra $A$ (over the complex numbers) plays the role of `functions' on our noncommutative space and need not be commutative.

A differential calculus on $A$ consists of $n$ forms $\Omega^n A$ for $n\ge 0$, an associative product
$\wedge:\Omega^nA \tens \Omega^mA\to\Omega^{n+m}A$ and an \textit{exterior derivative} $\extd:
\Omega^nA\to \Omega^{n+1}A$ satisfying the rules

1)\quad $\Omega^0A=A$ (i.e.\ the zero forms are just the `functions')\newline
2)\quad $\extd^2=0$\newline
3)\quad $\extd(\xi \wedge\eta)=\extd\xi \wedge\eta+(-1)^{|\xi|}\,\xi\wedge\extd\eta$ where $|\xi|=n$ if $\xi\in\Omega^nA$\newline
4) $\Omega$ is generated by degree 0,1.

These are the rules for a standard \textit{differential graded algebra}. Note that we \textit{do not} assume graded commutativity, which would be $\xi\wedge\eta=(-1)^{|\xi|\,|\eta|}\,\eta\wedge\xi$.

A vector bundle is expressed as a (projective) $A$-module. If $E$ is a left $A$-module we define a left connection $\nabla_{E}$ on $E$  to be a map $\nabla_{E}:E\to \Omega^1A \tens_A  E$
obeying the left Leibniz rule
\[ \nabla_{E}(ae)=\extd a \tens_A  e+a\nabla_{E}(e).\]
We say that we have a bimodule connection if $E$ is a bimodule and there is a bimodule map
 \[ \sigma_{E}: E \tens_A  \Omega^1A \to \Omega^1A \tens_A  E,\quad \nabla_{E}(ea)=(\nabla_{E}e)a+\sigma_{E}(e \tens_A \extd a).\]
 If $\sigma_{E}$ is well-defined then it is uniquely determined, so its existence is a property of a left connection on a bimodule. There is a natural tensor product of
 bimodule connections $(E,\nabla_{E})\tens(F,\nabla_{F})$ built on the tensor product $E \tens_A  F$ and
 \[ \nabla_{E \tens_A  F}(e \tens_A  f)=\nabla_{E}e \tens_A  f+ (\sigma_{E}\tens\id)(e \tens_A \nabla_{F}f).\]
 There is necessarily an associated $\sigma_{E \tens_A  F}$. We denote by $\CE$ the monoidal category of $A$-bimodules with $ \tens_A $. We denote by $\CD$ the monoidal category of pairs $(E,\nabla_E)$ of bimodules and bimodule connections over $A$. Its morphisms are bimodule maps that intertwine the connections. 
 
 In the case of a connection on $E=\Omega^nA$ we define the torsion as $T_\nabla=\wedge\nabla-\extd$. In the case $E=\Omega^1A$ we define a metric as $g\in \Omega^1A\tens_A\Omega^1A$ and
 in case of a bimodule connection the metric-compatibility tensor is $\nabla g\in \Omega^1A\tens_A\Omega^1A\tens_A\Omega^1A$. We also require $g$ to have an inverse $(\ ,\ ):\Omega^1A\tens_A\Omega^1A\to A$ with the usual bimodule map properties and this requires that $g$ is central. 
 
\subsection{Conjugates and star operations} \label{bcviosiuc}
 We suppose that $A$ is a star algebra, i.e.\ that there is a conjugate linear map $a\mapsto a^*$ so that
 $(ab)^*=b^*a^*$ and $a^{**}=a$. Also suppose that this extends to a star operation on the forms, so that
 $\extd(\xi^*)=(\extd\xi)^*$ and $(\xi\wedge\eta)^*=(-1)^{|\xi||\eta|}\eta^*\wedge\xi^*$.

 Next,  for any $A$-bimodule $E$, we consider its conjugate bimodule $\overline{E}$ with elements denoted by $\overline{e}\in \overline{E}$, where $e\in E$ and new right and left actions of $A$,
$\overline{e}a=\overline{a^*e}$ and $a\overline{e}=\overline{ea^*}$.  There is a canonical bimodule map $\Upsilon:\overline{E\tens_A F}\to \overline{F}\tens_A \overline{E}$ given by $\Upsilon(\overline{e\tens f})=\overline{f}\tens \overline{e}$. Also, if $\phi:E\to F$ is a bimodule map we have a bimodule map $\bar\phi:\bar E\to \bar F$ by  $\bar\phi(\bar e)=\overline{\phi(e)}$.  These constructions are examples of a general notion of a bar category\cite{BegMa3}  but for our purposes the reader should view the conjugate notation as a useful way to keep track of conjugates for noncommutative geometry, and as a book-keeping device to avoid problems. It allows, for example, conjugate linear functions to be viewed as linear functions to the conjugate of the original map's codomain. Bimodules form a bar category as explained and so does the category of pairs $(E,\nabla)$. Here $\bar E$ acquires a right handed connection $\bar\nabla(\bar e)=(\id\tens\star^{-1})\Upsilon\overline{\nabla e}$ which we convert to a left connection
 \[ \nabla\bar e= (\star^{-1}\tens\id)\Upsilon\overline{\sigma^{-1}\nabla e}.\]  Here $\star:\Omega^1A \to \overline{\Omega^1A}$ is the $*$-operation viewed formally as a linear map.

In general we say that $E$ is a $\star$-object if there is a linear operation $\star:E\to \overline{E}$ (which we can also write as $\star(e)=\overline{e^*}$ where $e\mapsto e^*$ is antilinear) such that $\bar\star\star(e)=\overline{\overline e}$ for all $e\in E$.  Also given $\star$-objects $E$, $F$ we say a morphism $\phi:E\to F$ is $*$-preserving $\bar\phi$ commutes with $\star$.  If $E$ is a star-object then we define a connection as $*$-preserving if
 \[ (\id\tens\star)\nabla(\star^{-1}\overline e)=(\star^{-1}\tens\id)\Upsilon\overline{\sigma^{-1}\nabla e}\]
and in this case $(E,\nabla)$ becomes a $\star$-object in this bar category. Cleary $\Omega^1A$ itself is an example of a star-object and so is $\Omega A$ in every degree. The product $\wedge$ is an example of an anti-$*$-preserving map (i.e. with a minus sign) on products of degree 1.

In the $*$-algebra case we say a metric $g\in \Omega^1A\tens_A\Omega^1A$ is `real' in the sense \begin{eqnarray*}
\Upsilon^{-1}(\star\tens\star)\,g=\overline{g}
\end{eqnarray*}
If $g_{ij}$ is real symmetric as a matrix valued function, then as the phrase `reality property' suggests, this is true classically. We can also work with general metrics equivalently
as `Hermitian metrics' $G=(\star\tens\id)g\in \overline{\Omega^1A}\tens_A\Omega^1A$ and this is `real' precisely when $\Upsilon^{-1}(\id\tens\mathrm{bb})G=\overline{G}$, where $\mathrm{bb}:
\Omega^1A \to \overline{\overline{\Omega^1A}}$ is the `identity' map to the double conjugate $\xi\mapsto \overline{\overline{\xi}}$. 
  In this context it is more natural to formulate metric compatibility using
the `Hermitian-metric compatibility tensor'
\begin{eqnarray} \label{hermcompat}
(\bar\nabla\tens\id+\id\tens\nabla)G\in \overline{\Omega^1A}\tens_A\Omega^1A\tens_A\Omega^1A\ .
\end{eqnarray}
If  $\nabla$ is $*$-preserving, then vanishing of this coincides with the regular notion of metric compatibility of the corresponding $g$.

\subsection{Bundles with extended morphisms}

We will sometimes want to refer to extended morphisms and their  covariant derivative, which will be particularly needed when we come to 
study quantum curvature in a sequel. For the moment suffice it to say that if $A$ is an algebra with DGA $\Omega A$, the category
of pairs $(E,\nabla_E)$ where $E$ is a left $A$-module equipped with a  left covariant derivatives, has an extended notion of morphism
as a  left module map $\theta:E\to \Omega^nA\tens_A F$ for any degree $n$. We say that $\theta\in {\rm Mor}_n(E,F)$, where the set of extended morphisms between two objects is now a graded vector space.  Composition of such an extended morphism with another, $\phi:F\to \Omega^mA\tens_A G$, is given by the following formula
\begin{eqnarray*}
\phi\circ\theta\,=\,(\id\wedge\phi)\theta:E\to \Omega^{n+m}A\tens_A G 
\end{eqnarray*}

\begin{proposition}\label{bcisbvc}
If $\theta:E\to \Omega^nA\tens_A F$ is an extended morphism from $(E,\nabla_E)$ to $(F,\nabla_F)$, then so is
\begin{eqnarray*}
\nabla(\theta)\,=\,\nabla^{[n]}_F\circ \theta-(\id\wedge\theta)\nabla_E:
E\to \Omega^{n+1}A\tens_A F
\end{eqnarray*} 
where $\nabla_F^{[n]}=\extd\tens\id+(-1)^{n}\id\wedge\nabla_F:\Omega^n(A)\tens_A F\to \Omega^{n+1}A\tens_A F$.
\end{proposition}
\proof  For $a\in A$ and $e\in E$,
and setting $\theta(e)=\xi\tens f$,
\begin{eqnarray*}
\nabla^{[n]}_F\circ \theta(a.e)-(\id\wedge\theta)\nabla_E(a.e) &=& 
\nabla^{[n]}_F(a.\xi\tens f)-(\id\wedge\theta)(\extd a\tens e+a.\nabla_E(e)) \cr
&=& \extd a\wedge \xi\tens f+
 a.\nabla^{[n]}_F(\xi\tens f)-\extd a\wedge\theta( e)+a.\nabla_E(e) \cr
 &=& a.(\nabla^{[n]}_F\circ \theta(e)-(\id\wedge\theta)\nabla_E(e))\ .\quad\square
\end{eqnarray*}

\begin{example}
The $A$-bimodule $A$ can be given the usual connection $\extd$. 
Consider the bimodule $\Omega^1 A$
with left connection $\nabla$, and the morphism
$\tau\in \mathrm{Mor}_1(\Omega^1 A,A)$ given by 
$\xi\mapsto \xi\tens_A 1$. 
Then $\nabla(\tau)$ is given by
\begin{eqnarray*}
\nabla(\tau)(\xi) &=& \extd\xi\tens 1 -(\id\wedge\tau)\nabla\,=\,(\extd\xi-\wedge\nabla\xi)\tens 1\ .
\end{eqnarray*}
This means that $-\nabla(\tau)\in  \mathrm{Mor}_2(\Omega^1 A,A)$ is the torsion of the connection on $\Omega^1A$.
\end{example}

\subsection{Imposing $\lambda^2=0$}\label{order2}

We will be working in the setting of a typically noncommutative algebra $A_\lambda$ and related structures expanded in a formal power series in $\lambda$ but truncated to different orders of approximation. It is intuitively clear what this means but one way to make it precise is as follows. 

For a field $k$, let $\CC$ be a $k$-linear Abelian category and let $k_n=k[\lambda]/(\lambda^{n+1})$. The quotient here simply means that we set $\lambda^{n+1}=0$.  For $V\in \CC$ we let $V[n]=k_n\tens_k V$ (so this consists on $n$ copies of $V$ labelled by powers of $\lambda$). 
The category $\CC[n]$ consists of such objects with morphisms those of $\CC$ extended $\lambda$-linearly to become linear over $k_n$. If $n>m$ there is a functor $\pi:\CC[n]\to \CC[m]$ given by the quotient $\lambda^{m+1}=0$ defining a map $k_n\to k_m$. If $\CC$ is monoidal then so is $\CC[n]$ with $V[n]\tens W[n]=(k_n\tens_k V)\tens_{k_n}(k_n\tens _kW)\cong k_n\tens_k (V\tens W)= (V\tens W)[n]$ and this is such that we have  a monoidal functor $\CC\to \CC[n]$ for any $n$. We denote by  $\tens_n$ the tensor product in the category $\CC[n]$.

In our case we are interested in categories where the underlying objects are vector spaces, so $\CC={\rm Vec}$. Let $A\in {\rm Vec}$ be an associative algebra. Deforming $A$ to first order then means
equipping $A_1=A[1]$ with an associative  product $A_1\tens_1 A_1\to A_1$ so that $A_1/(\lambda)=A$.  We will use a subscript 1 on categories (typically related to $A_1$)  to indicate that we are working in the deformed theory to order $\lambda$. Thus $\CE_1$ denotes the category of $A_1$-bimodules over $\C[\lambda]/(\lambda^2)$, and $\CD_1$ the category of pairs $(E,\nabla)$ as in Section~2.1 but over $A_1$. Aside from $A_1$ we will not explicitly denote the change of base on objects, for clarity.  

We are going to work over $k=\C$ with suitable reality conditions but it should be clear that constructions that do not depend on the $*$-involution work with care over most fields. 

\section{Semiquantization of bundles} 

This section constructs a monoidal functor $Q$ that quantises geometric data on a smooth manifold $M$ to first order in a deformation parameter $\lambda$. Here $A=C^\infty(M)$ is our initial algebra and its first order quantisation for us means a map $A_1\tens_1 A_1\to A_1$ as explained in Section~\ref{order2}. However, for readability purposes we will also continue to speak in more conventional terms of powerseries in $\lambda$ with errors $O(\lambda^2)$ being ignored. In an application where $\lambda$ was actually a number, the dropping of these higher powers would need to be justified by the physics.

\subsection{Quantising the algebra and modules}

The data we suppose is an antisymmetric  bivector $\omega$ on $M$ along with a linear connection $\nabla$  subject to the following `Poisson compatibility' \cite{BegMa1} 
\begin{eqnarray} \label{poissoncomp}
\extd(\omega^{ij})  - \omega^{kj}\,\nabla_k(\extd x^i)
-\omega^{ik}\, \nabla_k(\extd x^j) \,=\,0\ .
\end{eqnarray}

\begin{lemma}\label{compatT} Let $\omega$ be an antisymmetric bivector and $\nabla$ a linear connection, with torsion tensor $T$. Then  $\omega$ obeys (\ref{poissoncomp}) {\em if and only if}
\[ \omega^{ij}_{; m}+\omega^{ik}T^j_{km}-\omega^{jk}T^i_{km}=0.\]
In this case $\omega$ is a Poisson tensor if and only if
\[ \sum_{\textrm{cyclic}\ (i,j,k)}\omega^{im}\,\omega^{jp}\,T^k_{mp}=0.\]
\end{lemma}
\proof The first part is essentially in \cite{BegMa1} but given here more generally. For the first part the explicit version of (\ref{poissoncomp}) in terms of Christoffel symbols is
\begin{equation}\label{compGamma} 
 \omega^{ij}{}_{,m}+\omega^{kj}\Gamma^i_{km}+\omega^{ik}\Gamma^j_{km}=0. \end{equation}
We write the expression on the left as 
\[\omega^{ij}{}_{,m}+\omega^{kj}\Gamma^i_{mk}+\omega^{ik}\Gamma^j_{mk}+\omega^{kj}T^i_{km}+\omega^{ik}T^j_{km}\]
and we recognise the first three terms as the covariant derivative. For the second part, we put (\ref{compGamma}) into the following condition for a Poisson tensor:
\begin{eqnarray}\label{cyclic}
\sum_{\textrm{cyclic}\ (i,j,k)} \omega^{im}\,\omega^{jk}_{\phantom{jk},m}\,=\,0\ .\quad\square
\end{eqnarray}
 
\smallskip
For example, any manifold with a torsion free connection and $\omega$ a covariantly constant antisymmetric bivector will do.  This happens for example
in the case of a K\"ahler manifold, so our results include these. 

The action of the bivector on a pair of functions is denoted $\{\ ,\ \}$ as usual. If $\omega$ is a Poisson tensor then this is a Poisson bracket and from the Fedosov and Kontsevich  there is an associative multiplication
for functions
\begin{eqnarray} \label{fed12}
a\bullet b &=& ab+\lambda\{a,b\}/2+
O(\lambda^2).\ 
\end{eqnarray}
We take the same formula in any case and denote by $A_\lambda$ any (possibly not associative) quantisation with this leading order part, which means we fix our associative algebra $A_1$ over $\C[\lambda]/(\lambda^2)$ and leave higher order unspecified. We will normally assume that $\omega$ is a Poisson tensor because that will be desirable at higher order, but strictly speaking the results in the present paper do not really require this. 

Similarly, in \cite{BegMa1} we found the commutator of a function $a$ and a 1-form
$\xi\in\Omega^1(M)$
\begin{eqnarray}\label{reppo}
[a,\xi]_\bullet\,=\,\lambda\,\omega^{ij}\, a_{,i}\, (\nabla_j\xi)
+O(\lambda^2)\ ,
\end{eqnarray}
so we could define the deformed product of a function $a$ and a 1-form
$\xi$ as
\begin{eqnarray}
a\bullet \xi\,=\,a\,\xi+\lambda\,\omega^{ij}\, a_{,i}\, (\nabla_j\xi)/2
+O(\lambda^2)\ ,\cr
\xi\bullet a\,=\,a\,\xi-\lambda\,\omega^{ij}\, a_{,i}\, (\nabla_j\xi)/2
+O(\lambda^2)\ .
\end{eqnarray}
Again we can drop the corrections and  regard these as defining a bimodule structure $\Omega^1A_1\tens_1 A_1\to \Omega^1A_1$ and $A_1\tens_1\Omega^1A_1\to \Omega^1A_1$ where $\Omega^1A_1$ in this context is over $\C[\lambda]/(\lambda^2)$. 

Now let $(E,\nabla_E)$ be a classical bundle and covariant derivative on it, and define, for $e\in E$,
\begin{eqnarray}  \label{nbdhksav}
a\bullet e\,=\,a\,\xi+\tfrac{\lambda}{2} \,\omega^{ij}\, a_{,i}\, (\nabla_{Ej}e)
+O(\lambda^2)\ ,\cr
e\bullet a\,=\,a\,\xi-\tfrac{\lambda}{2} \,\omega^{ij}\, a_{,i}\, (\nabla_{Ej}e)
+O(\lambda^2)\ .
\end{eqnarray}
A brief check reveals that the following associative laws hold to  errors in $O(\lambda^2)$:
\begin{eqnarray}
\quad (a\bullet b)\bullet e\,=\,a\bullet (b\bullet e)\ ,\ (a\bullet e)\bullet b\,=\,a\bullet (e\bullet b)\ ,\ 
(e\bullet a)\bullet b\,=\,e\bullet (a\bullet b)\ ,
\end{eqnarray}
so we have a bimodule structure $E\tens_1 A_1\to E$ and $A_1\tens_1 E\to E$. We consider the following categories of modules over $A_1$:

\begin{tabular}{ccc}Name & Objects & Morphisms \\$\tilde\CE_1$ & bimodules over $A_1$ & left module maps \\$\CE_1$ & bimodules over $A_1$ & bimodule maps \\
$\CD_1$  & bimodules and connection & bimodule maps intertwining the connections
\end{tabular}

\begin{lemma} \label{bcuoisiutcxsty}   We define the functor $Q:\tilde{\CD_0}\to \tilde\CE_1$ sending
objects to objects according to $(\ref{nbdhksav})$ and sending bundle maps $T:E\to F$ to left module maps 
\[ Q(T)=T+\tfrac{\lambda}{2} \,\omega^{ij}\,\nabla_{F_i}\circ\nabla_j(T),\]
where $\nabla_j(T)= \nabla_{Fj}\circ T -T\circ\nabla_{Ej}$ as explained in the Preliminaries. The functor restricts to $Q:\CD_0\to \CE_1$ as $Q(T)=T$. In general
we have
\begin{eqnarray*}
Q(T\circ S) &=& Q(T)\circ Q(S) +\tfrac{\lambda}{2} \,\omega^{ij}\,\nabla_{i}( T) \circ \nabla_{j}(S) \ . 
\end{eqnarray*}
\end{lemma}
\proof Take $T_0:E\to F$ a bundle map. We aim for the bimodule properties
\begin{eqnarray}
(T_0+\lambda\,T_1)(a\bullet e) &=& a\bullet(T_0+\lambda\,T_1)(e)\ ,\cr
(T_0+\lambda\,T_1)(e\bullet a) &=& (T_0+\lambda\,T_1)(e)\bullet a\ ,
\end{eqnarray}
which to errors in $O(\lambda^2)$ is
\begin{eqnarray}
T_0(a\bullet e)+\lambda\,T_1(a\,e)  &=& a\bullet T_0(e)+ \lambda\,a\,T_1(e)\ ,\cr
T_0(e\bullet a)+\lambda\,T_1(e\,a)  &=& T_0(e)\bullet  a+ \lambda\,T_1(e)\,a\ .
\end{eqnarray}
Using the formula (\ref{nbdhksav}) for the deformed product gives our conditions as
\begin{eqnarray} \label{bhuivuygc}
T_0(\omega^{ij}\, a_{,i}\, (\nabla_{Ej}e)/2)+T_1(a\,e)  &=&
\omega^{ij}\, a_{,i}\, (\nabla_{Fj}T_0(e))/2  + a\,T_1(e)\ ,\cr
-\,T_0(\omega^{ij}\, a_{,i}\, (\nabla_{Ej}e)/2)+T_1(e\,a)  &=& -\,\omega^{ij}\, a_{,i}\, (\nabla_{Fj}T_0(e))/2+ T_1(e)\,a\ .
\end{eqnarray}
It is not possible to satisfy both parts of (\ref{bhuivuygc}) unless $T_0$ preserves the covariant derivatives, i.e.\ 
\begin{eqnarray}
\nabla_{Fj}T_0(e)\,=\, T_0(\nabla_{Ej}e)\ 
\end{eqnarray}
and in this case we set $T_1=0$ as a solution and $Q(T_0)=T_0$.

More generally, we solve only the first part of (\ref{bhuivuygc}), i.e. a left module map for $(A_1,\bullet)$, which needs
\begin{eqnarray} \label{bhuivuyxsqgc}
T_1(a\,e) -a\,T_1(e) &=&
\omega^{ij}\, a_{,i}\, \big(\nabla_{Fj}T_0(e)  -T_0(\nabla_{Ej}e)\big)/2 \ ,
\end{eqnarray}
Define a module map $\nabla_j(T_0):E\to F$ by
\begin{eqnarray}
\nabla_j(T_0)(e) \,=\, \nabla_{Fj}T_0(e)  -T_0(\nabla_{Ej}e)
\end{eqnarray}
now we have
the following choice:
\begin{eqnarray}
T_1\,=\, \omega^{ij}\,\nabla_{F_i}\circ \nabla_j(T_0)/2
\end{eqnarray}
which solves (\ref{bhuivuyxsqgc}). This gives  $Q(T_0)$. For compositions,
\begin{eqnarray*}
Q(T\circ S) &=& T\circ S+\tfrac{\lambda}{2} \,\omega^{ij}\,\nabla_{i}\circ \nabla_{j}(T \circ S) \cr
&=& T\circ S+\tfrac{\lambda}{2} \,\omega^{ij}\,\nabla_{i}\circ (\nabla_{j}(T) \circ S+T \circ \nabla_{j}(S)) \cr
&=& Q(T)\circ S +\tfrac{\lambda}{2} \,\omega^{ij}\,\nabla_{i}\circ T \circ \nabla_{j}(S) \cr
&=& Q(T)\circ S +\tfrac{\lambda}{2} \,\omega^{ij}\,\nabla_{i}( T) \circ \nabla_{j}(S) 
 +\tfrac{\lambda}{2} \,\omega^{ij}\,T \circ \nabla_{i}\circ  \nabla_{j}(S)      \ .\quad\square
\end{eqnarray*}

Note that we distinguish between $T:Q(E)\to Q(F)$ which is the map
$Q(e)\mapsto Q(T(e))$ and $Q(T):Q(E)\to Q(F)$ which is the given quantisation.

\subsection{Quantising the tensor product}
As we have now described how to deform the algebra and bimodules, we can now take the ``fiberwise" tensor product of two bimodules in the deformed case. This is taken to be $Q(E)\tens_1 Q(F)$, where similarly to the definition of $\tens_0$ in Section~\ref{kuycvgy}, we take $e\tens_1 a\bullet f=
e\bullet a\tens_1 f$ for all $a\in A_1$. Note that as $Q$ is the identity on objects, we could have written $E\tens_1 F$ above, but we used the $Q$ to emphasise that the bimodules are taken with the deformed actions. Now we seem to have two ways to quantise the tensor product $E\tens F$, but these are related by a \textit{natural transformation} $q$ in (\ref{vcusvdsavb}), this being the definition of $Q$ being a \textit{monoidal functor}:
\begin{eqnarray}\label{vcusvdsavb}
\xymatrix{
E\tens F \ar@{.>}[r]^{Q\circ\tens_0} \ar@{.>}[rd]_{Q\tens_1 Q} 
&Q(E\tens_0 F) \\
&
Q(E)\tens_1 Q(F) \ar[u]_{q_{E,F}}    }
\end{eqnarray}
We require the following diagram to commute:
\begin{eqnarray*}\label{vcusvdsavbfea}
\xymatrix{
Q(E)\tens_1 Q(F)\tens_1 Q(G) 
\ar[r]^{\id\tens_1 q_{F,G}} \ar[rd]_{q_{E,F}\tens_1 \id} & 
Q(E)\tens_1 Q(F\tens_0G)   \ar[r]^{q_{E,F\tens_0 G}}  & 
Q(E\tens_0F\tens_0G)  \\
 & Q(E\tens_0F)\tens_1 Q(G)    \ar[ru]_{q_{E\tens_0F,G}} 
  }
\end{eqnarray*}

\begin{proposition} \label{q} The functor $Q:\CD_0\to \CE_1$ is monoidal with associated natural transformation $q:Q\tens_1 Q\implies Q\circ\tens_0$ given by
\begin{eqnarray*}
q_{V,W}(Q(v)\tens_1 Q(w)) &=& Q(v\tens_0 w)+ \tfrac{\lambda}{2} Q(\omega^{ij}\,\nabla_{Vi} v
\tens_0 \nabla_{Wj} w)\ .
\end{eqnarray*}
In general we have, for $T:E\to V$ and $S:F\to W$,
\begin{eqnarray*}
q_{V,W}\,(T\tens_1\id_W) &=& (T\tens \id_W+\tfrac\lambda2\,\omega^{ij}\,\nabla_i(T)\tens\nabla_{Wj})\,q_{E,W}\ ,\cr
q_{V,W}\,(\id_V\tens_1 S) &=& (\id_V\tens S+\tfrac\lambda2\,\omega^{ij}\,\nabla_{Vi}\tens\nabla_j(S))\,q_{V,F}\ .
\end{eqnarray*}
As usual, we have $\nabla_i(T)=\nabla_{Vi}\circ T-T\circ\nabla_{Ei}$, etc.
\end{proposition}
\proof 
We want a natural morphism
$q_{V,W}:Q(V)\tens_1 Q(W)\to Q(V\tens_0 W)$ but we suppress $Q$ since it is the identity on objects, just viewed with a different $\bullet$ action.  For the proposed $q$ to be well-defined we need
\begin{eqnarray*}
q_{V,W}(v\bullet a\tens_1 w) &=& q_{V,W}(v\tens_1 a\bullet w)\ ,
\end{eqnarray*}
so from (\ref{nbdhksav}),
\begin{eqnarray*}
q_{V,W}((a\,v-\tfrac{\lambda}{2} \,\omega^{ij}\, a_{,i}\, (\nabla_{Vj}v))\tens_1 w) &=& q_{V,W}(v\tens_1 (a\,w+\tfrac{\lambda}{2} \,\omega^{ij}\, a_{,i}\, (\nabla_{Wj}w)))\ ,
\end{eqnarray*}
which is satisfied by the formula for $q_{V,W}$. 

Next, we require each $q_{V,W}$ to be a bimodule map over $A_1$. Thus,
\begin{eqnarray}
q_{V,W}(v\tens_1 (w\bullet a)) &=& v\tens_0 w\,a
- v\tens_0 \tfrac{\lambda}{2} \,\omega^{ij}\, a_{,i}\, \nabla_{Wj}w +\, \tfrac{\lambda}{2} \,\omega^{ij}\,\nabla_{Vi} v
\tens_0 \nabla_{Wj} (wa) \cr
&=& v\tens_0 w\,a
 + \tfrac{\lambda}{2} \,\omega^{ij}\,(\nabla_{Vi} v
\tens_0 \nabla_{Wj} w)\,a+\, \tfrac{\lambda}{2} \,\omega^{ij}\,(\nabla_{V\tens_0 Wi} (v
\tens_0 w))\,a_{,j}\cr
q_{V,W}(v\tens_1 w)\bullet a&=&(v\tens_0 w+\tfrac{\lambda}{2} \omega^{ij}\nabla_{Vi}v\tens_0\nabla_{Wj}w)\bullet a\cr
&=&(v\tens_0 w)\bullet a +\tfrac{\lambda}{2} \omega^{ij}(\nabla_{Vi}v\tens_0\nabla_{Wj}w)a\cr
&=&v\tens_0 w a-\tfrac{\lambda}{2} \omega^{ij}a_{,i}\nabla_{V\tens_0 Wj}(v\tens_0w)+\tfrac{\lambda}{2} \omega^{ij}(\nabla_{Vi}v\tens_0\nabla_{Wj}w)a\ .
\end{eqnarray}
using the quantum right module structure on $V$  etc (i.e. regarding it as $Q(V)$) from (\ref{nbdhksav}). We do not need to use the $\bullet$ product or non-trivial terms in $q$ if an expression already has a $\lambda$, as we are working to errors in $O(\lambda^2)$.  Our two expressions agree using antisymmetry of $\omega$. Similarly on the other side,
\begin{eqnarray}
q_{V,W}((a\bullet v)\tens_1 w) &=& a\,v\tens_0 w
+\tfrac{\lambda}{2} \,\omega^{ij}\, a_{,i}\, (\nabla_{Vj}v)\tens_0 w  +\, \tfrac{\lambda}{2} \,\omega^{ij}\,\nabla_{Vi} (a\,v)
\tens_0 \nabla_{Wj} w \cr
&=& v\tens_0 w\,a + \tfrac{\lambda}{2} \,a\,\omega^{ij}\,\nabla_{Vi} v
\tens_0 \nabla_{Wj} w+\,  \tfrac{\lambda}{2} \,\omega^{ij}\,a_{,i}\,\nabla_{V\tens_0 Wj} (v
\tens_0 w)\cr
a\bullet q_{V,W}(v\tens_1w)&=&a\bullet (v\tens_0 w+\tfrac{\lambda}{2} \omega^{ij}\nabla_{Vi}v\tens_0\nabla_{Wj}w)\cr
&=&a v\tens_0 w+\tfrac{\lambda}{2} \omega^{ij}a_{,i}\nabla_{V\tens_0Wj}(v\tens_0w)+\tfrac{\lambda}{2} \omega^{ij}\nabla_{Vi}v\tens_0\nabla_{Wj}w\ .
\end{eqnarray}

Next, we check that $q_{V,W}$ is functorial. Let $T:V\to Z$ be a morphism in $\CD_0$ (so intertwining the covariant derivatives) and recall that $Q(T)$ is just $T$. Then
\begin{eqnarray}
q_{Z,W}(Tv\tens_1 w)&=&=Tv\tens_0 w+\tfrac{\lambda}{2} \omega^{ij}\nabla_{Zi}(Tv)\tens_0\nabla_{Wj}w\cr
&=&Tv\tens_0w+\tfrac{\lambda}{2} \omega^{ij}(T\circ\nabla_{Vi}v)\tens_0\nabla_{Wj}w=(T\tens\id)q_{V,W}(v\tens_1 w)\ ,\end{eqnarray}
and similarly for functoriality on the other side. 

Finally, it remains to check that $q_{V\tens_0 W,Z}\circ(q_{V,W}\tens\id)=q_{V,W\tens_0Z}\circ(\id\tens q_{W,Z})$ where the associators
implicit here are all trivial to order $\lambda$. This is immediate from the formulae for $q$ working to order $\lambda$. Our $q$ are clearly
also invertible to this order by the same formula with $-\lambda$.  \endproof

 Now we discuss conjugate modules and star operations. For vector bundles with connection on real manifolds, we define covariant derivatives of conjugates in the obvious manner, $\nabla_{\bar E i}(\overline{e})=\overline{\nabla_{Ei} e}$. A star operation on a vector bundle will be conjugate linear bundle map to itself $e\mapsto e^*$, and will be compatible with a connection if $\nabla_{\bar Ei}(e^*)=(\nabla_{Ei}e)^*$. It will be convenient to take the linear map to the conjugate bundle $\star:E\to \overline{E}$ defined by $\star(e)=\overline{e^*}$.

\begin{proposition} \label{kihjvfcyiouy}
Over $\C$, the functor $Q:\CD_0\to\CE_1$ is a bar functor. Hence, if $\star:E\to \overline{E}$ is a star object and compatible with the connection, then $Q(\star):Q(E)\to\overline{Q(E)}$ is also a star object. 
\end{proposition}
\proof  To show we have a functor, we begin by identifying $Q(\overline{E})$ and $\overline{Q(E)}$ (Remember that $Q(E)$ is simply $E$ but with a different module structure, so $Q(e)\in Q(E)$ is simply $e\in E$ as sets.)
 To do this we need to show that, for all $a\in A_1$ and $e\in E$,
\begin{eqnarray*}
a.\overline{Q(e)} &=& \overline{Q(e).a^*} =\overline{Q(e\bullet a^*)} \ ,\quad
a.Q(\overline{e}) = Q(a\bullet \overline{e})\ ,
\end{eqnarray*}
so we need to show that $a\bullet \overline{e}=\overline{e\bullet a^*}$. 
Now
\begin{eqnarray*}
a\bullet \overline{e} &=& a. \overline{e}+\tfrac{\lambda}{2}\,\omega^{ij}\,a_{,i}.\nabla_j(\overline{e})  
= a. \overline{e}+\tfrac{\lambda}{2}\,\omega^{ij}\,a_{,i}.\overline{\nabla_j e} \cr
&=& \overline{e.a^*}+\tfrac{\lambda}{2}\,\omega^{ij}\,\overline{\nabla_j e.a_{,i}{}^*}   
= \overline{e.a^*}-\overline{\tfrac{\lambda}{2}\,\omega^{ij}\,\nabla_j e.a_{,i}{}^*}  \cr
&=& \overline{e.a^*}+\overline{\tfrac{\lambda}{2}\,\omega^{ij}\,\nabla_i e.a_{,j}{}^*} 
= \overline{e\bullet a^*}\ .
\end{eqnarray*}
Now we have to check the morphisms $T:E\to F$,
that $Q(\overline{T})=\overline{Q(T)}$. 
\begin{eqnarray*}
\overline{Q(T)}(\overline{Q(e)}) &=& \overline{Q(T)(Q(e))} = \overline{Q(T(e)+\tfrac\lambda2\,\omega^{ij}\,\nabla_i(\nabla_j(T)e))}\ ,\cr
Q(\overline{T})(Q(\overline{e}) &=& \overline{T}(\overline{e}) + \tfrac\lambda2\,\omega^{ij}\,\nabla_i(
\nabla_j(\overline{T})(\overline{e}))
\end{eqnarray*}
Now we check what $\nabla_j(\overline{T})$ is:
\begin{eqnarray*}
\nabla_j(\overline{T})(\overline{e}) &=& \nabla_j(\overline{T(e)})-\overline{T}(\nabla_j(\overline{e})) \cr
&=& \overline{\nabla_j(T(e))-T(\nabla_j(e))}=\overline{\nabla_j(T)e}\ .
\end{eqnarray*}
Then we have, as $\lambda$ is imaginary,
\begin{eqnarray*}
Q(\overline{T})(Q(\overline{e}) &=& \overline{T(e)} + \tfrac\lambda2\,\omega^{ij}\,\nabla_i(
\overline{\nabla_j(T)e})  \cr
&=&  \overline{T(e)} + \tfrac\lambda2\,\omega^{ij}\,
\overline{\nabla_i(\nabla_j(T)e)}  \cr
&=&  \overline{T(e)- \tfrac\lambda2\,\omega^{ij}\,
\nabla_i(\nabla_j(T)e)} \ .
\end{eqnarray*}
Thus we have
\begin{eqnarray*}
(\overline{Q(T)}- Q(\overline{T})   )(\overline{Q(e)}) &=&\overline{Q(\lambda\,\omega^{ij}\,\nabla_i(\nabla_j(T)e))}\ ,
\end{eqnarray*}
so if $T$ is a morphism in $\CD_0$ we get $\overline{Q(T)} = Q(\overline{T}) $. 

Now we show that the natural transformation $q$ is compatible with the natural transformation $\Upsilon$ in the bar category. This means that the following diagram commutes:
\begin{eqnarray}\label{vcusvdsavbbdfsii}
\xymatrix{
\overline{Q(E)\tens_1 Q(F)} \ar[r]^{\overline{q_{E,F}}} \ar[d]_{\Upsilon_{Q(E),Q(F)}} 
& \overline{Q(E\tens_0 F)}  \ar[r]^{ = }  &   Q(\overline{E\tens_0 F}) \ar[d]_{\Upsilon_{E,F}}   \\
\overline{ Q(F)}\tens_1 \overline{Q(E)} \ar[r]^{=} &
Q(\overline{ F}) \tens_1 Q(\overline{E} )  \ar[r]^{q_{\overline{ F},\overline{ E}}}    & Q(\overline{F} \tens_0 \overline{E} )
}
\end{eqnarray}
Now we have
\begin{eqnarray*}
\Upsilon_{E,F}(\overline{q_{E,F}}(\overline{Q(e)\tens_1 Q(f)})) &=& 
\Upsilon_{E,F}(\overline{Q(e\tens_0 f+\tfrac\lambda2\,\omega^{ij}\,\nabla_i e\tens_0\nabla_j f)}) \cr
&=& \Upsilon_{E,F}(Q(\overline{e\tens_0 f+\tfrac\lambda2\,\omega^{ij}\,\nabla_i e\tens_0\nabla_j f})) \cr
&=& Q(\overline{f}\tens_0 \overline{e} -\tfrac\lambda2\,\omega^{ij}\,\overline{ \nabla_j f  } \tens_0 \overline{  \nabla_i e }) \cr
&=& Q(\overline{f}\tens_0 \overline{e} +\tfrac\lambda2\,\omega^{ij}\,\nabla_i\overline{ f  } \tens_0 \nabla_j  \overline{ e }) \ , \cr
q_{\overline{ F},\overline{ E}}   \Upsilon_{Q(E),Q(F)}    (\overline{Q(e)\tens_1 Q(f)}) &=& 
q_{\overline{ F},\overline{ E}}  (\overline{Q(f)}\tens_1 \overline{Q(e)})  \cr
&=&  q_{\overline{ F},\overline{ E}}  (Q(\overline{f})\tens_1 Q(\overline{e}))  \cr
&=& Q(\overline{f}\tens_0 \overline{e} +\tfrac\lambda2\,\omega^{ij}\,\nabla_i\overline{ f  } \tens_0 \nabla_j  \overline{ e }) \ .\quad\square
\end{eqnarray*}

\subsection{Quantising the quantising covariant derivative} \label{cnioipubohv}

We now want to extend the functor $Q$ above to a functor $Q:\CD_0\to \CD_1$. 

\begin{theorem} \label{functor} Let $(E,\nabla_E)$ be a classical bundle and connection. Then $E$ with the bimodule structure $\bullet$ over $A_1$ has bimodule covariant derivative
\begin{eqnarray*}
\nabla_{Q(E)} &=& q^{-1}_{\Omega^1,E}\nabla_E-\tfrac{\lambda}{2} \,\omega^{ij}\,\extd x^k\tens_1[\nabla_{Ek},\nabla_{Ej}]\nabla_{Ei} \cr
\quad \sigma_{Q(E)}(e\tens_1\xi)\, &=&  \xi\tens_1 e  + \lambda\, \omega^{ij}\,\nabla_j\xi\tens_1\nabla_{Ei}e   + \, \lambda\, \omega^{ij}\,\xi_{j}\,\extd x^k \tens_1 [\nabla_{Ek},\nabla_{Ei}]e  \ ,
\end{eqnarray*}
where we view contraction with $\omega$ as a map $\Omega^1(M) \tens_0\Omega^2(M)\to \Omega^1(M)$.  Moreover, $Q(E,\nabla_E)=(Q(E),\nabla_{Q(E)})$ is a monoidal functor $Q:\CD_0\to \CD_1$ via $q$ as in Proposition~\ref{q}.
\end{theorem}
\proof
We start by considering
\begin{eqnarray*}
q^{-1}\nabla_E(a\,\la_1\,e) &=& q^{-1}\nabla_E\big(
a\,e+\lambda\,\omega^{ij}\,a_{,i}\,\nabla_je/2\big)  \cr
&=&  q^{-1}\big(\extd a\tens_0 e + a\,\extd x^k\tens_0\nabla_k e
 +\lambda\,\extd(a_{,i}\,\omega^{ij})\tens_0 \nabla_je/2 \cr
&& +\ \lambda\,\omega^{ij}\,a_{,i}\,\extd x^k\tens_0\nabla_k\nabla_je/2\big) \cr
&=&  \extd a\tens_1 e + a\,\extd x^k\tens_1\nabla_k e
 +\lambda\,\extd(a_{,i}\,\omega^{ij})\tens_1 \nabla_je/2 \cr
&& +\ \lambda\,\omega^{ij}\,a_{,i}\,\extd x^k\tens_1\nabla_k\nabla_je/2 \cr
&& -\, \lambda\,\omega^{ij}\, \nabla_i\extd a\tens_1 \nabla_j e/2
- \lambda\,\omega^{ij}\, \nabla_i(a\,\extd x^k)\tens_1\nabla_j\nabla_k e/2 \cr
&=&  \extd a\tens_1 e + a\,\extd x^k\tens_1\nabla_k e
 +\lambda\,\extd(a_{,i}\,\omega^{ij})\tens_1 \nabla_je/2 \cr
&& -\, \lambda\,\omega^{ij}\, \nabla_i\extd a\tens_1 \nabla_j e/2
- \lambda\,\omega^{ij}\,a\, \nabla_i(\extd x^k)\tens_1\nabla_j\nabla_k e/2\cr
&& +\ \lambda\,\omega^{ij}\,a_{,i}\,\extd x^k\tens_1[\nabla_k,\nabla_j]e/2  \ ,
\end{eqnarray*}
and
\begin{eqnarray*}
a\,\la_1\,q^{-1}\nabla_E(e) &=& a\,\la_1\,q^{-1}\big(
\extd x^k\tens_0\nabla_ke\big)  \cr
&=&  a\,\la_1\,\big(\extd x^k\tens_1\nabla_ke-\lambda\,\omega^{ij}\,\nabla_i(\extd x^k)\tens_1
\nabla_j\nabla_ke/2\big) \cr
&=& a\,\extd x^k\tens_1\nabla_ke-\lambda\,\omega^{ij}\,a\,\nabla_i(\extd x^k)\tens_1
\nabla_j\nabla_ke/2  \cr
&&+\, \lambda\,\omega^{ij}\,a_{,i}\, \nabla_j(\extd x^k)\tens_1\nabla_ke/2\ ,
\end{eqnarray*}
and then
\begin{eqnarray}  \label{chdksauyigc}
&& q^{-1}\nabla_E(a\,\la_1\,e) -a\,\la_1\,q^{-1}\nabla_E(e)  \cr
&=&  \extd a\tens_1 e 
 +\lambda\,\extd(a_{,i}\,\omega^{ij})\tens_1 \nabla_je/2 +\lambda\,\omega^{ij}\,a_{,i}\,\extd x^k\tens_1[\nabla_k,\nabla_j]e/2 \cr
&& -\, \lambda\,\omega^{ij}\, \nabla_i\extd a\tens_1 \nabla_j e/2
-\lambda\,\omega^{ij}\,a_{,i}\, \nabla_j(\extd x^k)\tens_1\nabla_ke/2 \cr
&=&  \extd a\tens_1 e 
 +\lambda\,\extd(a_{,i}\,\omega^{ij})\tens_1 \nabla_je/2 +\lambda\,\omega^{ij}\,a_{,i}\,\extd x^k\tens_1[\nabla_k,\nabla_j]e)/2 \cr
&& -\, \lambda\,\omega^{ij}\, \nabla_i(a_{,k}\,\extd x^k)\tens_1 \nabla_j e/2
-\lambda\,\omega^{ik}\,a_{,i}\, \nabla_k(\extd x^j)\tens_1\nabla_je/2 \cr
&=&  \extd a\tens_1 e 
 +\lambda\,a_{,i}\,\extd(\omega^{ij})\tens_1 \nabla_je/2 + \lambda\,\omega^{ij}\,a_{,i}\,\extd x^k\tens_1[\nabla_k,\nabla_j]e/2 \cr
&& -\, \lambda\,\omega^{ij}\, a_{,k}\,\nabla_i(\extd x^k)\tens_1 \nabla_j e/2
-\lambda\,\omega^{ik}\,a_{,i}\, \nabla_k(\extd x^j)\tens_1\nabla_je/2 \cr
&=&  \extd a\tens_1 e 
 +\lambda\,a_{,i}\,\extd(\omega^{ij})\tens_1 \nabla_je/2  +\lambda\,\omega^{ij}\,a_{,i}\,\extd x^k\tens_1[\nabla_k,\nabla_j]e/2\cr
&& -\, \lambda\,\omega^{kj}\, a_{,i}\,\nabla_k(\extd x^i)\tens_1 \nabla_j e/2
-\lambda\,\omega^{ik}\,a_{,i}\, \nabla_k(\extd x^j)\tens_1\nabla_je/2 \cr
&=&  \extd a\tens_1 e 
 +\lambda\,a_{,i}\,\big(\extd(\omega^{ij})  - \omega^{kj}\,\nabla_k(\extd x^i)
-\omega^{ik}\, \nabla_k(\extd x^j)   \big)\tens_1\nabla_je/2 \cr
&&+\, \lambda\,\omega^{ij}\,a_{,i}\,\extd x^k\tens_1[\nabla_k,\nabla_j]e/2 \ .
\end{eqnarray}
If we have the condition (\ref{poissoncomp})
then the last long bracket in (\ref{chdksauyigc}) vanishes, giving
\begin{eqnarray}  \label{chdksauyigc}
\quad q^{-1}\nabla_E(a\,\la_1\,e) -a\,\la_1\,q^{-1}\nabla_E(e)  
&=&  \extd a\tens_1 e +\lambda\,\omega^{ij}\,a_{,i}\,\extd x^k\tens_1[\nabla_{Ek},\nabla_{Ej}]e/2 \ .
\end{eqnarray}
Now we can set the first order quantisation of the 
left covariant derivative to be 
\begin{eqnarray*}
Q(\nabla_E)(e)\,=\,
q^{-1}_{\Omega^1,E}\nabla_E(e)-\tfrac{\lambda}{2} \,\omega^{ij}\,\extd x^k\tens_1[\nabla_{Ek},\nabla_{Ej}]\nabla_{Ei}(e)
\end{eqnarray*}
which we can write as stated.

Next, we want to see about a bimodule connection. We compute
\begin{eqnarray*}
\sigma(e\tens_1\extd a) &=& \extd a\tens_1 e + \nabla[e,a]+[a,\nabla e] \cr
&=&  \extd a\tens_1 e + \lambda\, \nabla(\omega^{ij}\,\nabla_i(e)\,a_{,j})+[a,\extd x^k \tens_1 \nabla_k e] \cr
&=&  \extd a\tens_1 e + \lambda\, \extd(\omega^{ij}\,a_{,j})\tens\nabla_i(e)
  + \lambda\, \omega^{ij}\,a_{,j}\,\extd x^k \tens_1 \nabla_k\nabla_i(e)  \cr
&& +\, \lambda\,\omega^{ij}\,a_{,i}\,\nabla_j(\extd x^k) \tens_1 \nabla_k e
+\lambda\,\omega^{ij}\,a_{,i}\,\extd x^k \tens_1 \nabla_j \nabla_k e \cr
&=&  \extd a\tens_1 e + \lambda\, \extd(\omega^{ij})\,a_{,j}\tens\nabla_i(e) + \lambda\, \omega^{ij}\,a_{,jk}\,\extd x^k\tens_1\nabla_i(e) \cr
&&  + \, \lambda\, \omega^{ij}\,a_{,j}\,\extd x^k \tens_1 \nabla_k\nabla_i(e)  \cr
&& +\, \lambda\,\omega^{ij}\,a_{,i}\,\nabla_j(\extd x^k) \tens_1 \nabla_k e
+\lambda\,\omega^{ji}\,a_{,j}\,\extd x^k \tens_1 \nabla_i \nabla_k e \cr
&=&  \extd a\tens_1 e + \lambda\, \extd(\omega^{ij})\,a_{,j}\tens\nabla_i(e) + \lambda\, \omega^{ij}\,a_{,jk}\,\extd x^k\tens_1\nabla_i(e) \cr
&&  + \, \lambda\, \omega^{ij}\,a_{,j}\,\extd x^k \tens_1 [\nabla_k,\nabla_i](e) +
 \lambda\,\omega^{jk}\,a_{,j}\,\nabla_k(\extd x^i) \tens_1 \nabla_i e\cr
 &=&  \extd a\tens_1 e  + \lambda\, \omega^{ij}\,\nabla_j(a_{,k}\,\extd x^k)\tens_1\nabla_i(e)
   - \lambda\, \omega^{ij}\,a_{,k}\,\nabla_j(\extd x^k)\tens_1\nabla_i(e) \cr
&&  + \, \lambda\, \omega^{ij}\,a_{,j}\,\extd x^k \tens_1 [\nabla_k,\nabla_i](e)  \cr
&& +\,  \lambda\,a_{,j}\big( \extd(\omega^{ij})-
\omega^{kj}\,\nabla_k(\extd x^i)\big) \tens_1 \nabla_i e\cr
 &=&  \extd a\tens_1 e  + \lambda\, \omega^{ij}\,\nabla_j(a_{,k}\,\extd x^k)\tens_1\nabla_i(e) \cr
&&  + \, \lambda\, \omega^{ij}\,a_{,j}\,\extd x^k \tens_1 [\nabla_k,\nabla_i](e)  \cr
&& +\,  \lambda\,a_{,j}\big( \extd(\omega^{ij})-
\omega^{kj}\,\nabla_k(\extd x^i)-\omega^{ik}\,\nabla_k(\extd x^j)\big) \tens_1 \nabla_i e\ ,
\end{eqnarray*}
and under condition (\ref{poissoncomp}) this becomes the given formula. This constructs the quantised covariant derivative. The following two lemmas then verify the desired categorical properties so as to complete the proof. \endproof

The following two lemmas complete the proof of the theorem. The first one by showing that the covariant derivatiive is natural with respect to the tensor product structure -- i.e.\ that the quantisation of the classical tensor product covariant derivative is the tensor product of the quantised covariant derivatives (using the $\sigma$ map). This is summarised by
\begin{eqnarray}\label{vcusvdsavbbdfs}
\xymatrix{
Q(E)\tens_1 Q(F) \ar[r]^{q} \ar[d]_{\nabla_{Q(E)\tens_1 Q(F)}} 
&Q(E\tens_0 F)\ar[d]_{\nabla_{Q(E\tens_0 F)}}      \\
\Omega^1A_1\tens_1 Q(E)\tens_1 Q(F) \ar[r]_{\id\tens q} &
\Omega^1A_1\tens_1  Q(E\tens_0 F)    }
\end{eqnarray}

\begin{lemma}\label{qfunct} For all $e\in E$ and $f\in F$,
\[ (\id\tens q_{E,F})\left(\nabla_{Q(E)}e\tens_1 f+(\sigma_{Q(E)}\tens\id)(e\tens_1\nabla_{Q(F)}f)\right)=\nabla_{Q(E\tens_0 F)}q_{E,F}(e\tens_1 f)\]
\end{lemma}
\proof 
Begin with
\begin{eqnarray*}
&& \nabla_{Q(E\tens_0 F)} q_{E,F}(e\tens_1 f) \cr
 &=& \nabla_{Q(E\tens_0 F)} (e\tens_0 f) +\tfrac{\lambda}{2}\,
 \nabla_{Q(E\tens_0 F)} \big(\omega^{ij}\, \nabla_i e\tens_0\nabla_j f\big) \cr
&=& q^{-1}_{\Omega^1,E}\nabla(e\tens f)-\tfrac{\lambda}{2} \,\omega^{ij}\,\extd x^k\tens_1[\nabla_{k},\nabla_{j}]\nabla_{i}(e\tens f)  + \tfrac{\lambda}{2}\,
 \nabla_{Q(E\tens_0 F)} \big(\omega^{ij}\, \nabla_i e\tens_0\nabla_j f\big) \cr
&=&q^{-1}_{\Omega^1,E\tens F}(\extd x^k\tens (\nabla_k e\tens f) + \extd x^k\tens (e\tens \nabla_k f)) +
\tfrac{\lambda}{2}\,\extd (\omega^{ij})\tens_1 \big(\nabla_i e\tens_0\nabla_j f\big) \cr
&& -\, \tfrac{\lambda}{2} \,\omega^{ij}\,\extd x^k\tens_1[\nabla_{k},\nabla_{j}](\nabla_{i}e\tens f+e\tens \nabla_{i}f) 
+ \tfrac{\lambda}{2}\,\omega^{ij}\,\extd x^k\tens_1
 \nabla_k \big( \nabla_i e\tens_0\nabla_j f\big)\cr
&=&\extd x^k\tens_1 (\nabla_k e\tens f) + \extd x^k\tens_1 (e\tens \nabla_k f)  \cr
&& -\, \tfrac{\lambda}{2} \,\omega^{ij}\,\nabla_i(\extd x^k)\tens_1 \nabla_j\big(
(\nabla_k e\tens f) + (e\tens \nabla_k f)\big)+ \tfrac{\lambda}{2}\,\omega^{ij}\,\extd x^k\tens_1
 \nabla_k \big( \nabla_i e\tens_0\nabla_j f\big)\cr
&& -\, \tfrac{\lambda}{2} \,\omega^{ij}\,\extd x^k\tens_1[\nabla_{k},\nabla_{j}](\nabla_{i}e\tens f+e\tens \nabla_{i}f) 
 + \tfrac{\lambda}{2}\,\extd (\omega^{ij})\tens_1 \big(\nabla_i e\tens_0\nabla_j f\big)\ .
\end{eqnarray*}
Also
\begin{eqnarray*}
&&(\id\tens q_{E,F})(\nabla_{Q(E)}e\tens_1 f) \cr
&=& (\id\tens q_{E,F})\big(\big( q^{-1}_{\Omega^1,E}\nabla_E(e)-\tfrac{\lambda}{2} \,\omega^{ij}\,\extd x^k\tens_1[\nabla_{k},\nabla_{j}]\nabla_{i}(e)\big) \tens_1 f\big) \cr
&=&(\id\tens q_{E,F})\big(\big( \extd x^k \tens_1\nabla_k(e)-\tfrac{\lambda}{2} \,\omega^{ij}\,\nabla_i(\extd x^k)\tens_1 \nabla_{j}\nabla_{k}(e)\big) \tens_1 f\big) \cr
&&  -  \, \tfrac{\lambda}{2} \,\omega^{ij}\,\extd x^k\tens_1\big([\nabla_{k},\nabla_{j}]\nabla_{i}(e) \tens_0 f\big) \cr
&=& \extd x^k \tens_1(\nabla_k(e)\tens_0 f)-\tfrac{\lambda}{2} \,\omega^{ij}\,\nabla_i(\extd x^k)\tens_1 \big(\nabla_{j}\nabla_{k}(e) \tens_0 f\big) \cr
&&  -  \, \tfrac{\lambda}{2} \,\omega^{ij}\,\extd x^k\tens_1\big([\nabla_{k},\nabla_{j}]\nabla_{i}(e) \tens_0 f\big)  +   \tfrac{\lambda}{2} \, \extd x^k \tens_1\big(  \omega^{ij}(\nabla_i\nabla_k(e)\tens_0 \nabla_jf)\big) \ .
\end{eqnarray*}
Then
\begin{eqnarray*}
&& \nabla_{Q(E\tens_0 F)} q_{E,F}(e\tens_1 f)-(\id\tens q_{E,F})(\nabla_{Q(E)}e\tens_1 f) \cr
&=&   \extd x^k\tens_1 (e\tens \nabla_k f)  - \tfrac{\lambda}{2} \,\omega^{ij}\,\nabla_i(\extd x^k)\tens_1 \big(
\nabla_k e\tens\nabla_j f \big)  \cr
&& -\, \tfrac{\lambda}{2} \,\omega^{ij}\,\nabla_i(\extd x^k)\tens_1 \nabla_j\big(e\tens \nabla_k f\big)+ \tfrac{\lambda}{2}\,\omega^{ij}\,\extd x^k\tens_1
 \nabla_k \big( \nabla_i e\tens_0\nabla_j f\big)\cr
&& -\, \tfrac{\lambda}{2} \,\omega^{ij}\,\extd x^k\tens_1(\nabla_{i}e\tens[\nabla_{k},\nabla_{j}] f) 
 + \tfrac{\lambda}{2}\,\extd (\omega^{ij})\tens_1 \big(\nabla_i e\tens_0\nabla_j f\big)\cr
 && -\, \tfrac{\lambda}{2} \,\omega^{ij}\,\extd x^k\tens_1(\nabla_{k}\nabla_{j}e\tens \nabla_{i}f)
 - \tfrac{\lambda}{2} \,\omega^{ij}\,\extd x^k\tens_1(e\tens [\nabla_{k},\nabla_{j}]\nabla_{i}f)  \cr 
 && -\,  \lambda\, \extd x^k \tens_1\big(  \omega^{ij}(\nabla_i\nabla_k(e)\tens_0 \nabla_jf)\big) \cr
 &=&   \extd x^k\tens_1 (e\tens \nabla_k f)  - \tfrac{\lambda}{2} \,\omega^{ij}\,\nabla_i(\extd x^k)\tens_1 \big(
\nabla_k e\tens\nabla_j f \big)  \cr
&& -\, \tfrac{\lambda}{2} \,\omega^{ij}\,\nabla_i(\extd x^k)\tens_1 \nabla_j\big(e\tens \nabla_k f\big)+ \tfrac{\lambda}{2}\,\omega^{ij}\,\extd x^k\tens_1
 \big( \nabla_k \nabla_i e\tens_0\nabla_j f\big)\cr
&& +\, \tfrac{\lambda}{2} \,\omega^{ij}\,\extd x^k\tens_1(\nabla_{i}e\tens \nabla_{j}\nabla_{k} f) 
 + \tfrac{\lambda}{2}\,\extd (\omega^{ij})\tens_1 \big(\nabla_i e\tens_0\nabla_j f\big)\cr
 && -\, \tfrac{\lambda}{2} \,\omega^{ij}\,\extd x^k\tens_1(\nabla_{k}\nabla_{j}e\tens \nabla_{i}f)
 - \tfrac{\lambda}{2} \,\omega^{ij}\,\extd x^k\tens_1(e\tens [\nabla_{k},\nabla_{j}]\nabla_{i}f)  \cr 
  && -\,  \lambda\, \extd x^k \tens_1\big(  \omega^{ij}(\nabla_i\nabla_k(e)\tens_0 \nabla_jf)\big) \cr
  &=&   \extd x^k\tens_1 (e\tens \nabla_k f)  - \tfrac{\lambda}{2} \,\omega^{ij}\,\nabla_i(\extd x^k)\tens_1 \big(
\nabla_k e\tens\nabla_j f \big)  \cr
&& -\, \tfrac{\lambda}{2} \,\omega^{ij}\,\nabla_i(\extd x^k)\tens_1 \big(\nabla_j e\tens \nabla_k f\big)+ \tfrac{\lambda}{2}\,\omega^{ij}\,\extd x^k\tens_1
 \big( \nabla_k \nabla_i e\tens_0\nabla_j f\big)\cr
&& +\, \tfrac{\lambda}{2} \,\omega^{ij}\,\extd x^k\tens_1(\nabla_{i}e\tens \nabla_{j}\nabla_{k} f) 
 + \tfrac{\lambda}{2}\,\extd (\omega^{ij})\tens_1 \big(\nabla_i e\tens_0\nabla_j f\big)\cr
 && -\, \tfrac{\lambda}{2} \,\omega^{ij}\,\extd x^k\tens_1(\nabla_{k}\nabla_{j}e\tens \nabla_{i}f)
 - \tfrac{\lambda}{2} \,\omega^{ij}\,\extd x^k\tens_1(e\tens [\nabla_{k},\nabla_{j}]\nabla_{i}f)  \cr 
 && -\, \tfrac{\lambda}{2} \,\omega^{ij}\,\nabla_i(\extd x^k)\tens_1 \big(e\tens\nabla_j \nabla_k f\big)   -  \lambda\, \extd x^k \tens_1\big(  \omega^{ij}(\nabla_i\nabla_k(e)\tens_0 \nabla_jf)\big) \cr
   &=&   \extd x^k\tens_1 (e\tens \nabla_k f)  - \tfrac{\lambda}{2} \,\omega^{ij}\,\nabla_i(\extd x^k)\tens_1 \big(
\nabla_k e\tens\nabla_j f \big)  \cr
&& -\, \tfrac{\lambda}{2} \,\omega^{ij}\,\nabla_i(\extd x^k)\tens_1 \big(\nabla_j e\tens \nabla_k f\big)+ \lambda\,\omega^{ij}\,\extd x^k\tens_1
 \big( [\nabla_k, \nabla_i] e\tens_0\nabla_j f\big)\cr
&& +\, \tfrac{\lambda}{2} \,\omega^{ij}\,\extd x^k\tens_1(\nabla_{i}e\tens \nabla_{j}\nabla_{k} f) 
 + \tfrac{\lambda}{2}\,\extd (\omega^{ij})\tens_1 \big(\nabla_i e\tens_0\nabla_j f\big)\cr
 && -\,  \tfrac{\lambda}{2} \,\omega^{ij}\,\extd x^k\tens_1(e\tens [\nabla_{k},\nabla_{j}]\nabla_{i}f)   - \tfrac{\lambda}{2} \,\omega^{ij}\,\nabla_i(\extd x^k)\tens_1 \big(e\tens\nabla_j \nabla_k f\big)  \ .
\end{eqnarray*}
Next
\begin{eqnarray*}
&&  (\sigma_{Q(E)}\tens\id)(e\tens_1\nabla_{Q(F)}f)   \cr
&=& (\sigma_{Q(E)}\tens\id)\big(e\tens_1q^{-1}_{\Omega^1,E}\nabla(f)-\tfrac{\lambda}{2} \,\omega^{ij}\,e\tens_1(\extd x^k\tens_1[\nabla_{k},\nabla_{j}]\nabla_{i}(f))\big)   \cr
&=& (\sigma_{Q(E)}\tens\id)\big(e\tens_1\extd x^k\tens_1 \nabla_k(f)-\tfrac{\lambda}{2} \,\omega^{ij}\,
e\tens_1\nabla_i( \extd x^k)\tens_1 \nabla_j\nabla_k(f) \cr
&& - \, \tfrac{\lambda}{2} \,\omega^{ij}\,e\tens_1\extd x^k\tens_1[\nabla_{k},\nabla_{j}]\nabla_{i}(f)\big)   \cr
&=& \extd x^k\tens_1 e\tens_1 \nabla_k(f)-\tfrac{\lambda}{2} \,\omega^{ij}\,
\nabla_i( \extd x^k)\tens_1 e\tens_1 \nabla_j\nabla_k(f) \cr
&& - \, \tfrac{\lambda}{2} \,\omega^{ij}\,\extd x^k\tens_1 e\tens_1[\nabla_{k},\nabla_{j}]\nabla_{i}(f) + \lambda\,\omega^{ij}\,\nabla_j(\extd x^k)\tens_1 \nabla_i e \tens_1 \nabla_k f \cr
&&+\, \lambda\,\omega^{ik}\,\extd x^l \tens_1 [\nabla_l,\nabla_i]e \tens_1 \nabla_k f\ ,
\end{eqnarray*}
and so
\begin{eqnarray*}
&&  (\id\tens q_{E,F})(\sigma_{Q(E)}\tens\id)(e\tens_1\nabla_{Q(F)}f)   \cr
&=& \extd x^k\tens_1 (e\tens \nabla_k(f))-\tfrac{\lambda}{2} \,\omega^{ij}\,
\nabla_i( \extd x^k)\tens_1 (e\tens \nabla_j\nabla_k(f)) \cr
&& - \, \tfrac{\lambda}{2} \,\omega^{ij}\,\extd x^k\tens_1 (e\tens[\nabla_{k},\nabla_{j}]\nabla_{i}(f)) + \lambda\,\omega^{ij}\,\nabla_j(\extd x^k)\tens_1 ( \nabla_i e \tens \nabla_k f) \cr
&&+\, \lambda\,\omega^{ik}\,\extd x^l \tens_1 ([\nabla_l,\nabla_i]e \tens \nabla_k f)
+  \tfrac{\lambda}{2} \,\omega^{ij}\,\extd x^k\tens_1 (\nabla_i e\tens \nabla_j\nabla_k(f))
\end{eqnarray*}
From this we have
\begin{eqnarray*}
&& \nabla_{Q(E\tens_0 F)} q_{E,F}(e\tens_1 f)-(\id\tens q_{E,F})(\nabla_{Q(E)}e\tens_1 f) -(\id\tens q_{E,F})(\sigma_{Q(E)}\tens\id)(e\tens_1\nabla_{Q(F)}f)  \cr
   &=&     - \tfrac{\lambda}{2} \,\omega^{ij}\,\nabla_i(\extd x^k)\tens_1 \big(
\nabla_k e\tens\nabla_j f \big)   -  \tfrac{\lambda}{2} \,\omega^{ij}\,\nabla_i(\extd x^k)\tens_1 \big(\nabla_j e\tens \nabla_k f\big)    \cr
&& -\,  \lambda\,\omega^{ij}\,\nabla_j(\extd x^k)\tens_1 ( \nabla_i e \tens \nabla_k f) 
 + \tfrac{\lambda}{2}\,\extd (\omega^{ij})\tens_1 \big(\nabla_i e\tens_0\nabla_j f\big)\cr
 &=& \tfrac{\lambda}{2}\,\big(\extd (\omega^{ij})-\omega^{kj}\,\nabla_k(\extd x^i)+\omega^{k   i}\,\nabla_k(\extd x^{j})
 \big)\tens_1 \big(\nabla_i e\tens_0\nabla_j f\big)
\end{eqnarray*}
and this vanishes by (\ref{poissoncomp}).\endproof

\medskip

 The second lemma checks functoriality under morphisms.

\begin{lemma}\label{functorT}
If $T:E\to F$ is a module map and commutes with the covariant derivative, then $Q(\nabla_F)Q(T)=(\id\tens_1 Q(T))\,Q(\nabla_E)$. 
\end{lemma}
\proof In this case $Q(T)=T$, and 
\begin{eqnarray*}
(\id\tens_1 Q(T))\,\nabla_{Q(E)} &=& (\id\tens_1 Q(T))\,q^{-1}_{\Omega^1,E}\nabla_E(e)-\tfrac{\lambda}{2} \,\omega^{ij}\,\extd x^k\tens_1 T [\nabla_{Ek},\nabla_{Ej}]\nabla_{Ei}(e) \cr
 &=& \,q^{-1}_{\Omega^1,E}(\id\tens T)\,\nabla_E(e)-\tfrac{\lambda}{2} \,\omega^{ij}\,\extd x^k\tens_1 T [\nabla_{Ek},\nabla_{Ej}]\nabla_{Ei}(e) \ ,\cr
 Q(\nabla_F)T &=& q^{-1}_{\Omega^1,E}\nabla_E T(e)-\tfrac{\lambda}{2} \,\omega^{ij}\,\extd x^k\tens_1[\nabla_{Ek},\nabla_{Ej}]\nabla_{Ei}T(e)\ .
\end{eqnarray*}
These are equal, as $T$ commutes with the covariant derivative. \endproof

\begin{lemma} \label{cluiktycxtu}
Over $\C$, $\nabla_Q$ is star preserving to order $\lambda$, i.e.\ the following diagram commutes:
\begin{eqnarray}\label{vcusvdfsii}
\xymatrix{
\overline{Q(E)} = Q(\overline{E})  \ar[d]_{\overline{\nabla_Q } } 
& Q(E)  \ar[l]_{  \star  }    \ar[r]^{ \nabla_{Q} }  &   
Q(\Omega^1(M)) \tens_1 Q(E)   \ar[d]_{\star\tens_1\star}   \\
\overline{ Q(\Omega^1(M)) \tens_1 Q(E)} \ar[r]^{  \overline{\sigma_{QE}^{-1} }  } &
\overline{ Q(E) \tens_1 Q(\Omega^1(M))}   \ar[r]^{\Upsilon}    &
\overline{ Q(\Omega^1(M)) } \tens_1 \overline{ Q(E) }
}
\end{eqnarray}
 \end{lemma}\proof
Begin with, using $q$ a natural transformation and (\ref{vcusvdsavbbdfsii})
\begin{eqnarray*}
&&\overline{\sigma_{QE} } \,  \Upsilon^{-1}(\star\tens_1\star)\nabla_Q(Q(e))  \cr &=& 
\overline{\sigma_{QE} } \Upsilon^{-1}(\star\tens_1\star)(q^{-1}(\extd x^p\tens_0\nabla_p e)
-\tfrac{\lambda}{2} \,\omega^{ij}\,\extd x^k\tens_1[\nabla_{k},\nabla_{j}]\nabla_{i}e)  \cr
 &=& 
\overline{\sigma_{QE} } \Upsilon^{-1}(q^{-1}(\overline{ \extd x^p}\tens_0\overline{\nabla_p e^*})
-\tfrac{\lambda}{2} \,\omega^{ij}\,\overline{\extd x^k}\tens_1\overline{[\nabla_{k},\nabla_{j}]\nabla_{i}e^*})  \cr
&=& \overline{\sigma_{QE}\,q^{-1} }   \Upsilon^{-1}(\overline{ \extd x^p}\tens_0\overline{\nabla_p e^*})
-\tfrac{\lambda}{2} \,\omega^{ij}\,\overline{\sigma_{QE} }(
\overline{[\nabla_{k},\nabla_{j}]\nabla_{i}e^*\tens_1\extd x^k }   ) \cr
&=& \overline{\sigma_{QE}\,q^{-1} }  (\overline{\nabla_p e^*\tens_0 \extd x^p})
-\tfrac{\lambda}{2} \,\omega^{ij}\,\overline{\sigma_{QE} }(
\overline{[\nabla_{k},\nabla_{j}]\nabla_{i}e^*\tens_1\extd x^k }   ) \cr
&=& \overline{q^{-1}q\,\sigma_{QE}\,q^{-1} (\nabla_p e^*\tens_0 \extd x^p)}
-\tfrac{\lambda}{2} \,\omega^{ij}\,
\overline{\extd x^k  \tens_1 [\nabla_{k},\nabla_{j}]\nabla_{i}e^* }    \cr
&=& \overline{q^{-1}(\extd x^p \tens_0 \nabla_p e^*+\lambda\,\omega^{ip}\, \extd x^k \tens_0
[\nabla_k,\nabla_i]  \nabla_p e^*  )}
-\tfrac{\lambda}{2} \,\omega^{ij}\,
\overline{\extd x^k  \tens_1 [\nabla_{k},\nabla_{j}]\nabla_{i}e^* }    \cr
&=& \overline{q^{-1}(\extd x^p \tens_0 \nabla_p e^*+\lambda\,\omega^{ij}\, \extd x^k \tens_0
[\nabla_k,\nabla_i]  \nabla_j e^*  )}
-\tfrac{\lambda}{2} \,\omega^{ij}\,
\overline{\extd x^k  \tens_1 [\nabla_{k},\nabla_{j}]\nabla_{i}e^* }    \cr
&=& \overline{q^{-1}(\extd x^p \tens_0 \nabla_p e^*-\lambda\,\omega^{ij}\, \extd x^k \tens_0
[\nabla_k,\nabla_j]  \nabla_i e^*  )}
-\tfrac{\lambda}{2} \,\omega^{ij}\,
\overline{\extd x^k  \tens_1 [\nabla_{k},\nabla_{j}]\nabla_{i}e^* }    \cr
&=& \overline{q^{-1}(\extd x^p \tens_0 \nabla_p e^*-\tfrac\lambda2\,\omega^{ij}\, \extd x^k \tens_0
[\nabla_k,\nabla_j]  \nabla_i e^*  )}\cr
&=& \overline{ \nabla_Q(e^*)  }\ .\quad\square
\end{eqnarray*}

\subsection{Properties of the generalised braiding}
There is an extension of $\sigma_{Q(E)}:E\tens_1 \Omega^1 A_1\to \Omega^1 A_1\tens_1 E$ (as in Theorem~\ref{functor}) to $\sigma_{Q(E)}:E\tens_1 \Omega^n A_1\to \Omega^n A_1\tens_1 E$ given by
\begin{eqnarray*}
 \sigma_{Q(E)}(e\tens_1\xi)\, &=&  \xi\tens_1 e  + \lambda\, \omega^{ij}\,\nabla_j\xi\tens_1\nabla_{Ei}e   + \, \lambda\, \omega^{ij}\,\extd x^k\wedge (\partial_j \,\righthalfcup\,\xi) \tens_1 [\nabla_{Ek},\nabla_{Ei}]e  \ 
\end{eqnarray*}
where $\righthalfcup$ is the interior product. Now we use
\begin{eqnarray} \label{bcuysioiuyf}
[\nabla_k, \nabla_i]\xi &=& -\,R^{s}_{\phantom{s}nki}\,\extd x^n\wedge(\partial_s\, \righthalfcup\,\xi)
\end{eqnarray}
to give, where $E$ is $\Omega^mA_1$ with the compatible covariant derivative,
\begin{eqnarray*}
\sigma_{Q}(\eta\tens_1\xi)\, &=&  \xi\tens_1 \eta  + \lambda\, \omega^{ij}\,\nabla_j\xi\tens_1\nabla_{i}\eta   + \, \lambda\, \omega^{ij}\,\extd x^k\wedge (\partial_j \,\righthalfcup\,\xi) \tens_1 [\nabla_{k},\nabla_{i}]\eta  \cr
&=&   \xi\tens_1 \eta  + \lambda\, \omega^{ij}\,\nabla_j\xi\tens_1\nabla_{i}\eta    -  \,  \lambda\, \omega^{ij}\,R^{s}_{\phantom{s}nki}\,\extd x^k\wedge (\partial_j \,\righthalfcup\,\xi) \tens_1   \extd x^n\wedge(\partial_s\, \righthalfcup\,\eta)\ .
\end{eqnarray*}
Now we can take this twice to give
\begin{eqnarray*}
\sigma_{Q}^2(\eta\tens_1\xi)\, 
&=&   \eta\tens_1 \xi  + \lambda\, \omega^{ij}\,\nabla_{i}\eta\tens_1\nabla_j\xi
+  \lambda\, \omega^{ij}\,\nabla_j\eta\tens_1\nabla_{i}\xi \cr 
&&   -  \,  \lambda\, \omega^{ij}\,R^{s}_{\phantom{s}nki}\,
   \extd x^n\wedge(\partial_s\, \righthalfcup\,\eta)  \tens_1   \extd x^k\wedge (\partial_j \,\righthalfcup\,\xi) \cr
&&   -  \,  \lambda\, \omega^{ij}\,R^{s}_{\phantom{s}nki}\,\extd x^k\wedge (\partial_j \,\righthalfcup\,\eta) \tens_1   \extd x^n\wedge(\partial_s\, \righthalfcup\,\xi) \cr
&=&  \eta\tens_1 \xi  -   \lambda\, (\omega^{ij}\,R^{s}_{\phantom{s}nki} + \omega^{is}\,R^{j}_{\phantom{s}kni}
) \, \extd x^n\wedge(\partial_s\, \righthalfcup\,\eta)  \tens_1   \extd x^k\wedge (\partial_j \,\righthalfcup\,\xi)
\end{eqnarray*}

\begin{proposition}
The generalised braiding obeys the braid relation, in the sense that for any $(E,\nabla_E)$ giving $\sigma_{Q(E)}$ and $(\Omega^1(M),\nabla)$ giving $\sigma_Q$
\begin{eqnarray*}
(\sigma_Q\tens\id)(\id\tens\sigma_{Q(E)})(\sigma_{Q(E)}\tens\id)=
(\id\tens\sigma_{Q(E)})(\sigma_{Q(E)}\tens\id)(\id\tens\sigma_Q)
\end{eqnarray*}
as a map  $E\tens_1 \Omega^1A_1\tens_1  \Omega^1A_1\to \Omega^1A_1\tens_1  \Omega^1A_1 \tens_1 E$.
\end{proposition}
\proof  If we write
\begin{eqnarray*}
\sigma_{Q(E)}(e\tens_1\xi) &=& \xi\tens_1 e+\lambda\,T(\xi)\tens_1 T'(e)\ ,\cr
\sigma_{Q}(\eta\tens_1\xi) &=& \xi\tens_1 \eta+\lambda\,S(\xi)\tens_1 S'(\eta)\ ,
\end{eqnarray*}
then both maps above give the following on being applied to $e\tens_1\xi\tens_1\eta$:
\begin{eqnarray*}
&& \eta\tens_1\xi\tens_1 e+\lambda \, \eta\tens_1T(\xi)\tens_1 T'(e)
+\lambda\,T(\eta)\tens_1\xi\tens_1 T'(e)+\lambda\, S(\eta)\tens_1S'(\xi)\tens_1 e\ .\square
\end{eqnarray*}

\subsection{Quantising other connections relative to a given $(E,\nabla_E)$} \label{jytcxfufgx}

Classically a general covariant derivative is given by $\nabla_S=\nabla_E+S$, where
$S:E\to \Omega^1 (M)\tens E$ is a vector bundle map. This is adding a left module map to a left covariant derivative, giving another covariant derivative on the same bundle.

\begin{proposition}  \label{quantQS} For any  bundle map $S$,
\[
\nabla_{QS}=\nabla_{Q(E)}+q^{-1}Q(S),\quad \sigma_{QS}(e\tens_1\xi)=\sigma_{Q}(e\tens_1\xi)+\lambda\,\omega^{ij}\,\xi_{i}\,\nabla_j(S)(e)\]
defines a bimodule connection on $Q(E)$.
\end{proposition}
\proof Here 
\begin{eqnarray*}
q\,\nabla_{QS}=q\,\nabla_{Q(E)}+S+\tfrac\lambda2\,\omega^{ij}\,\nabla_i\circ\nabla_j(S)\ .
\end{eqnarray*}
From the general equation
\begin{eqnarray*}
\sigma(e\tens_1\extd a) &=& \extd a\tens_1 e + \nabla[e,a]+[a,\nabla e]
\end{eqnarray*}
we see that, to  order $\lambda$
\begin{eqnarray*}
\sigma_{QS}(e\tens_1\extd a) &=& \sigma_{Q}(e\tens_1\extd a)+S([e,a])+[a,S(e)]  \cr
 &=&  \sigma_{Q}(e\tens_1\extd a)-\lambda\,\omega^{ij}\,a_{,i}\,S(\nabla_j e)+\lambda\,\omega^{ij}\,a_{,i}\,\nabla_jS(e)\cr
  &=&  \sigma_{Q}(e\tens_1\extd a)+\lambda\,\omega^{ij}\,a_{,i}\,\nabla_j(S)(e)\ . \quad\square
\end{eqnarray*}

\medskip\noindent
Now we look at the tensor products and reality of such general connections:

\begin{proposition} \label{jdolobvlo}
Given $S:E\to \Omega^1 (M)\tens E$ and 
 $T:F\to \Omega^1 (M)\tens F$, define $H:E\tens F\to \Omega^1 (M)\tens E\tens F$ by (where $\tau$ is transposition)
\begin{eqnarray*}
H=S\tens\id_F+(\tau\tens\id_F)(\id_E\tens T)\ .
\end{eqnarray*}
Then $q^2$ of the tensor product of $\nabla_{QS}$ and $\nabla_{QT}$ is given by
\begin{eqnarray*}
q\,\nabla_{QH}\,q+\lambda\, \mathrm{rem}:Q(E)\tens_1 Q(F)\to Q(\Omega^1(M)\tens E\tens F)
\end{eqnarray*}
where, using $T(f)=\extd x^k\tens T_k(f)$,
\begin{eqnarray*}
\mathrm{rem}(e\tens_1 f) = \omega^{ij}\,( \extd x^k \tens [\nabla_{Ek},\nabla_{Ei}]e
- \nabla_i(S)(e)  )  \tens T_j(f)\ .
\end{eqnarray*}
\end{proposition}
\proof  Using $q^2$ for applying $q$ twice (by Proposition~\ref{q}  order does not matter),
\begin{eqnarray*}
q^2(\nabla_{QS}\tens_1\id_F) &=& q^2((\nabla_{Q} + q^{-1}\,Q(S))\tens_1\id_F) \cr
&=& q^2(\nabla_{Q} \tens_1\id_F) + q^2( q^{-1}\,Q(S)\tens_1\id_F)  \cr
&=& q^2(\nabla_{Q} \tens_1\id_F) + q( Q(S)\tens_1\id_F)  \cr
&=& q^2(\nabla_{Q} \tens_1\id_F) + ( Q(S)\tens\id_F+\tfrac\lambda2\,\omega^{ij}\,\nabla_i(Q(S))\tens \nabla_j)q  \cr
&=& q^2(\nabla_{Q} \tens_1\id_F) + ( Q(S)\tens\id_F+\tfrac\lambda2\,\omega^{ij}\,\nabla_i(S)\tens \nabla_j)q  \ ,
\end{eqnarray*}
where we have used Proposition~\ref{q}. We also need
\begin{eqnarray*}
&&q^2(\sigma_{QS}\tens_1\id_F)(\id_E\tens_1\nabla_{QT}) \cr
 &=& q(q\sigma_{QS}q^{-1}\tens_1\id_F)(q\tens_1\id_F)(\id_E\tens_1\nabla_{QT}) \cr
&=& (q\sigma_{QS}q^{-1}\tens\id_F+\tfrac\lambda2\,\omega^{ij}\,\nabla_i(q\sigma_{QS}q^{-1})\tens\nabla_j)q^2(\id_E\tens_1\nabla_{QT}) \cr
&=& (q\sigma_{QS}q^{-1}\tens\id_F+\tfrac\lambda2\,\omega^{ij}\,\nabla_i(q\sigma_{QS}q^{-1})\tens\nabla_j)q(\id_E\tens_1q\nabla_{QT}) \cr
&=& (q\sigma_{QS}q^{-1}\tens\id_F+\tfrac\lambda2\,\omega^{ij}\,\nabla_i(q\sigma_{QS}q^{-1})\tens\nabla_j)q(\id_E\tens_1(q\nabla_{Q}+Q(T))) \ .
\end{eqnarray*}
Now $\lambda\,\sigma_{QS}=\lambda\,\tau$ to  order $\lambda$, where $\tau$ is transposition, so 
$\tfrac\lambda2\,\nabla_i(q\sigma_{QS}q^{-1})=0$. Then, where we set $q\sigma_{QS}q^{-1}=q\sigma_{Q}q^{-1}+\lambda\,S'$,
\begin{eqnarray*}
&&q^2(\sigma_{QS}\tens_1\id_F)(\id_E\tens_1\nabla_{QT}) \cr
&=& (q\sigma_{QS}q^{-1}\tens\id_F)q\big(\id_E\tens_1q\nabla_{Q}+\id_E\tens_1Q(T)\big) \cr
&=& ((q\sigma_{Q}q^{-1}+\lambda\,S')\tens\id_F)q(\id_E\tens_1q\nabla_{Q})   \cr
&& +\,  ((q\sigma_{Q}q^{-1}+\lambda\,S')\tens\id_F)(\id_E\tens Q(T)+\tfrac\lambda2\,\omega^{ij}\,\nabla_i\tens\nabla_j(Q(T)))\,q \cr
&=& ((q\sigma_{Q}q^{-1}+\lambda\,S')\tens\id_F)q(\id_E\tens_1q\nabla_{Q})   \cr
&& +\,  ((q\sigma_{Q}q^{-1}+\lambda\,S')\tens\id_F)(\id_E\tens Q(T)+\tfrac\lambda2\,\omega^{ij}\,\nabla_i\tens\nabla_j(T))\,q \ .
\end{eqnarray*}
It follows that the contribution of $S$ and $T$ to $q^2$ of the tensor product derivative is
\begin{eqnarray}  \label{bgycoudwuyrdx}
&&( Q(S)\tens\id_F)\,q+\tfrac\lambda2\,\omega^{ij}\,\nabla_i(S)\tens \nabla_j  +  \lambda\,(S'\tens\id_F)(\id_E\tens\nabla)   \cr
&& +\,  (q\sigma_{Q}q^{-1}\tens\id_F)(\id_E\tens Q(T))\,q   +  \tfrac\lambda2\,\omega^{ij}\,(\tau\tens\id_F)(\nabla_i\tens\nabla_j(T))  \cr
&& +\,  \lambda\,(S'\tens\id_F)(\id_E\tens T)    \cr
&=& ( Q(S)\tens\id_F)\,q+\tfrac\lambda2\,\omega^{ij}\,\nabla_i(S)\tens \nabla_j  +  \lambda\,(S'\tens\id_F)(\id_E\tens\nabla_T)   \cr
&& +\,  (q\sigma_{Q}q^{-1}\tens\id_F)(\id_E\tens Q(T))\,q   +  \tfrac\lambda2\,\omega^{ij}\,(\tau\tens\id_F)(\nabla_i\tens\nabla_j(T))  \ .
\end{eqnarray}
From Theorem~\ref{functor} we have
\begin{eqnarray*}
q\, \sigma_{Q(E)}q^{-1}(e\tens\xi)\, &=&  \xi\tens e   +  \lambda\, \omega^{ij}\,\xi_{j}\,\extd x^k \tens [\nabla_{k},\nabla_{i}]e  \ ,
\end{eqnarray*}
and for the moment we write this as $q\, \sigma_{Q(E)}q^{-1}=\tau+\lambda\,\sigma'$. 
Then (\ref{bgycoudwuyrdx}) becomes
\begin{eqnarray}  \label{bgvdfsdx}
&& ( S\tens\id_F)\,q+\tfrac\lambda2\,\omega^{ij}\,\nabla_i(S)\tens \nabla_j  +  \lambda\,(S'\tens\id_F)(\id_E\tens\nabla_T)   \cr
&& +\,  (\tau\tens\id_F)(\id_E\tens Q(T))\,q   +  \tfrac\lambda2\,\omega^{ij}\,(\tau\tens\id_F)(\nabla_i\tens\nabla_j(T))  \cr
&&+\, \tfrac\lambda2\,\omega^{ij}\,\nabla_i\circ\nabla_j(S)\tens\id_F
+\lambda\,(\sigma'\tens\id_F)(\id_E\tens T) \cr
&=&  ( S\tens\id_F)\,q+\tfrac\lambda2\,\omega^{ij}\,\nabla_i(S)\tens \nabla_j  +  \lambda\,(S'\tens\id_F)(\id_E\tens\nabla_T)   \cr
&& +\,  (\tau\tens\id_F)(\id_E\tens T)\,q   +  \tfrac\lambda2\,\omega^{ij}\,(\tau\tens\id_F)(\nabla_i\tens\nabla_j(T))  \cr
&&+\, \tfrac\lambda2\,\omega^{ij}\,\nabla_i\circ\nabla_j(S)\tens\id_F
+\lambda\,(\sigma'\tens\id_F)(\id_E\tens T) \cr
&& +\, \tfrac\lambda2\,\omega^{ij}\, (\tau\tens\id_F)(\id_E\tens \nabla_i\circ\nabla_j(T))\ .
\end{eqnarray}
Now we use $H$ given above
with
\begin{eqnarray*}
\nabla_j(H) &=& \nabla_j(S)\tens\id_F + (\tau\tens\id_F)(\id_E\tens\nabla_j(T))
\end{eqnarray*}
 to write (\ref{bgvdfsdx}) as
\begin{eqnarray*}
 Q(H)\,q+\lambda\,\omega^{ij}\,\nabla_i(S)\tens \nabla_j  +  \lambda\,(S'\tens\id_F)(\id_E\tens\nabla_T)   + \lambda\,(\sigma'\tens\id_F)(\id_E\tens T) \ .
\end{eqnarray*}
Now from Proposition~\ref{quantQS}
\begin{eqnarray*}
(S'\tens\id_F)(\id_E\tens\nabla_T)(e\tens f) &=& (S'\tens\id_F)(e\tens \extd x^k\tens(\nabla_k f+T_k(f))) \cr
&=& \omega^{ij}\,\nabla_j(S)(e)\tens(\nabla_if+T_i(f))
\end{eqnarray*}
so we rewrite  (\ref{bgvdfsdx}) as
\begin{eqnarray}
 Q(H)\,q +\lambda\,\omega^{ij}\,\nabla_j(S)\tens T_i   + \lambda\,(\sigma'\tens\id_F)(\id_E\tens T) \ .
\end{eqnarray}
Finally for $T(f)=\extd x^i\tens T_i(f)$,
\begin{eqnarray*}
(\sigma'\tens\id_F)(\id_E\tens T)(e\tens f) &=& (\sigma'\tens\id_F)(e\tens \extd x^p\tens T_p(f)) \cr
&=& \lambda\,\omega^{ij}\, \extd x^k \tens [\nabla_k,\nabla_i]e\tens T_j(f)\ .\quad\hfill\square
\end{eqnarray*}

\begin{lemma} \label{oiuyfvcftyytrs}
Over $\C$ and if $S$ is real, the difference in going clockwise minus anticlockwise round the diagram
\begin{eqnarray}\label{vcusfsii}
\xymatrix{
\overline{Q(E)} = Q(\overline{E})  \ar[d]_{\overline{\nabla_QS } } 
& Q(E)  \ar[l]_{  \star  }    \ar[r]^{ \nabla_{QS} }  &   
Q(\Omega^1(M)) \tens_1 Q(E)   \ar[d]_{\star\tens_1\star}   \\
\overline{ Q(\Omega^1(M)) \tens_1 Q(E)}&
\overline{ Q(E) \tens_1 Q(\Omega^1(M))}   \ar[l]_{  \overline{\sigma_{QS}}  }  &
\overline{ Q(\Omega^1(M)) } \tens_1 \overline{ Q(E) }   \ar[l]_{\Upsilon^{-1}}  
}
\end{eqnarray}
starting at $Q(e)\in Q(E)$ is
\begin{eqnarray*}
 \overline{\lambda\,\omega^{ij}\, \nabla_j(S)    ( S_i(e^*))    }  -\,  \overline{\lambda\,\omega^{ij}\,\nabla_i(\nabla_j(S))(e^*) }
 +  \overline{\lambda\,    \omega^{ij}\,\extd x^k\tens [\nabla_k,\nabla_i]S_j(e^*)     }  \ .
\end{eqnarray*}
 \end{lemma}
\proof From lemma~\ref{cluiktycxtu} the diagram commutes for $S=0$. We look only at the difference
from the $S=0$ to the general $S$ case. Going anticlockwise from $Q(E)$ we get 
\begin{eqnarray*}
\overline{q^{-1}Q(S)}\star(e) &=& \overline{q^{-1}Q(S)(e^*)}
\end{eqnarray*}
Going clockwise is more complicated, as two of the arrows involve $S$. If we set $q\sigma_{QS}q^{-1}=q\sigma_{Q}q^{-1}+\lambda\,S'$ as in the proof of Proposition~\ref{jdolobvlo}, then to  order $\lambda$ we get the clockwise contributions
\begin{eqnarray*}
&& \overline{\sigma_{Q}} \,\Upsilon^{-1}(\star\tens_1\star)q^{-1}Q(S)(Q(e))+\overline{\lambda\,S'}\, \Upsilon^{-1}(\star\tens\star)\,\nabla_S(e)  \cr
&=& \overline{\sigma_{Q}} \,\Upsilon^{-1}\,q^{-1}(\star\tens_0\star)Q(S)(Q(e))+\overline{\lambda\,S'}\, \Upsilon^{-1}(\star\tens\star)\,\nabla_S(e)  \cr
&=& \overline{\sigma_{Q}\,q^{-1}} \,\Upsilon^{-1}(\star\tens_0\star)Q(S)(Q(e))+\overline{\lambda\,S'}\, \Upsilon^{-1}(\star\tens\star)\,\nabla_S(e)  \cr
&=& \overline{q^{-1}\,q\,\sigma_{Q}\,q^{-1}} \,\Upsilon^{-1}(\star\tens_0\star)Q(S)(Q(e))+\overline{\lambda\,S'}\, \Upsilon^{-1}(\star\tens\star)\,\nabla_S(e)  \cr
&=& \overline{q^{-1}\,q\,\sigma_{Q}\,q^{-1}} \,\Upsilon^{-1}(\star\tens_0\star)( S(e)+\tfrac\lambda2\,\omega^{ij}\,\nabla_i(\nabla_j(S)(e)) )
\cr && + \, \overline{\lambda\,S'}\, \Upsilon^{-1}(\star\tens\star)\,(\extd x^p\tens\nabla_p e + S(e)) \cr
&=& \overline{q^{-1}\,\tau} \,\Upsilon^{-1}(\star\tens_0\star)( S(e)+\tfrac\lambda2\,\omega^{ij}\,\nabla_i(\nabla_j(S)(e)) )\cr 
&& +\,  \overline{\lambda\,q^{-1}\,\sigma'} \,\Upsilon^{-1}(\star\tens_0\star)( S(e) )\cr 
&& + \, \overline{\lambda\,S'}\, \Upsilon^{-1}(\star\tens\star)\,(\extd x^p\tens\nabla_p e + S(e)) \ ,
\end{eqnarray*}
Where we have put
$q\, \sigma_{Q}q^{-1}=\tau+\lambda\,\sigma'$. As the classical connections preserve $\star$ and $\lambda$ is imaginary, we get the following for the clockwise contributions, where $S(e)=\extd x^p\tens S_p(e)$,
\begin{eqnarray*}
&& \overline{q^{-1}\,( S(e^*)-\tfrac\lambda2\,\omega^{ij}\,\nabla_i(\nabla_j(S)(e^*)) )}
 +  \overline{\lambda\,q^{-1}\,\sigma'} \,\Upsilon^{-1}(\star\tens_0\star)( \extd x^p\tens S_p(e) )\cr 
&& + \, \overline{\lambda\,S'}\, \Upsilon^{-1}(\star\tens\star)\,(\extd x^p\tens\nabla_p e + S(e)) \cr
&=& \overline{q^{-1}\,( S(e^*)-\tfrac\lambda2\,\omega^{ij}\,\nabla_i(\nabla_j(S)(e^*)) )}
 +  \overline{\lambda\,q^{-1}\,\sigma'} \,\Upsilon^{-1}( \overline{\extd x^p}\tens \overline{ S_p(e^*)} )\cr 
&& + \, \overline{\lambda\,S'}\, \Upsilon^{-1}(\star\tens\star)\,(\extd x^p\tens\nabla_p e + S(e)) \cr
&=& \overline{q^{-1}\,( S(e^*)-\tfrac\lambda2\,\omega^{ij}\,\nabla_i(\nabla_j(S)(e^*)) )}
 +  \overline{\lambda\,q^{-1}\,\sigma'( S_p(e^*) \tens \extd x^p     )}  \cr 
&& + \, \overline{\lambda\,S'}\, \Upsilon^{-1}(\star\tens\star)\,(\extd x^p\tens(\nabla_p e + S_p(e))   ) \cr
&=& \overline{q^{-1}\,( S(e^*)-\tfrac\lambda2\,\omega^{ij}\,\nabla_i(\nabla_j(S)(e^*)) )}
 +  \overline{\lambda\,q^{-1}\,(    \omega^{ip}\,\extd x^k\tens [\nabla_k,\nabla_i]S_p(e^*)     )}  \cr 
&& + \, \overline{\lambda\,S' (     (\nabla_p e^* + S_p(e^*))  \tens\extd x^p      )}\ \cr
&=& \overline{q^{-1}\,( S(e^*)-\tfrac\lambda2\,\omega^{ij}\,\nabla_i(\nabla_j(S)(e^*)) )}
 +  \overline{\lambda\,q^{-1}\,(    \omega^{ip}\,\extd x^k\tens [\nabla_k,\nabla_i]S_p(e^*)     )}  \cr 
&& + \, \overline{\lambda\,\omega^{pj}\, \nabla_j(S)    (\nabla_p e^* + S_p(e^*))    }\ \ .
\end{eqnarray*}
Then the difference, clockwise minus anticlockwise, is to  order $\lambda$,
\begin{eqnarray*}
&& -\,  \overline{\lambda\,\omega^{ij}\,\nabla_i(\nabla_j(S)(e^*)) }
 +  \overline{\lambda\,    \omega^{ip}\,\extd x^k\tens [\nabla_k,\nabla_i]S_p(e^*)     }  \cr 
&& + \, \overline{\lambda\,\omega^{pj}\, \nabla_j(S)    (\nabla_p e^* + S_p(e^*))    }   \cr
&=& -\,  \overline{\lambda\,\omega^{ij}\,\nabla_i(\nabla_j(S)(e^*)) }
 +  \overline{\lambda\,    \omega^{ij}\,\extd x^k\tens [\nabla_k,\nabla_i]S_j(e^*)     }  \cr 
&& + \, \overline{\lambda\,\omega^{ij}\, \nabla_j(S)    (\nabla_i e^* + S_i(e^*))    }   \cr
&=& \overline{\lambda\,\omega^{ij}\, \nabla_j(S)    ( S_i(e^*))    }  -\,  \overline{\lambda\,\omega^{ij}\,\nabla_i(\nabla_j(S))(e^*) }
 +  \overline{\lambda\,    \omega^{ij}\,\extd x^k\tens [\nabla_k,\nabla_i]S_j(e^*)     }  \ .\quad\square
\end{eqnarray*}

Note that Lemma~\ref{oiuyfvcftyytrs} shows that $\nabla_{QS}(\star)$ is $\lambda$ times a module map (i.e.\ it involves no derivatives of $e$). This means that $\nabla_{QS}(\star)$ is also a right module map, and thus it is automatically star-compatible at order $\lambda$ in the sense described in \cite{BegMa3}. We will also need the following Lemma.

\begin{lemma}  \label{hilskljhvcj} Let $(E,\nabla_E)$ be a bundle with connection and  $e\in E$ such that $\nabla(e)=0$. Then $e$ is central in the quantised bimodule, the quantised derivative $\nabla_Q(e)=0$. If in addition $T(e)=0$ for some $T:E\to\Omega^1(M)\tens_0 E$ then $\nabla_{QT}(e)=0$. 
\end{lemma}
\vspace{-9pt}
\proof  
For the quantised connection,
\begin{eqnarray*}
\nabla_{Q}(e) &=& q^{-1}\nabla(e)    -     \tfrac{\lambda}{2} \,\omega^{ij}\,\extd x^k\tens_1[\nabla_{k},\nabla_{j}]\nabla_{i}(e) =0\ .
\end{eqnarray*}
Then classically 
$\nabla_T(e)=\nabla(e)+T(e)=0$. Now
\begin{eqnarray*}
\nabla_j(T)(e)=\nabla_j(T(e))-T(\nabla_j(e))=0\ ,
\end{eqnarray*}
so in the quantised case
\begin{eqnarray*}
Q(T)(e) &=& T(e)+\tfrac\lambda2 \, \omega^{ij} \, \nabla_i\circ\nabla_j(T)(e)=0\ .
\end{eqnarray*}
Now we have, for the quantisation $\nabla_{QT}$ of $\nabla_T$,
\begin{eqnarray*}
\nabla_{QT}(e)=\nabla_{Q}(e)+q^{-1}\,Q(T)(e)=0\ .\quad\square
\end{eqnarray*}

\section{Semiquantisation of the exterior algebra}

In noncommutative geometry the notion of `differential structure' is largely encoded as a differential graded algebra extending the quantisation of functions to differential forms. The main result in this section is that by adapting the semiquantisation functor of Section~3 we have from the same data a canonical semiquantisation of forms of all degree and their wedge product. This is Theorem~\ref{dga}. 

\subsection{Quantizing the wedge product}

Our starting point for the quantum wedge product is the associative product which is given by the functor $Q$, where we assume that the connection $\nabla$ on $\Omega^1A$ extends to all orders in the natural way, is the composition
\begin{eqnarray*}
Q(\Omega^1(M)) \tens_1 Q(\Omega^1(M)) \stackrel{q} \longrightarrow Q(\Omega^1(M)
\tens_0 \Omega^1(M)) \stackrel{Q(\wedge)} \longrightarrow Q(\Omega^2(M))\ .
\end{eqnarray*}
This gives the formula for $\wedge_Q$;
\begin{eqnarray}  \label{piuvgcghjk}
\xi\wedge_Q\eta \,=\,\xi\wedge\eta+\tfrac{\lambda}{2} \omega^{ij}\,\nabla_i\xi\wedge\nabla_j\eta \ .
\end{eqnarray}

To look at the Leibniz rule we need the following result:

\begin{lemma}  \label{vghsuioi}
\begin{eqnarray*}\extd(\xi\wedge_Q\eta)-(\extd\xi)\wedge_Q\eta
-(-1)^{|\xi|}\xi\wedge_Q\extd\eta&=&-\lambda H^{ji}\wedge(\del_i\righthalfcup\xi)\wedge\nabla_j\eta
\cr && +\ \lambda(-1)^{|\xi|}H^{ij}\wedge\nabla_i\xi\wedge(\del_j\righthalfcup\eta)
\end{eqnarray*}
where 
\[ H^{ij}:=\tfrac14\omega^{is}\left( T^j_{nm;s}-2 R^j{}_{nms}\right)\extd x^m\wedge\extd x^n.\]
\end{lemma}
\proof Using (\ref{poissoncomp}) in the following form
\begin{eqnarray*} \label{poissoncacomp}
\extd(\omega^{ij})  - \omega^{kj}\,\nabla_k(\extd x^i)
-\omega^{ik}\, \nabla_k(\extd x^j) \,=\,0\ ,
\end{eqnarray*}
and also using
\begin{eqnarray*}
\extd\zeta\,=\,\extd x^k\wedge\nabla_k\zeta+ \tfrac12\,T^s_{kn}\,\extd x^k\wedge\extd x^n\wedge (\partial_s\, \righthalfcup\,\zeta)
\end{eqnarray*}
and relabeling indices, we find
\begin{eqnarray*}
\extd(\omega^{ij}\,\nabla_i\xi\wedge\nabla_j\eta) &=& \omega^{ij}\,\nabla_i(\extd x^k)\wedge \nabla_k\xi\wedge\nabla_j\eta + \omega^{ij}\, \nabla_j(\extd x^k)\wedge \nabla_i\xi\wedge\nabla_k\eta \cr
&& +\, \omega^{ij}\,\extd x^k\wedge\nabla_k \nabla_i\xi\wedge\nabla_j\eta+ 
(-1)^{|\xi|}\,\omega^{ij}\,\nabla_i\xi\wedge\extd x^k\wedge\nabla_k\nabla_j\eta  \cr
&&  +\, \omega^{ij}\,
\tfrac12\,T^s_{kn}\,\extd x^k\wedge\extd x^n\wedge (\partial_s\, \righthalfcup\,\nabla_i\xi)\wedge\nabla_j\eta  \cr
&& + \,
(-1)^{|\xi|}\,\omega^{ij}\,\nabla_i\xi\wedge  \tfrac12\,T^s_{kn}\,\extd x^k\wedge\extd x^n\wedge (\partial_s\, \righthalfcup\,\nabla_j\eta )  \cr
&=& \omega^{ij}\,\nabla_i(\extd x^k\wedge \nabla_k\xi)\wedge\nabla_j\eta +(-1)^{|\xi|}\, \omega^{ij}\, \nabla_i\xi\wedge\nabla_j(\extd x^k\wedge \nabla_k\eta) \cr
&& +\, \omega^{ij}\,\extd x^k\wedge[\nabla_k, \nabla_i]\xi\wedge\nabla_j\eta+ 
(-1)^{|\xi|}\,\omega^{ij}\,\nabla_i\xi\wedge\extd x^k\wedge[\nabla_k,\nabla_j]\eta  \cr
&&  +\, \omega^{ij}\,
\tfrac12\,T^s_{kn}\,\extd x^k\wedge\extd x^n\wedge (\partial_s\, \righthalfcup\,\nabla_i\xi)\wedge\nabla_j\eta  \cr
&& + \,
(-1)^{|\xi|}\,\omega^{ij}\,\nabla_i\xi\wedge  \tfrac12\,T^s_{kn}\,\extd x^k\wedge\extd x^n\wedge (\partial_s\, \righthalfcup\,\nabla_j\eta )\ .
\end{eqnarray*}
From this we get
\begin{eqnarray*}
&&\extd(\omega^{ij}\,\nabla_i\xi\wedge\nabla_j\eta)-\omega^{ij}\,\nabla_i\extd\xi\wedge\nabla_j\eta
-(-1)^{|\xi|}\, \omega^{ij}\,\nabla_i\xi\wedge\nabla_j\extd\eta \cr
&=&-\, \omega^{ij}\,\nabla_i(\tfrac12\,T^s_{kn}\,\extd x^k\wedge\extd x^n\wedge (\partial_s\, \righthalfcup\,\xi))\wedge\nabla_j\eta \cr
&& - \, (-1)^{|\xi|}\, \omega^{ij}\, \nabla_i\xi\wedge\nabla_j(\tfrac12\,T^s_{kn}\,\extd x^k\wedge\extd x^n\wedge (\partial_s\, \righthalfcup\,\eta)) \cr
&& +\, \omega^{ij}\,\extd x^k\wedge[\nabla_k, \nabla_i]\xi\wedge\nabla_j\eta+ 
(-1)^{|\xi|}\,\omega^{ij}\,\nabla_i\xi\wedge\extd x^k\wedge[\nabla_k,\nabla_j]\eta  \cr
&&  +\, \omega^{ij}\,
\tfrac12\,T^s_{kn}\,\extd x^k\wedge\extd x^n\wedge (\partial_s\, \righthalfcup\,\nabla_i\xi)\wedge\nabla_j\eta  \cr
&& + \,
(-1)^{|\xi|}\,\omega^{ij}\,\nabla_i\xi\wedge  \tfrac12\,T^s_{kn}\,\extd x^k\wedge\extd x^n\wedge (\partial_s\, \righthalfcup\,\nabla_j\eta )\ .
\end{eqnarray*}
As $\nabla_i(v\, \righthalfcup\,\xi)=\nabla_i(v)\, \righthalfcup\,\xi+v\, \righthalfcup\,\nabla_i\xi$, we get
\begin{eqnarray*}
&&\extd(\omega^{ij}\,\nabla_i\xi\wedge\nabla_j\eta)-\omega^{ij}\,\nabla_i\extd\xi\wedge\nabla_j\eta
-(-1)^{|\xi|}\, \omega^{ij}\,\nabla_i\xi\wedge\nabla_j\extd\eta \cr
&=&-\, \omega^{ij}\,\tfrac12\,T^s_{kn;i}\,\extd x^k\wedge\extd x^n\wedge (\partial_s\, \righthalfcup\,\xi)\wedge\nabla_j\eta \cr
&& - \, (-1)^{|\xi|}\, \omega^{ij}\, \nabla_i\xi\wedge\tfrac12\,T^s_{kn;j}\,\extd x^k\wedge\extd x^n\wedge (\partial_s\, \righthalfcup\,\eta) \cr
&& +\, \omega^{ij}\,\extd x^k\wedge[\nabla_k, \nabla_i]\xi\wedge\nabla_j\eta+ 
(-1)^{|\xi|}\,\omega^{ij}\,\nabla_i\xi\wedge\extd x^k\wedge[\nabla_k,\nabla_j]\eta  \ .
\end{eqnarray*}
Now we use (\ref{bcuysioiuyf})
to give
\begin{eqnarray*}
&&\extd(\omega^{ij}\,\nabla_i\xi\wedge\nabla_j\eta)-\omega^{ij}\,\nabla_i\extd\xi\wedge\nabla_j\eta
-(-1)^{|\xi|}\, \omega^{ij}\,\nabla_i\xi\wedge\nabla_j\extd\eta \cr
&=&-\, \omega^{ij}\,\tfrac12\,T^s_{kn;i}\,\extd x^k\wedge\extd x^n\wedge (\partial_s\, \righthalfcup\,\xi)\wedge\nabla_j\eta \cr
&& - \, (-1)^{|\xi|}\, \omega^{ij}\, \nabla_i\xi\wedge\tfrac12\,T^s_{kn;j}\,\extd x^k\wedge\extd x^n\wedge (\partial_s\, \righthalfcup\,\eta) \cr
&& -\, \omega^{ij}\,\extd x^k\wedge R^{s}_{\phantom{s}nki}\,\extd x^n\wedge(\partial_s\, \righthalfcup\,\xi) \wedge\nabla_j\eta  \cr
&& -  \, 
(-1)^{|\xi|}\,\omega^{ij}\,\nabla_i\xi\wedge\extd x^k\wedge R^{s}_{\phantom{s}nkj}\,\extd x^n\wedge(\partial_s\, \righthalfcup\,\eta)
  \ . \qquad\square
\end{eqnarray*}

We see that $\wedge_Q$ will not in general obey the Leibniz rule for the undeformed $\extd$. We have a choice of persisting with a modified Leibniz rule perhaps linking up to examples such as \cite{Delius,FaddPyat}, or modifying the wedge product. We choose the second more conventional option:

\begin{lemma} \label{hcuioiutydfdyu}
For vector field $v$, $\xi\in \Omega^n(M)$ and covariant derivative $\nabla$,
\begin{eqnarray*}
v\, \righthalfcup\,\extd\xi +\extd(v\, \righthalfcup\,\xi)&=& v\, \righthalfcup\,\nabla\xi+\wedge(\nabla(v)\,\righthalfcup\,\xi)+v^j\,T^k_{ji}\,\extd x^i \wedge(\partial_k\, \righthalfcup\,\xi)\ .
\end{eqnarray*}
(The left hand side here is the usual Lie derivative).
\end{lemma}
\proof First we start with a 1-form $\xi$, when
\begin{eqnarray*}
v\, \righthalfcup\,\extd(\xi_i\,\extd x^i)+\extd(v\, \righthalfcup\,\xi_i\,\extd x^i) &=& v^j\,\xi_{i,j}\,\extd x^i-v^i\,\xi_{i,j}\,\extd x^j
+ v^i\,\xi_{i,j}\,\extd x^j+v^i_{\phantom{i},j}\,\xi_{i}\,\extd x^j \cr
&=& v^j\,(\xi_{i,j} - \Gamma_{ji}^k\,\xi_k  )\,\extd x^i+(v^i_{\phantom{i},j}  + \Gamma^i_{jk}\,v^k  )\,\xi_{i}\,\extd x^j +v^j\,T_{ji}^k\,\xi_k  \,\extd x^i \ .
\end{eqnarray*}
Now we extend this by induction, for $\xi\in\Omega^1(M)$,
\begin{eqnarray*}
v\, \righthalfcup\,\extd(\xi\wedge\eta)&=& v\, \righthalfcup\,(\extd\xi\wedge\eta-\xi\wedge\extd\eta)\cr
&=&  (v\, \righthalfcup\,\extd\xi)\wedge\eta + \extd\xi\wedge(v\, \righthalfcup\,\eta)     -(v\, \righthalfcup\,\xi)\,\extd\eta
+\xi\wedge (v\, \righthalfcup\,\extd\eta)\ ,\cr
\extd(v\, \righthalfcup\,(\xi\wedge\eta)) &=& \extd((v\, \righthalfcup\,\xi)\wedge\eta)-
\extd(\xi\wedge (v\, \righthalfcup\,\eta)) \cr&=& \extd(v\, \righthalfcup\,\xi)\wedge\eta
+(v\, \righthalfcup\,\xi)\,\extd(\eta)-
\extd \xi\wedge (v\, \righthalfcup\,\eta)+
\xi\wedge \extd(v\, \righthalfcup\,\eta)\ .
\end{eqnarray*}
Then, assuming the $\eta\in\Omega^n(M)$ and that the result works for $n$,
\begin{eqnarray*}
&& v\, \righthalfcup\,\extd(\xi\wedge\eta)+\extd(v\, \righthalfcup\,(\xi\wedge\eta))\cr
&=&  \big(\extd(v\, \righthalfcup\,\xi)+ v\, \righthalfcup\,\extd\xi  \big)\wedge\eta
+ \xi\wedge \big(v\, \righthalfcup\,\extd\eta   +\extd(v\, \righthalfcup\,\eta)     \big) \cr
&=&  \big(v\, \righthalfcup\,\nabla\xi+\wedge(\nabla(v)\,\righthalfcup\,\xi)+v^j\,T^k_{ji}\,\extd x^i \wedge(\partial_k\, \righthalfcup\,\xi)\big)\wedge\eta \cr
&& +\, \xi\wedge \big(  v\, \righthalfcup\,\nabla\eta+\wedge(\nabla(v)\,\righthalfcup\,\eta)   +v^j\,T^k_{ji}\,\extd x^i \wedge(\partial_k\, \righthalfcup\,\eta) \big)\ .\quad\square
\end{eqnarray*}

\begin{proposition} \label{wedge1leib}
Let $H^{ij}$ be as in Lemma~\ref{vghsuioi}. Then
\[ \xi\wedge_1\eta=\xi\wedge_Q\eta+\lambda\,(-1)^{|\xi|+1}\, H^{ij}\wedge (\partial_i \, \righthalfcup\, \xi)\wedge(\partial_j \, \righthalfcup\, \eta)\ \]
is associative to order $\lambda$ and the Leibniz rule holds to order $\lambda$ if and only if
\[  H^{ij}=H^{ji},\quad  \extd H^{ij}+\Gamma^i_{rp}\,\extd x^p\wedge H^{rj} +\Gamma^j_{rp}\,\extd x^p\wedge H^{ir}=0\,\quad\forall i,j.\]
\end{proposition}
\proof We write
\[ \xi\wedge_1\eta=\xi\wedge_Q\eta+\lambda\, \xi\widehat\wedge\eta\]
where 
\begin{eqnarray*}
\xi\,\widehat\wedge\,\eta\,=\,(-1)^{|\xi|+1}\, H^{ij}\wedge (\partial_i \, \righthalfcup\, \xi)\wedge(\partial_j \, \righthalfcup\, \eta)\ ,
\end{eqnarray*}
and for the moment $H^{ij}$ is an arbitrary collection of 2-forms (the first part holds in general).  For the first part, we compute
\begin{eqnarray*}
(\xi\,\widehat\wedge\,\eta)\wedge\zeta &=&(-1)^{|\xi|+1}\, H^{ij}\wedge (\partial_i \, \righthalfcup\, \xi)\wedge(\partial_j \, \righthalfcup\, \eta)\wedge\zeta\ ,\cr
(\xi\wedge\eta)\,\widehat\wedge\,\zeta&=&(-1)^{|\xi|+|\eta|+1}\, H^{ij}\wedge (\partial_i \, \righthalfcup\, (\xi\wedge\eta))\wedge(\partial_j \, \righthalfcup\, \zeta)\cr
&=&(-1)^{|\xi|+|\eta|+1}\, H^{ij}\wedge (\partial_i \, \righthalfcup\, \xi)\wedge\eta\wedge(\partial_j \, \righthalfcup\, \zeta)\cr
&& +\, (-1)^{|\eta|+1}\, H^{ij}\wedge \xi\wedge(\partial_i \, \righthalfcup\, \eta)\wedge(\partial_j \, \righthalfcup\, \zeta)\ ,
\end{eqnarray*}
and
\begin{eqnarray*}
\xi\wedge (\eta\,\widehat\wedge\,\zeta) &=&(-1)^{|\eta|+1}\,\xi\wedge H^{ij}\wedge (\partial_i \, \righthalfcup\, \eta)\wedge(\partial_j \, \righthalfcup\, \zeta)\ ,\cr
\xi\,\widehat\wedge\, (\eta\wedge\zeta)&=&(-1)^{|\xi|+1}\, H^{ij}\wedge (\partial_i \, \righthalfcup\, \xi)\wedge(\partial_j \, \righthalfcup\, (\eta\wedge\zeta))\cr
&=&(-1)^{|\xi|+|\eta|+1}\, H^{ij}\wedge (\partial_i \, \righthalfcup\, \xi)\wedge\eta\wedge(\partial_j \, \righthalfcup\, \zeta)\cr
&&+\, (-1)^{|\xi|+1}\, H^{ij}\wedge (\partial_i \, \righthalfcup\, \xi)\wedge(\partial_j \, \righthalfcup\, \eta)\wedge\zeta\ . 
\end{eqnarray*}
Hence
\begin{eqnarray*}
(\xi\,\widehat\wedge\,\eta)\wedge\zeta + (\xi\wedge\eta)\,\widehat\wedge\,\zeta
=\xi\wedge (\eta\,\widehat\wedge\,\zeta) + \xi\,\widehat\wedge\, (\eta\wedge\zeta)\ ,
\end{eqnarray*}
which given that $\wedge_Q$ is necessarily associative to order $\lambda$ by functoriality gives the result stated.

Next, using again the given definition of $\xi\,\widehat\wedge\,\eta$,
\begin{eqnarray*}
\extd(\xi\,\widehat\wedge\,\eta) &=& (-1)^{|\xi|+1}\,  \extd H^{ij}\wedge (\partial_i \, \righthalfcup\, \xi)\wedge(\partial_j \, \righthalfcup\, \eta)
+(-1)^{|\xi|+1}\, H^{ij}\wedge \extd(\partial_i \, \righthalfcup\, \xi)\wedge(\partial_j\, \righthalfcup\, \eta) \cr
&&+\, H^{ij}\wedge (\partial_i \, \righthalfcup\, \xi)\wedge\extd(\partial_j \, \righthalfcup\, \eta) \cr
&=& (-1)^{|\xi|}\,  (\Gamma^i_{rt}\,\extd x^t\wedge H^{rj} +  \Gamma^j_{rt}\,\extd x^t\wedge H^{ir}-G^{ij})\wedge (\partial_i \, \righthalfcup\, \xi)\wedge(\partial_j \, \righthalfcup\, \eta)
\cr && +\,   (-1)^{|\xi|+1}\, H^{ij}\wedge \extd(\partial_i \, \righthalfcup\, \xi)\wedge(\partial_j\, \righthalfcup\, \eta) 
+ H^{ij}\wedge (\partial_i \, \righthalfcup\, \xi)\wedge\extd(\partial_j \, \righthalfcup\, \eta) \ .
\end{eqnarray*}
From Lemma~\ref{hcuioiutydfdyu} we use
\begin{eqnarray*}
\partial_j\, \righthalfcup\,\extd\xi +\extd(\partial_j\, \righthalfcup\,\xi)&=& \partial_j\, \righthalfcup\,\nabla\xi+\wedge(\nabla(\partial_j)\,\righthalfcup\,\xi)+T^k_{jt}\,\extd x^t \wedge(\partial_k\, \righthalfcup\,\xi)\cr
&=& \nabla_j\xi+
\extd x^t\wedge \Gamma^s_{tj}\, (\partial_s\,\righthalfcup\,\xi)+T^s_{jt}\,\extd x^t \wedge(\partial_s\, \righthalfcup\,\xi) \cr
&=& \nabla_j\xi+
\extd x^t\wedge \Gamma^s_{jt}\, (\partial_s\,\righthalfcup\,\xi)
\end{eqnarray*}
to give
\begin{eqnarray*}
\extd(\xi\,\widehat\wedge\,\eta) 
&=& (-1)^{|\xi|}\,  (\Gamma^i_{rt}\,\extd x^t\wedge H^{rj} +  \Gamma^j_{rt}\,\extd x^t\wedge H^{ir}-G^{ij})\wedge (\partial_i \, \righthalfcup\, \xi)\wedge(\partial_j \, \righthalfcup\, \eta)
\cr && +\,   (-1)^{|\xi|+1}\, H^{ij}\wedge
\big( \nabla_i\xi+
\extd x^t\wedge \Gamma^s_{it}\, (\partial_s\,\righthalfcup\,\xi)-\partial_i\, \righthalfcup\,\extd\xi  \big) \wedge(\partial_j\, \righthalfcup\, \eta) \cr
&&
+\,  H^{ij}\wedge (\partial_i \, \righthalfcup\, \xi)\wedge\big( \nabla_j\eta+
\extd x^t\wedge \Gamma^s_{jt}\, (\partial_s\,\righthalfcup\,\eta)-\partial_j\, \righthalfcup\,\extd\eta  \big) \ .
\end{eqnarray*}
Comparing these fragments, we find
\begin{eqnarray*}
 \extd(\xi\,\widehat\wedge\,\eta)- \extd(\xi)\,\widehat\wedge\,\eta & - & (-1)^{|\xi|}\, \xi\,\widehat\wedge\,\extd(\eta) \cr
&=& H^{ij}\,\wedge (\partial_i\, \righthalfcup\, \xi)\wedge\nabla_j\eta- (-1)^{|\xi|}\, 
H^{ij}\,\wedge \nabla_i\xi \wedge (\partial_j\, \righthalfcup\, \eta)\cr
&& -\,  (-1)^{|\xi|}\,  G^{ij}\wedge (\partial_i \, \righthalfcup\, \xi)\wedge(\partial_j \, \righthalfcup\, \eta)\ ,
\end{eqnarray*}
where $G^{ij}:=\extd H^{ij}+\Gamma^i_{rp}\,\extd x^p\wedge H^{rj} +\Gamma^j_{rp}\,\extd x^p\wedge H^{ir}$. Again, this expression holds for any collection $H^{ij}$.

Now comparing with Lemma~\ref{vghsuioi} and taking $H^{ij}$ as defined there, we see that the Leibniz rule holds with respect to $\wedge_1$ if and only if $H^{ij}$ is symmetric and
$G^{ij}=0$. To see that these have to hold separately, one may take $\eta$ in degree 0 so that the interior product $\del_j\righthalfcup\eta=0$. \endproof
 
This gives conditions on the curvature and torsion contained in $H^{ij}$ to obtain a differential graded algebra to order $\lambda$. 

\subsection{Results on curvature, torsion and the tensor $N$}

Here we do some calculations in Riemannian geometry with torsion in order to simplify our two conditions in Proposition~\ref{wedge1leib} on the tensor $H^{ij}$.  We use \cite{DelBianchi} and \cite{Str} for the Bianchi identities with torsion;
\begin{eqnarray*} \label{bianchitorbchals}
(B1) && \sum_{\mathrm{cyclic\,permutations}(abc)} \Big(T^k_{bc;a}-R^k_{\phantom{k}abc}-T^k_{ai}\,T^i_{bc}\Big)\ =\ 0\ ,\cr
(B2) && \sum_{\mathrm{cyclic\,permutations}(abc)} \Big(R^k_{\phantom{k}jbc;a}-R^k_{\phantom{k}jai}\,T^i_{bc}\Big)\ =\ 0\ .
\end{eqnarray*}
We also have to bring out a technical point of the semicolon equals covariant derivative notation, which only occurs if we use it more than once. For some tensor $K$ (with various indices), we have $K_{;i}=\nabla_i K$, but $K_{;ij}\neq\nabla_j\nabla_i K$. This is because in $K_{;ij}$ we take the $j$th covariant derivative of
$K_{;i}$ including $i$ in the indices we take the covariant derivative with respect to, so we get an extra term $-\Gamma_{ji}^p\,K_{;p}$ which does not appear in $\nabla_j\nabla_i K$. 
Thus in the presence of torsion we get $K_{;ij}-K_{;ji}=[\nabla_j,\nabla_i ]K-T_{ji}^p\,K_{;p}$, where the commutator $[\nabla_j,\nabla_i ]$ gives the curvature. 

\begin{lemma}
Given the compatibility condition (\ref{compatT}), the 2-forms $H^{ij}$ in Lemma~\ref{vghsuioi} obey $H^{ij}=H^{ji}$.
\end{lemma}
\proof  Differentiate the compatibility condition to get
\begin{eqnarray*}
0 \,&=&\, \omega^{ij}_{\phantom{ij}; mn}+\omega^{ik}_{\phantom{ik};n}\,T^j_{km}+\omega^{kj}_{\phantom{ik};n}\,T^i_{km}
+\omega^{ik}\,T^j_{km;n}+\omega^{kj}\,T^i_{km;n} \cr
&=&\, \omega^{ij}_{\phantom{ij}; mn}-(\omega^{is}\,T^k_{sn}+\omega^{sk}\,T^i_{sn})\,T^j_{km}-(\omega^{ks}\,T^j_{sn}+\omega^{sj}\,T^k_{sn})\,T^i_{km}
+\omega^{ik}\,T^j_{km;n}+\omega^{kj}\,T^i_{km;n} 
\end{eqnarray*}
which we rearrange as
\begin{eqnarray*}
\omega^{ij}_{\phantom{ij}; mn}\,&=&\, \omega^{is}\,T^k_{sn}\,T^j_{km}+\omega^{sj}\,T^k_{sn}\,T^i_{km}+\omega^{sk}(T^i_{sn}\,T^j_{km}+T^i_{sm}\,T^j_{kn})
-\omega^{ik}\,T^j_{km;n}-\omega^{kj}\,T^i_{km;n} 
\end{eqnarray*}
Now use
\begin{eqnarray*}
\omega^{ij}_{\phantom{ij}; mn}-\omega^{ij}_{\phantom{ij}; nm} &=& \omega^{sj} \,R^i_{\phantom{i}snm}
+ \omega^{is} \,R^j_{\phantom{i}snm}-T^p_{nm}\,\omega^{ij}{}_{;p}   \cr
&=& \omega^{sj} \,R^i_{\phantom{i}snm}
+ \omega^{is} \,R^j_{\phantom{i}snm}+T^p_{nm}\,(\omega^{ik}T^j_{kp}+\omega^{kj}T^i_{kp})\ ,
\end{eqnarray*}
where we have used the compatibility condition again, to get
\begin{eqnarray*}
 \omega^{sj} \,R^i_{\phantom{i}snm}
+ \omega^{is} \,R^j_{\phantom{i}snm}\,&=&\, \omega^{is}(T^k_{sn}\,T^j_{km}-T^k_{sm}\,T^j_{kn})+\omega^{sj}(T^k_{sn}\,T^i_{km}-T^k_{sm}\,T^i_{kn}) \cr
&&\,
-\,  \omega^{is}(T^j_{sm;n}-T^j_{sn;m})-\omega^{sj}(T^i_{sm;n} -T^i_{sn;m} ) \cr
&&-\, T^p_{nm}\,(\omega^{ik}T^j_{kp}+\omega^{kj}T^i_{kp})   \cr
&=& \, \omega^{is}(T^k_{sn}\,T^j_{km}-T^k_{sm}\,T^j_{kn})+\omega^{sj}(T^k_{sn}\,T^i_{km}-T^k_{sm}\,T^i_{kn}) \cr
&&\,
-\,  \omega^{is}(T^j_{sm;n}-T^j_{sn;m})-\omega^{sj}(T^i_{sm;n} -T^i_{sn;m} ) \cr
&&-\, T^k_{nm}\,(\omega^{is}T^j_{sk}+\omega^{sj}T^i_{sk})   \cr
&=& \, \omega^{is}(T^k_{sn}\,T^j_{km}-T^k_{sm}\,T^j_{kn}-T^k_{nm}\,T^j_{sk})    -  \omega^{is}(T^j_{sm;n}-T^j_{sn;m}) \cr
&& +\,\omega^{sj}(T^k_{sn}\,T^i_{km}-T^k_{sm}\,T^i_{kn}-T^k_{nm}\,T^i_{sk}) 
-\omega^{sj}(T^i_{sm;n} -T^i_{sn;m} ) \ ,
\end{eqnarray*}
which we rearrange to give
\begin{eqnarray*}
0 &=& \, \omega^{is}\big((T^k_{ns}\,T^j_{mk}+T^k_{sm}\,T^j_{nk}+T^k_{mn}\,T^j_{sk})    - (T^j_{sm;n}+T^j_{ns;m})
 +  R^j_{\phantom{i}smn} \big)  \cr
&& +\,\omega^{sj}\big((T^k_{ns}\,T^i_{mk}+T^k_{sm}\,T^i_{nk}+T^k_{mn}\,T^i_{sk}) 
- (T^i_{sm;n} + T^i_{ns;m} ) + R^i_{\phantom{i}smn}\big)\ .
\end{eqnarray*}
Using (B1) gives the symmetry of $H^{ij}$. 
\begin{eqnarray*}
0 &=& \, \omega^{is}( T^j_{mn;s}-R^j{}_{mns}-R^j{}_{nsm} )   +   \omega^{sj}( T^i_{mn;s}-R^i{}_{mns}-R^i{}_{nsm} )\ .\quad \square
\end{eqnarray*}

\begin{lemma}   \label{bchisaokjhgc}
Given the compatibility condition (\ref{compatT}), the 2-forms $H^{ij}$  in Lemma~\ref{vghsuioi}
obey $ \extd H^{ij}+\Gamma^i_{rp}\,\extd x^p\wedge H^{rj} +\Gamma^j_{rp}\,\extd x^p\wedge H^{ir}=0$. 
\end{lemma}
\proof  To calculate $ \extd  H^{ij}$ it is important to note that the $i,j$ are fixed indices, and are not summed with the vector or covector basis. This is the reason for the extra Christoffel symbols entering the following expression:
\begin{eqnarray*}
\nabla_p(H^{ij})
&=&\tfrac14\,\nabla_p(\omega^{is}\big(T^j_{nm;s}-2\,R^j_{\phantom{p} nms }
\big)\,\extd x^m\wedge\extd x^n)  \cr
&=& \tfrac14\omega^{is}_{\phantom{is};p}\big(T^j_{nm;s}-2\,R^j_{\phantom{p} nms }
\big)\,\extd x^m\wedge\extd x^n \cr
&&+\, \tfrac14\omega^{is}\big(T^j_{nm;sp}-2\,R^j_{\phantom{p} nms;p }
\big)\,\extd x^m\wedge\extd x^n \cr
&& -\,\Gamma^i_{pr}\,H^{rj}-\Gamma^j_{pr}\,H^{ir}\ .
\end{eqnarray*}
Thus we have, using the compatibility condition,
\begin{eqnarray*}
&&\extd x^p\wedge \nabla_p(H^{ij}) +\Gamma^i_{rp}\,\extd x^p\wedge   H^{rj} +  \Gamma^j_{rp}\,\extd x^p\wedge   H^{ir} \cr
&=& \tfrac14\omega^{is}_{\phantom{is};p}\big(T^j_{nm;s}-2\,R^j_{\phantom{p} nms }
\big)\,\extd x^p\wedge \extd x^m\wedge\extd x^n \cr
&&+\, \tfrac14\omega^{is}\big(T^j_{nm;sp}-2\,R^j_{\phantom{p} nms;p }
\big)\,\extd x^p\wedge \extd x^m\wedge\extd x^n \cr
&& -\,T^i_{pr}\,\extd x^p\wedge H^{rj}-T^j_{pr}\,\extd x^p\wedge H^{ir}\cr
&=& -\,\tfrac14\,(\omega^{it}T^s_{tp}+\omega^{ts}T^i_{tp})\,\big(T^j_{nm;s}-2\,R^j_{\phantom{p} nms }
\big)\,\extd x^p\wedge \extd x^m\wedge\extd x^n \cr
&&+\, \tfrac14\omega^{is}\big(T^j_{nm;sp}-2\,R^j_{\phantom{p} nms;p }
\big)\,\extd x^p\wedge \extd x^m\wedge\extd x^n \cr
&& -\,T^i_{pr}\,\extd x^p\wedge H^{rj}-T^j_{pr}\,\extd x^p\wedge H^{ir}\ .
\end{eqnarray*}
Using (B1) and then differentiating, we see that
\begin{eqnarray*} \label{bianchitorbchals}
 \sum_{\mathrm{cyclic}(pmn)} \big(T^j_{nm;p}-R^j{}_{pnm}-T^j_{pr}\,T^r_{nm}\big) &=& 0 \ ,\cr
  \sum_{\mathrm{cyclic}(pmn)} \big(T^j_{nm;ps}-R^j{}_{pnm;s}-T^j_{pr;s}\,T^r_{nm}
  -T^j_{pr}\,T^r_{nm;s}  \big) &=& 0 \ .
\end{eqnarray*}
Since the 3-form has cyclic symmetry in $(pmn)$,
\begin{eqnarray*}
&&\extd x^p\wedge \nabla_p(H^{ij}) +\Gamma^i_{rp}\,\extd x^p\wedge   H^{rj} +  \Gamma^j_{rp}\,\extd x^p\wedge   H^{ir} \cr
&=& -\,\tfrac14\,(\omega^{it}T^s_{tp}+\omega^{ts}T^i_{tp})\,\big(T^j_{nm;s}-2\,R^j_{\phantom{p} nms }
\big)\,\extd x^p\wedge \extd x^m\wedge\extd x^n \cr
&&+\, \tfrac14\omega^{is}\big(T^j_{nm;sp}
- T^j_{nm;ps} + R^j{}_{pnm;s} + T^j_{pr;s}\,T^r_{nm}  +T^j_{pr}\,T^r_{nm;s} 
-2\,R^j_{\phantom{p} nms;p }
\big)\,\extd x^p\wedge \extd x^m\wedge\extd x^n \cr
&& -\,T^i_{pr}\,\extd x^p\wedge H^{rj}-T^j_{pr}\,\extd x^p\wedge H^{ir}   \cr
&=& -\,\tfrac14\,(\omega^{it}T^s_{tp}+\omega^{ts}T^i_{tp})\,\big(T^j_{nm;s}-2\,R^j_{\phantom{p} nms }
\big)\,\extd x^p\wedge \extd x^m\wedge\extd x^n \cr
&&+\, \tfrac14\omega^{is}\big(T^r_{nm}\,R^{j}{}_{rps}
- T^j_{rm} \,R^{r}{}_{nps} - T^j_{nr} \,R^{r}{}_{mps} - T^j_{nm;r} \,T^{r}_{ps}
+ R^j{}_{pnm;s} \cr
&&  + \  T^j_{pr;s}\,T^r_{nm}  +T^j_{pr}\,T^r_{nm;s} 
-2\,R^j_{\phantom{p} nms;p }
\big)\,\extd x^p\wedge \extd x^m\wedge\extd x^n \cr
&& -\,T^i_{pr}\,\extd x^p\wedge H^{rj}-T^j_{pr}\,\extd x^p\wedge H^{ir}   \cr
&=& -\,\tfrac14\,\big(  \omega^{is}T^r_{sp}\,\(T^j_{nm;r}-2\,R^j_{\phantom{p} nmr })+\omega^{rs}T^i_{rp} \,\(T^j_{nm;s}-2\,R^j_{\phantom{p} nms })
\big)\,\extd x^p\wedge \extd x^m\wedge\extd x^n \cr
&&+\, \tfrac14\omega^{is}\big(T^r_{nm}\,R^{j}{}_{rps}
- T^j_{rm} \,R^{r}{}_{nps} - T^j_{nr} \,R^{r}{}_{mps} - T^j_{nm;r} \,T^{r}_{ps}
+ R^j{}_{pnm;s}  \cr
&&   +   \  T^j_{pr;s}\,T^r_{nm}  +T^j_{pr}\,T^r_{nm;s} 
-2\,R^j_{\phantom{p} nms;p }
\big)\,\extd x^p\wedge \extd x^m\wedge\extd x^n \cr
&& -\, \tfrac14\omega^{rs}( T^j_{nm;s}-2 R^j{}_{nms})\,T^i_{pr}\,\extd x^p\wedge \extd x^m\wedge\extd x^n  -  \tfrac14\omega^{is} \, T^j_{pr}\, \left( T^r_{nm;s}-2 R^r{}_{nms}\right)\, \extd x^p\wedge \extd x^m\wedge\extd x^n  \cr
&=& -\,\tfrac14\,\big(  \omega^{is}T^r_{sp}\,\(-2\,R^j_{\phantom{p} nmr })
\big)\,\extd x^p\wedge \extd x^m\wedge\extd x^n \cr
&&+\, \tfrac14\omega^{is}\big(T^r_{nm}\,R^{j}{}_{rps}
- T^j_{rm} \,R^{r}{}_{nps} - T^j_{nr} \,R^{r}{}_{mps} 
+ R^j{}_{pnm;s} + T^j_{pr;s}\,T^r_{nm}  
-2\,R^j_{\phantom{p} nms;p }
\big)\,\extd x^p\wedge \extd x^m\wedge\extd x^n \cr
&&  - \,  \tfrac14\omega^{is} \, T^j_{pr}\, \left( -2 R^r{}_{nms}\right)\, \extd x^p\wedge \extd x^m\wedge\extd x^n  \cr
&=& \tfrac14\omega^{is}\big(T^r_{nm}\,R^{j}{}_{rps}
- T^j_{rm} \,R^{r}{}_{nps} - T^j_{nr} \,R^{r}{}_{mps} 
+ R^j{}_{pnm;s} + T^j_{pr;s}\,T^r_{nm}  
-2\,R^j_{\phantom{p} nms;p }   \cr
&&  + \,  2\,T^r_{sp}\,R^j_{\phantom{p} nmr }  +2\,T^j_{pr}\,  R^r{}_{nms}
\big)\,\extd x^p\wedge \extd x^m\wedge\extd x^n \ .
\end{eqnarray*}
Given the overall $\extd x^p\wedge \extd x^m\wedge\extd x^n$ factor, we can make the following substitutions:
\newline\noindent $-\,R^j{}_{ nms;p} \mapsto -\,R^j{}_{ pns;m} \mapsto R^j{}_{ psn;m} 
 \mapsto -\,R^j{}_{ psm;n}   \mapsto R^j{}_{ pms;n}   $
 \newline\noindent $T^r_{sp}\,R^j{}_{nmr } \mapsto T^r_{sm}\,R^j{}_{pnr }
  \mapsto -\,T^r_{ms}\,R^j{}_{pnr }  \mapsto T^r_{ns}\,R^j{}_{pmr }  \mapsto T^r_{ns}\,R^j{}_{pmr }
   \mapsto -\, T^r_{sn}\,R^j{}_{pmr }$
 \newline\noindent $T^j_{pr}\,  R^r{}_{nms} \mapsto T^j_{nr}\,  R^r{}_{mps} \mapsto -\, T^j_{mr}\,  R^r{}_{nps} \mapsto T^j_{rm}\,  R^r{}_{nps}$
   \newline\noindent Using these we can rewrite the previous equations, and then use (B2) to get
\begin{eqnarray*}
&&\extd x^p\wedge \nabla_p(H^{ij}) +\Gamma^i_{rp}\,\extd x^p\wedge   H^{rj} +  \Gamma^j_{rp}\,\extd x^p\wedge   H^{ir} \cr
&=& \tfrac14\omega^{is}\big(T^r_{nm}\,R^{j}{}_{rps}
- T^j_{rm} \,R^{r}{}_{nps} - T^j_{nr} \,R^{r}{}_{mps} 
 + T^j_{pr;s}\,T^r_{nm}  \cr
&& +  \,  R^j{}_{pnm;s} + R^j{}_{ psn;m} + R^j{}_{ pms;n} - T^r_{ms}\,R^j{}_{pnr } - T^r_{sn}\,R^j{}_{pmr }
  \cr
&&  + \,  T^j_{nr}\,  R^r{}_{mps} +  T^j_{rm}\,  R^r{}_{nps}  
\big)\,\extd x^p\wedge \extd x^m\wedge\extd x^n \cr
&=& \tfrac14\omega^{is}\big(T^r_{nm}\,R^{j}{}_{rps}
 + T^j_{pr;s}\,T^r_{nm}   +   T^r_{nm}\,R^j{}_{psr }
\big)\,\extd x^p\wedge \extd x^m\wedge\extd x^n \cr
&=& \tfrac14\,\omega^{is} \, T^r_{nm}\, \big(R^{j}{}_{rps}
 + T^j_{pr;s}   +   R^j{}_{psr }
\big)\,\extd x^p\wedge \extd x^m\wedge\extd x^n \ .
\end{eqnarray*}
We also need to calculate
\begin{eqnarray*}
&&\tfrac12\,T^t_{vu}\,\extd x^v\wedge\extd x^u\wedge (\partial_t\, \righthalfcup\,H^{ij}) \cr
&=& \tfrac12\,T^t_{vu}\,\extd x^v\wedge\extd x^u\wedge (\partial_t\, \righthalfcup\,(\tfrac14\omega^{is}\big(T^j_{nm;s}-2\,R^j_{\phantom{p} nms }
\big)\,\extd x^m\wedge\extd x^n))  \cr
&=& \tfrac18\,T^p_{vu}\,\omega^{is}\big(T^j_{np;s}-2\,R^j_{\phantom{p} nps }
 \big)\,\extd x^v\wedge\extd x^u\wedge\extd x^n  \cr
&&  -\,  \tfrac18\,T^p_{vu}\,\omega^{is}\big(T^j_{pn;s}-2\,R^j_{\phantom{p} pns }
 \big)\,\extd x^v\wedge\extd x^u\wedge \extd x^n  \cr
&=& \tfrac14\,T^p_{vu}\,\omega^{is}\big(T^j_{np;s}-R^j_{\phantom{p} nps } + \,R^j_{\phantom{p} pns }
  \big)\,\extd x^v\wedge\extd x^u\wedge\extd x^n   \cr
&=& \tfrac14\,T^r_{pm}\,\omega^{is}\big(T^j_{nr;s}-R^j_{\phantom{p} nrs } + \,R^j_{\phantom{p} rns }
\big)\,\extd x^p\wedge\extd x^m\wedge\extd x^n    \cr
&=& \tfrac14\,T^r_{mn}\,\omega^{is}\big(T^j_{pr;s}-R^j_{\phantom{p} prs } + \,R^j_{\phantom{p} rps }
\big)\,\extd x^p\wedge\extd x^m\wedge\extd x^n     \cr
&=& \tfrac14\,T^r_{mn}\,\omega^{is}\big(T^j_{pr;s}+R^j_{\phantom{p} psr } + \,R^j_{\phantom{p} rps }
\big)\,\extd x^p\wedge\extd x^m\wedge\extd x^n \ .
\end{eqnarray*}
The result follows by using
\begin{eqnarray*}
 \extd  H^{ij} &=& \extd x^k\wedge\nabla_k H^{ij}+ \tfrac12\,T^t_{kn}\,\extd x^k\wedge\extd x^n\wedge (\partial_t\, \righthalfcup\,H^{ij})\ .\quad\square
\end{eqnarray*}

\medskip\noindent The last two lemmas prove the following theorem:

\begin{theorem} \label{dga} Suppose that the compatibility condition (\ref{compatT}) holds.
Then the conditions on $H^{ij}$ in Proposition~\ref{wedge1leib} hold, i.e.\ we have a differential graded algebra $(\wedge_1,\extd)$ to order $\lambda$.
\end{theorem}

\begin{proposition} Over $\C$, the DGA above is a $*$-DGA.
\end{proposition}
\proof As both $\extd$ and $\star$ are undeformed, it is automatic that $\extd(\xi^*)=(\extd \xi)^*$. Next
\begin{eqnarray*}
\eta^*\wedge_1\xi^* &=& \eta^*\wedge\xi^*+\tfrac\lambda2\,\omega^{ij}\,\nabla_i\eta^*\wedge\nabla_j\xi^*+\lambda\,(-1)^{|\eta|+1}\, H^{ij}\wedge (\partial_i \, \righthalfcup\, \eta^*)\wedge(\partial_j \, \righthalfcup\, \xi^*)   \cr
&=& \eta^*\wedge\xi^*+\tfrac\lambda2\,\omega^{ij}\,\nabla_i\eta^*\wedge\nabla_j\xi^*+\lambda\,(-1)^{|\eta|+1}\, H^{ij}\wedge (\partial_i \, \righthalfcup\, \eta)^*\wedge(\partial_j \, \righthalfcup\, \xi)^*   \cr
&=& (-1)^{|\xi|\,|\eta|}  (\xi\wedge\eta)^*
+(-1)^{|\xi|\,|\eta|} \tfrac\lambda2\,\omega^{ij}\,(\nabla_j\xi\wedge\nabla_i\eta)^* \cr
&& +  \, \lambda\,(-1)^{|\eta|+1+(|\eta|-1)(|\xi|-1)}\, H^{ij}\wedge ((\partial_j \, \righthalfcup\, \xi)\wedge (\partial_i \, \righthalfcup\, \eta))^*   \cr
&=& (-1)^{|\xi|\,|\eta|}  (\xi\wedge\eta)^*
+(-1)^{|\xi|\,|\eta|} (\tfrac\lambda2\,\omega^{ij}\,\nabla_i\xi\wedge\nabla_j\eta)^* \cr
&& +  \, \lambda\,(-1)^{|\xi|\,|\eta|+|\xi|+2}\, (H^{ij}\wedge (\partial_j \, \righthalfcup\, \xi)\wedge (\partial_i \, \righthalfcup\, \eta))^*   \cr
&=& (-1)^{|\xi|\,|\eta|}  (\xi\wedge\eta)^*
+(-1)^{|\xi|\,|\eta|} (\tfrac\lambda2\,\omega^{ij}\,\nabla_i\xi\wedge\nabla_j\eta)^* \cr
&& +  \, (-1)^{|\xi|\,|\eta|}\, ( \lambda\,(-1)^{|\xi|+1}\,H^{ij}\wedge (\partial_j \, \righthalfcup\, \xi)\wedge (\partial_i \, \righthalfcup\, \eta))^*   \quad  \square
\end{eqnarray*}

\subsection{The quantum torsion of the quantising connection}
Here we consider the quantum torsion of the quantum connection given by applying 
Theorem~\ref{functor} to the Poisson-compatible connection $(\Omega^1,\nabla)$ itself. This is intimately tied
up with the quantum differential calculus above. 

\begin{lemma}
We have
\begin{eqnarray*}
\wedge_1( \sigma_{Q}+\id^{\tens 2})(\eta\tens_1\xi) 
 &=& \tfrac12 \, \lambda \,  \xi_{j}\,\eta_{p} \,
\omega^{ji}\, T^p_{nk;i}\,\extd x^k\wedge\extd x^n\ ,
\end{eqnarray*}
so the quantum torsion is a right module map if and only if $\omega^{ji}\, T^p_{nk;i}=0$. 
\end{lemma}
\proof  From Theorem~\ref{functor},
\begin{eqnarray*}
\wedge_1( \sigma_{Q}+\id^{\tens 2})(\eta\tens_1\xi)\, &=&  \xi\wedge_1 \eta +\eta\wedge_1\xi  + \lambda\, \omega^{ij}\,\nabla_j\xi\wedge_1\nabla_{i}\eta \cr
&&  + \, \lambda\, \omega^{ij}\,\xi_{j}\,\extd x^k \wedge_1 [\nabla_{k},\nabla_{i}]\eta  \ ,
\end{eqnarray*}
and using Proposition~\ref{wedge1leib} and the definition of $H^{ij}$ in Lemma~\ref{vghsuioi} we get
\begin{eqnarray*}
&&\wedge_1( \sigma_{Q}+\id^{\tens 2})(\eta\tens_1\xi) \cr
 &=&  2 \, \lambda \, H^{ij} \, \xi_{i}\,\eta_{j} +  \lambda\, \omega^{ij}\,\xi_{j}\,\extd x^k \wedge_1 [\nabla_{k},\nabla_{i}]\eta  \cr
&=& \tfrac12 \, \lambda \,  \xi_{i}\,\eta_{j} \,
\omega^{is}\left( T^j_{nm;s}-2\, R^j{}_{nms}\right)\extd x^m\wedge\extd x^n
 -  \lambda\, \omega^{ij}\,\xi_{j}\,R^p_{nki}\,\eta_p\,\extd x^k \wedge_1 \extd x^n   \cr
 &=& \tfrac12 \, \lambda \,  \xi_{j}\,\eta_{p} \,
\omega^{ji}\left( T^p_{nk;i}-2\, R^p{}_{nki}\right)\extd x^k\wedge\extd x^n
 -  \lambda\, \omega^{ij}\,\xi_{j}\,R^p_{nki}\,\eta_p\,\extd x^k \wedge \extd x^n   \ .\quad\square
\end{eqnarray*}

\begin{proposition}  \label{torsionQ}
The quantum torsion of the quantising connection on $\Omega^1(M)$ is
\begin{eqnarray*}
(\wedge_1\nabla_Q-\extd)(\xi)=(\wedge\nabla-\extd)(\xi)
+ \tfrac\lambda4\,(\partial_j  \righthalfcup  \nabla_i\xi)\,
\omega^{is}\, T^j_{nm;s}\,\extd x^m\wedge\extd x^n\ .
\end{eqnarray*}
\end{proposition}
\proof  Here all covariant derivatives are the quantising connection on $\Omega^1 (M)$:
\begin{eqnarray*}
\wedge_1\,\nabla_{Q}\xi &=& \wedge_1\,q^{-1}\nabla\xi-\tfrac{\lambda}{2} \,\omega^{ij}\,\extd x^k\wedge[\nabla_{k},\nabla_{j}]\nabla_{i}\xi \cr
&=& \wedge_1\,q^{-1}(\extd x^k\tens_0\nabla_k\xi)  -  \tfrac{\lambda}{2} \,\omega^{ij}\,\extd x^k\wedge[\nabla_{k},\nabla_{j}]\nabla_{i}\xi \cr
&=& \extd x^k\wedge\nabla_k\xi + \lambda\,H^{ij}\, (\partial_j  \righthalfcup  \nabla_i\xi)    -  \tfrac{\lambda}{2} \,\omega^{is}\,\extd x^m\wedge[\nabla_{m},\nabla_{s}]\nabla_{i}\xi \ ,
\end{eqnarray*}
and using the definition of $H^{ij}$ in Lemma~\ref{vghsuioi} we get
\begin{eqnarray*}
\wedge_1\,\nabla_{Q}\xi 
&=& \extd x^k\wedge\nabla_k\xi + \tfrac\lambda4\,(\partial_j  \righthalfcup  \nabla_i\xi)\,
\omega^{is}\left( T^j_{nm;s}-2 R^j{}_{nms}\right)\extd x^m\wedge\extd x^n \cr
&&   +  \,  \tfrac{\lambda}{2} \,\omega^{is}\, R^j{}_{nms}\,   \extd x^m\wedge\extd x^n (\partial_j  \righthalfcup  \nabla_i\xi)\ .\quad\square
\end{eqnarray*}

Finally,  in classical differential geometry one has for any linear connection that
$\nabla_k(\xi\wedge\eta)=\nabla_k(\xi)\wedge\eta+\xi\wedge\nabla_k(\eta)$. This works because the usual tensor product covariant derivative on $\Omega^1(M)\tens_0\Omega^1(M)$ preserves symmetry, so things in the kernel of $\wedge$ stay in the kernel. So given a quantising covariant derivative $\nabla$ on $\Omega^1 (M)$, we naturally get covariant derivatives on all the $\Omega^i(M)$, which we also call $\nabla$. We can then say that the wedge product  $\wedge_0:\Omega^n(M)\tens_0\Omega^n(M)\to \Omega^{n+m}(M)$ intertwines the covariant derivatives. 
We conclude by studying what happens when we quantise the covariant derivatives.

\begin{proposition}   \label{ouyccdyuk}
\begin{eqnarray*}
&&\kern-30pt (\id\tens_1 \wedge_1)\nabla_{Q\tens_1 Q} (\xi\tens_1\eta)  - \nabla_{Q}(\xi\, \wedge_1\, \eta) \cr
&=&   \lambda\,(-1)^{|\xi|}\, \extd x^k\tens_1 (\nabla_k  H^{ij}+\Gamma_{kp}^i\,H^{pj}+\Gamma_{kp}^j\,H^{ip})\wedge (\partial_i \, \righthalfcup\, \xi)\wedge(\partial_j \, \righthalfcup\, \eta)   \ .
\end{eqnarray*}
\end{proposition}
\noindent \textbf{Proof:}\quad From modifying (\ref{vcusvdsavbfea}) we get
\begin{eqnarray}\label{vcusvdvbbdfs}
\xymatrix{
Q(\Omega(M))\tens_1 Q(\Omega(M)) \ar[r]^{q} \ar[d]_{\nabla_{Q(\Omega(M))\tens_1 Q(\Omega(M))}} 
&Q(\Omega(M)\tens_0 \Omega(M))\ar[d]_{\nabla_{Q(\Omega(M)\tens_0 \Omega(M))}}      \\
Q(\Omega(M))\tens_1 Q(\Omega(M))\tens_1 Q(\Omega(M)) \ar[r]_{\id\tens q} &
Q(\Omega(M)) \tens_1  Q(\Omega(M)\tens_0 \Omega(M) ) \ar[d]_{\id\tens_1(\wedge)} \\
& Q(\Omega(M)) \tens_1  Q(\Omega(M)) 
   } 
\end{eqnarray}
As classically $\wedge$ intertwines the covariant derivatives,
\begin{eqnarray*}
(\id\tens_1 (\wedge q))\nabla_{Q(\Omega(M))\tens_1 Q(\Omega(M))}&=& \nabla_{Q(\Omega(M))}\, (\wedge q):\\ &&Q(\Omega(M))\tens_1 Q(\Omega(M)) \to  Q(\Omega(M)) \tens_1  Q(\Omega(M)) \ .
\end{eqnarray*}
In the notation of Proposition~\ref{wedge1leib} we now look at $\xi\wedge_1\eta=\xi\wedge_Q\eta+\lambda\, \xi \, \widehat\wedge \, \eta$ where
\[ \xi \, \widehat{\wedge} \, \eta=(-1)^{|\xi|+1}\, H^{ij}\wedge (\partial_i \, \righthalfcup\, \xi)\wedge(\partial_j \, \righthalfcup\, \eta)\   .\]
Then we get
\begin{eqnarray*}
 && \lambda\,(\id\tens_1 \widehat{\wedge})\nabla_{Q(\Omega(M))\tens_1 Q(\Omega(M))} (\xi\tens_1\eta)   \cr
 &=& \lambda\,(\id\tens_1 \widehat{\wedge})(\extd x^k\tens_1 (\nabla_k\xi\tens_1\eta + \xi\tens_1\nabla_k\eta))  \cr
  &=& \lambda\,\extd x^k\tens_1 (\nabla_k\xi \, \widehat\wedge \, \eta + \xi \, \widehat\wedge \, \nabla_k\eta)   \cr
    &=& \lambda\,(-1)^{|\xi|+1}\, \extd x^k\tens_1 H^{ij}\wedge (
     (\partial_i \, \righthalfcup\, \nabla_k\xi)\wedge(\partial_j \, \righthalfcup\, \eta) 
     +  (\partial_i \, \righthalfcup\, \xi)\wedge(\partial_j \, \righthalfcup\, \nabla_k\eta)  ) \ .
\end{eqnarray*}
Also, using $\nabla_i(v\, \righthalfcup\,\xi)=\nabla_i(v)\, \righthalfcup\,\xi+v\, \righthalfcup\,\nabla_i\xi$
\begin{eqnarray*}
&& \lambda\,\nabla_{Q(\Omega(M))}(\xi\, \widehat{\wedge}\, \eta) \cr
&=& \lambda\,(-1)^{|\xi|+1}\, \extd x^k\tens_1 \nabla_k(  H^{ij}\wedge (\partial_i \, \righthalfcup\, \xi)\wedge(\partial_j \, \righthalfcup\, \eta) )  \cr
    &=& \lambda\,(-1)^{|\xi|+1}\, \extd x^k\tens_1 H^{ij}\wedge (
     (\partial_i \, \righthalfcup\, \nabla_k\xi)\wedge(\partial_j \, \righthalfcup\, \eta) 
     +  (\partial_i \, \righthalfcup\, \xi)\wedge(\partial_j \, \righthalfcup\, \nabla_k\eta)  )  \cr
&&+\,  \lambda\,(-1)^{|\xi|+1}\, \extd x^k\tens_1 \nabla_k  H^{ij}\wedge (\partial_i \, \righthalfcup\, \xi)\wedge(\partial_j \, \righthalfcup\, \eta)   \cr
&&+\,  \lambda\,(-1)^{|\xi|+1}\, \extd x^k\tens_1 H^{ij}\wedge (
     (\Gamma_{ki}^p\partial_p \, \righthalfcup\, \xi)\wedge(\partial_j \, \righthalfcup\, \eta) 
     +  (\partial_i \, \righthalfcup\, \xi)\wedge(\Gamma_{kj}^p\,\partial_p \, \righthalfcup\, \eta)  )  \cr
&=& \lambda\,(\id\tens_1 \widehat{\wedge})\nabla_{Q(\Omega(M))\tens_1 Q(\Omega(M))} (\xi\tens_1\eta)   \cr
&&+\,  \lambda\,(-1)^{|\xi|+1}\, \extd x^k\tens_1 (\nabla_k  H^{ij}+\Gamma_{kp}^i\,H^{pj}+\Gamma_{kp}^j\,H^{ip})\wedge (\partial_i \, \righthalfcup\, \xi)\wedge(\partial_j \, \righthalfcup\, \eta)   \ .   \quad \square
\end{eqnarray*}

\subsection{Quantizing other linear connections relative to the background $(\Omega^1,\nabla)$}

Here we extend the above to other connections $\nabla_S=\nabla+S$ on $\Omega^1(M)$ different from the quantizing one $\nabla$, where $S(\xi)=\xi_p\,S^p_{nm}\,\extd x^n\tens \extd x^m$ for $\xi\in \Omega^1(M)$. Quantisation is achieved on the same quantum bundle as defined by $\nabla$, using Proposition~\ref{quantQS}

\begin{proposition}  \label{bhadsuiytcddxtuy}
The torsion of $\nabla_{QS}$ is given by
\begin{eqnarray*}
T_{\nabla_{QS}}(\xi) &=& T_{\nabla_S}(\xi)+ \tfrac\lambda4\, \xi_{p;i}\,
\omega^{ij}\, (T^p_{nm;j}
-2\,S^p_{nm;j})\,\extd x^m\wedge \extd x^n \cr
&&+\,\lambda\, \xi_p\,(S^p_{nm}\,H^{nm} +\, \tfrac12\,\omega^{ij}\,S^p_{nm;\hat\jmath i}\,\extd x^n\wedge \extd x^m) \ .
\end{eqnarray*}
Note that the hat on $\hat\jmath$ denotes that the $j$ index does not take part in the covariant differentiation in the $i$ direction. 
\end{proposition}
\proof  We have
\begin{eqnarray*}
T_{\nabla_{QS}}(\xi) = T_{\nabla_{Q}}(\xi)+\wedge_1\,q^{-1}S(\xi)+\wedge\,\tfrac\lambda2\,\omega^{ij}\,\nabla_i\circ\nabla_j(S)(\xi)\ .
\end{eqnarray*}
Now
\begin{eqnarray*}
\wedge_1\,q^{-1}S(\xi) &=& \xi_p\,S^p_{nm}\,\extd x^n\wedge \extd x^m
+\lambda\, \xi_p\,S^p_{nm}\,H^{nm} \ ,\cr
\wedge\,\tfrac\lambda2\,\omega^{ij}\,\nabla_i\circ\nabla_j(S)(\xi) &=& 
\wedge\,\tfrac\lambda2\,\omega^{ij}\,\nabla_i(\xi_p\,S^p_{nm;j}\,\extd x^n\tens \extd x^m) \cr
 &=& 
\tfrac\lambda2\,\omega^{ij}\,\nabla_i(\xi_p\,S^p_{nm;j}\,\extd x^n\wedge \extd x^m)\cr
 &=& 
\tfrac\lambda2\,\omega^{ij}\,(\xi_{p;i}\,S^p_{nm;j}+\xi_p\,S^p_{nm;\hat\jmath i})\,\extd x^n\wedge \extd x^m\ .
\end{eqnarray*}
By Proposition~\ref{torsionQ} we get
\begin{eqnarray*}
T_{\nabla_{QS}}(\xi) &=& T_{\nabla}(\xi)+ \tfrac\lambda4\, \xi_{j;i}\,
\omega^{is}\, T^j_{nm;s}\,\extd x^m\wedge\extd x^n+\wedge S(\xi)+\lambda\, \xi_p\,S^p_{nm}\,H^{nm} \cr
&&+\, \tfrac\lambda2\,\omega^{ij}\,(\xi_{p;i}\,S^p_{nm;j}+\xi_p\,S^p_{nm;\hat\jmath i})\,\extd x^n\wedge \extd x^m \cr
&=& T_{\nabla_S}(\xi)+ \tfrac\lambda4\, \xi_{p;i}\,
\omega^{ij}\, T^p_{nm;j}\,\extd x^m\wedge\extd x^n+\lambda\, \xi_p\,S^p_{nm}\,H^{nm} \cr
&&+\, \tfrac\lambda2\,\omega^{ij}\,(\xi_{p;i}\,S^p_{nm;j}+\xi_p\,S^p_{nm;\hat\jmath i})\,\extd x^n\wedge \extd x^m \cr
&=& T_{\nabla_S}(\xi)+ (\tfrac\lambda4\, \xi_{p;i}\,
\omega^{ij}\, T^p_{nm;j}
-\, \tfrac\lambda2\,\omega^{ij}\,\xi_{p;i}\,S^p_{nm;j})\,\extd x^m\wedge \extd x^n \cr
&&+\,\lambda\, \xi_p\,S^p_{nm}\,H^{nm} +\, \tfrac\lambda2\,\omega^{ij}\,\xi_p\,S^p_{nm;\hat\jmath i}\,\extd x^n\wedge \extd x^m \cr
&=& T_{\nabla_S}(\xi)+ \tfrac\lambda4\, \xi_{p;i}\,
\omega^{ij}\, (T^p_{nm;j}
-2\,S^p_{nm;j})\,\extd x^m\wedge \extd x^n \cr
&&+\,\lambda\, \xi_p\,(S^p_{nm}\,H^{nm} +\, \tfrac12\,\omega^{ij}\,S^p_{nm;\hat\jmath i}\,\extd x^n\wedge \extd x^m) \ .\quad\square
\end{eqnarray*}

\begin{corollary} \label{kjgccghgc}
In the case where $\nabla_S$ is torsion free,  the quantum torsion 
$T_{\nabla_{QS}}(\xi):= \tfrac\lambda2\, \xi_p\,A^p_{nm}\, \extd x^m\wedge \extd x^n$ has associated tensor
\begin{eqnarray*}
A^p_{nm} &=& \tfrac14\omega^{is}\,(S^p_{ij}+S^p_{ji})\,( T^j_{nm;s}-R^j{}_{nms}+R^j{}_{mns})  - \tfrac14\,\omega^{ij}\,(T^s_{nm}\,R^p{}_{sij}-T^p_{sm}\,R^s{}_{nij}+T^p_{sn}\,R^s{}_{mij})\ .
\end{eqnarray*}
\end{corollary}
\proof  
Begin with
\begin{eqnarray*}
\wedge\,\nabla_S(\extd x^p)=\wedge\,\nabla(\extd x^p)+S^p_{nm}\,\extd x^n\wedge \extd x^m
\end{eqnarray*}
so $0=T_\nabla(\extd x^p)+S^p_{nm}\,\extd x^n\wedge \extd x^m$ and from (\ref{vuyihgfxxy}) we deduce
\begin{eqnarray} \label{bhcuiiycxyr}
S^p_{nm}\,\extd x^n\wedge \extd x^m=\tfrac12\,T^p_{nm}\,\extd x^n\wedge \extd x^m\ .
\end{eqnarray}
Then Proposition~\ref{bhadsuiytcddxtuy} gives
\begin{eqnarray*}
T_{\nabla_{QS}}(\xi) &=& T_{\nabla_S}(\xi) +
\tfrac\lambda2\, \xi_p\,(2\,S^p_{nm}\,H^{nm} +\, \tfrac12\,\omega^{ij}\,T^p_{nm;\hat\jmath i}\,\extd x^n\wedge \extd x^m) \ ,
\end{eqnarray*}
and use the formula for the curvature of a tensor and the symmetry of $H^{nm}$.\quad$\square$

We see that quantisation introduces an element of torsion at order $\lambda$ in the quantisation $\nabla_{QS}$. We similarly
 look at Lemma~\ref{oiuyfvcftyytrs} to measure the deviation of $\nabla_{QS}$ from being star preserving and find an
 error of order $\lambda$:

\begin{lemma}   \label{bgyiuouy1}
Over $\C$ and if $S$ is real, the difference $D^a_{ijnm}\overline{\tfrac\lambda2\,\omega^{ij}\, \extd x^n\tens\extd x^m}$ in going clockwise minus anticlockwise round the diagram in 
Lemma~\ref{oiuyfvcftyytrs} starting from $Q(\extd x^a)$ is given by
\begin{eqnarray*}
D^a_{ijnm}=2\, S^a_{ip}\,S^p_{nm;j}-(S^b_{nm}\,R^a{}_{bij}
-S^a_{rm}\,R^r{}_{nij}   - S^a_{nr}\,R^r{}_{mij}   ) -2\,S^a_{jr}\,R^r{}_{mni} 
\end{eqnarray*}
\end{lemma}
\proof  Putting $e=e^*=\extd x^a$ in Lemma~\ref{oiuyfvcftyytrs}, and using $\nabla_j(S)(\xi)=\xi_p\,S^p_{nm;j}\,\extd x^n\tens \extd x^m$ we get
\begin{eqnarray*}
\nabla_j(S)(S_i(\extd x^a)) &=& \nabla_j(S)(S^a_{ir}\,\extd x^r) 
= S^a_{ip}\,S^p_{km;j}\,\extd x^k\tens \extd x^m\ ,\cr
\extd x^k\tens [\nabla_k,\nabla_i]S_j(\extd x^a)&=& \extd x^k\tens [\nabla_k,\nabla_i](S^a_{jr}\,\extd x^r)=
-\,S^a_{jr}\,R^r{}_{mki}\, \extd x^k\tens\extd x^m\ .
\end{eqnarray*}
Now we use the antisymmetry of $\omega^{ij}$ to get
\begin{eqnarray*}
\omega^{ij}\,\nabla_i(\nabla_j(S))(\extd x^a) &=& \tfrac12\,\omega^{ij}\, (S^b_{km}\,R^a{}_{bij}
-S^a_{rm}\,R^r{}_{kij}   - S^a_{kr}\,R^r{}_{mij}   )\extd x^k\tens \extd x^m\ .
\end{eqnarray*}

We note that this is equivalent to  the derivative of $\star:Q(\Omega^1(M))\to Q(\overline{\Omega^1(M)})$ being
\begin{eqnarray*}
\nabla_{QS}(\star)(Q(\extd x^a))=\tfrac\lambda2\,\omega^{ij} D^a_{ijkm} \, \extd x^k\tens\overline{\extd x^m}
\end{eqnarray*}
and we see that this is not necessarily zero. \endproof

Next we consider adding a correction, so 
\[ \nabla_1=\nabla_{QS}+\lambda K\]
where $K:\Omega^1(M) \to \Omega^1(M) \tens_0 \Omega^1(M)$ is given by
$K(\xi)=\xi_p\,K^p_{nm}\,\extd x^n\tens \extd x^m$. 

\begin{theorem}\label{starpresK}
Over $\C$ and for any real $S$, there is a unique real $K$ such that $\nabla_{QS}+\lambda K$ is star preserving (namely
$K^a_{nm} = \tfrac14\,\omega^{ij} D^a_{ijnm}$). Moreover, if $\nabla_S$ is torsion free, this unique $\nabla_{QS}+\lambda K$ is quantum torsion free.
\end{theorem}
\proof  We look at the following diagram:
\begin{eqnarray*}
\xymatrix{
\overline{Q(\Omega^1(M))} = Q(\overline{\Omega^1(M)})  \ar[d]_{\overline{\lambda\, K } } 
& Q(\Omega^1(M))  \ar[l]_{  \star  }    \ar[r]^{ \lambda\, K }  &   
Q(\Omega^1(M)) \tens_1 Q(\Omega^1(M))   \ar[d]_{\star\tens_1\star}   \\
\overline{ Q(\Omega^1(M)) \tens_1 Q(\Omega^1(M))}&
\overline{ Q(\Omega^1(M)) \tens_1 Q(\Omega^1(M))}   \ar[l]_{  \overline{\sigma_{QS}}  }  &
\overline{ Q(\Omega^1(M)) } \tens_1 \overline{ Q(\Omega^1(M)) }   \ar[l]_{\Upsilon^{-1}}  
}
\end{eqnarray*}
where at this order $\sigma_{QS}$ is simply transposition. Hence for $\nabla_{QS}+\lambda K$ the effect of adding $K$ is to add 
\begin{eqnarray*} 
-\,\overline{\lambda\,(K^a_{nm}+(K^a_{nm})^*)\,\extd x^n\tens\extd x^m}\ .
\end{eqnarray*}
to the difference in Lemma~\ref{bgyiuouy1}. This gives the unique value if we assume $K$ is real for the connection to be $*$-preserving. 
Adding $K$  also adds $\lambda\,\xi_a\,K^a_{nm}\,\extd x^n\wedge \extd x^m$ to the formula for the torsion in Proposition~\ref{bhadsuiytcddxtuy} so if
$K$ has the unique real value stated and if $\nabla_S$ is torsion free, and using (\ref{bhcuiiycxyr}), 
\begin{eqnarray*}
&& K^a_{nm}\,\extd x^n\wedge\extd x^m=\tfrac14\,\omega^{ij}D^a_{ijnm}\,\extd x^n\wedge\extd x^m \cr
&=& \tfrac14\,\omega^{ij}\big( S^a_{ip}\,T^p_{nm;j}-(\tfrac12\,T^b_{nm}\,R^a{}_{bij}
-S^a_{rm}\,R^r{}_{nij}   - S^a_{nr}\,R^r{}_{mij}   ) -2\,S^a_{jr}\,R^r{}_{mni}
\big)\,\extd x^n\wedge\extd x^m \cr
&=& \tfrac14\,\omega^{ij}\big( S^a_{ip}\,T^p_{nm;j}-(\tfrac12\,T^b_{nm}\,R^a{}_{bij}
-S^a_{rm}\,R^r{}_{nij}   + S^a_{mr}\,R^r{}_{nij}   ) -2\,S^a_{jr}\,R^r{}_{mni}
\big)\,\extd x^n\wedge\extd x^m \cr
&=& \tfrac14\,\omega^{ij}\big( S^a_{ip}\,T^p_{nm;j}-(\tfrac12\,T^b_{nm}\,R^a{}_{bij}
  + T^a_{mr}\,R^r{}_{nij}   ) +2\,S^a_{ir}\,R^r{}_{mnj}
\big)\,\extd x^n\wedge\extd x^m \cr
&=& \tfrac14\,S^a_{ip}\,\omega^{ij}\big( T^p_{nm;j} +2\,R^p{}_{mnj}
\big)\,\extd x^n\wedge\extd x^m 
 - \tfrac14\,\omega^{ij}(\tfrac12\,T^b_{nm}\,R^a{}_{bij}
  + T^a_{mr}\,R^r{}_{nij}   )\,\extd x^n\wedge\extd x^m \cr
  &=& -\,S^a_{ip}\,H^{ip}- \tfrac18\,\omega^{ij}(T^b_{nm}\,R^a{}_{bij}
  + 2\,T^a_{mr}\,R^r{}_{nij}   )\,\extd x^n\wedge\extd x^m\ .
\end{eqnarray*}
Now Corollary~\ref{kjgccghgc} gives $\nabla_{QS}+\lambda K$ torsion free. \quad$\square$

Thus requiring torsion free and star-preserving gives a unique `star-preserving and torsion-preserving' quantisation of any classical torsion-free connection $\nabla_S$ on $\Omega^1(M)$. 

\section{Semiquantization of Riemannian geometry}
We are now in position to semiquantise Riemannian geometry on our above datum $(\omega,\nabla)$. We need to proceed carefully, as there are various places where modifications arise, and there are typically two  connections involved. Throughout this section suppose  $g=g_{ij}\,\extd x^i\tens\extd x^j\in\Omega^{\tens 2}(M)$ is a Riemannian metric on $M$. We start with the quantum metric and the quantisation $\nabla_Q$ of the quantizing connection $\nabla$.

\subsection{The quantised metric} \label{bchijadflkjhov} We obtain to first order a quantum metric $g_1\in \Omega^{\tens_1 2}A$ characterised by quantum symmetry and centrality. The former is the statement that $g_1$ is in the kernel of
$\wedge_1:\Omega^{\tens_1 2}A_1 \to \Omega^{2}A_1$ and the latter is that $g_1$ commutes with the elements of the algebra, i.e.\ $a.g_1=g_1.a$ for all $a\in A_1$. Without this property, we could not simply apply the metric to a number of tensor products over the algebra (i.e.\ the fiberwise tensor product of bundles), and using the metric would become much more complicated. 

As with the wedge product, we start with a functorial part of the quantum metric 
\begin{equation}\label{gQ}
g_Q:=q^{-1}_{\Omega^1,\Omega^1}(g)= g_{ij}\extd x^i\tens_1\extd x^j+\tfrac{\lambda}{2} \omega^{ij}(g_{ms,i}-g_{ks}\Gamma^k_{im})\extd x^m \tens_1\Gamma^s_{jn}\extd x^n\ .
\end{equation}
and of the quantum connection
\begin{equation}\label{nablaQ}
\nabla_Q\extd x^i=-\left( \Gamma^i_{mn} +{\lambda\over 2}\omega^{sj}(\Gamma^i_{mk,s}\Gamma^k_{jn} -\Gamma^i_{kt}\Gamma^k_{sm}\Gamma^t_{jn} -\Gamma^i_{jk} R^k{}_{nms})\right)\extd x^m \tens_1 \extd x^n
\end{equation}
by application of our functor in Section~3.

\begin{lemma} \label{lemmagQ} If we have $\nabla g=0$, then we also have $\nabla_Q g_Q=0$ as an application of Theorem~\ref{functor}. Moreover, over $\C$, $g_Q$ is `real' and $\nabla_Q$ is $*$-preserving. \end{lemma}
\proof We consider the metric as a morphism
$\tilde g:C^\infty(M)\to \Omega^1(M)\tens_0\Omega^1(M)$ in $\CD_0$, where $g=\tilde g(1)$ and $\Omega^1(M)$, is  equipped with the background quantising connection (assumed now to be metric compatible).  Then $q^{-1}_{\Omega^1,\Omega^1}Q(\tilde g):Q(C^\infty(M))\to Q(\Omega^1(M))\tens_1 Q(\Omega^1(M))$ and we evaluate this on $1$ to give the element $g_Q\in \Omega^1A_1 \tens_1\Omega^1A_1$. In this case  the morphism property of $q_{\Omega^1(M),\Omega^1(M)}$  implies 
(suppressing $M$ for clarity)
\begin{eqnarray*}
 \nabla_{Q(\Omega^1)\tens_1Q(\Omega^1)}q^{-1}_{\Omega^1,\Omega^1}\circ Q(\tilde g)(1)=(\id\tens q^{-1}_{\Omega^1,\Omega^1})\nabla_{Q(\Omega^1\tens_0\Omega^1)}Q(\tilde g)(1)
\end{eqnarray*}
and the right hand side is zero since $\nabla_{\Omega^1\tens_0\Omega^1}g=0$.  One can also see this another way, which some readers may prefer: By Lemma~\ref{qfunct} (which is best summarised by the commuting diagram (\ref{vcusvdsavbbdfs})), as long as the corresponding $q$s are inserted, the tensor product of the quantised connections is the same as the quantisation of the tensor product connection. We take a special case of (\ref{vcusvdsavbbdfs}), remembering that $\Omega^1A_1=Q(\Omega^1(M))$.
\begin{eqnarray}\label{vcusvdsavbbd}
\xymatrix{
Q(\Omega^1 (M))\tens_1 Q(\Omega^1 (M)) \ar[r]^{q} \ar[d]_{\nabla_{Q\tens_1 Q}} 
&Q(\Omega^{\tens 2} (M))\ar[d]_{\nabla_{Q(\Omega^{\tens 2} (M))}}      \\
\Omega^1A_1\tens_1 Q(\Omega^1 (M))\tens_1 Q(\Omega^1 (M)) \ar[r]_{\qquad\id\tens q} &
\Omega^1A_1\tens_1  Q(\Omega^{\tens 2} (M))    }
\end{eqnarray}
Now we suppose that classically the quantising connection preserves the classical Riemannian metric $g\in \Omega^{\tens 2} (M)$, i.e.\ that
$\nabla_{\Omega^{\tens 2} (M)}g=0$. By Lemma~\ref{hilskljhvcj} we have $\nabla_{Q(\Omega^{\tens 2} (M))}g=0$, which also gives $g$ central in the quantised system. Also by (\ref{vcusvdsavbbd}) we see that $g_Q=q^{-1}g\in \Omega^1 A_1\tens_1\Omega^1 A_1$ is indeed preserved by the tensor product of the quantised connections $\nabla_{Q\tens_1 Q}$. Moreover, we already know from Lemma~\ref{cluiktycxtu} that $\nabla_Q$ preserves the star operation hence in this case we also have  Hermitian-metric compatibility with $g_Q$ in the sense  \begin{eqnarray*}
(\bar\nabla_{Q}\tens\id+\id\tens\nabla_{Q})(\star\tens\id)g_Q=0\ .
\end{eqnarray*}
Over $\C$, reality of $g_Q$ in the sense $\Upsilon^{-1}(\star\tens_1\star)\,g_Q=\overline{g_Q}$ reduces by (\ref{vcusvdsavbbdfsii})  to the classical statement for $g$ and $\star\tens_0\star$, which is trivial certainly if the classical coefficients $g_{ij}$ are real and symmetric. 
 \endproof

However, $g_Q$ is not necessarily `quantum symmetric'. We can correct for this by an adjustment at order $\lambda$.

\begin{proposition}\label{nablaQCR} Let $(\omega,\nabla)$ be a Poisson tensor with Poisson-compatible connection and
define the associated `generalised Ricci 2-form' and adjusted metric
\[ \CR=g_{ij}H^{ij},\quad g_1=g_Q-\lambda q^{-1}\CR\]
where on the right the 2-form is lifted to an antisymmetric tensor. Suppose that $\nabla g=0$.
\begin{enumerate}
\item If the lowered $T_{ijk}$ is totally antisymmetric  then $\extd\CR=0$. 
\item  $\wedge_1(g_1)=0$, and $q^2\nabla_Q g_1=-\lambda\nabla\CR$. Here $\nabla_Qg_1=0$ if and only if $\nabla \CR=0$.
\item Over $\C$, $g_1$ is `real' (and $\nabla_Q$ is star-preserving).
\end{enumerate}
\end{proposition}
\proof (1) We use Lemma~\ref{bchisaokjhgc} in the following formula,
\begin{eqnarray*}
\extd(g_{ij}\,H^{ij}) &=& g_{ij,p}\,\extd x^p\wedge H^{ij}+g_{ij}\,\extd H^{ij}   \cr
&=& g_{ij,p}\,\extd x^p\wedge H^{ij}-g_{ij}\,(\Gamma^i_{rp}\,\extd x^p\wedge H^{rj} +\Gamma^j_{rp}\,\extd x^p\wedge H^{ir}) \cr
&=& (g_{ij,p} - g_{rj}\,\Gamma^r_{ip}   -  g_{ir}\,  \Gamma^r_{jp})\,\extd x^p\wedge H^{ij}\ .
\end{eqnarray*}
If $\nabla$ preserves the metric we also have
\begin{eqnarray*}
0 &=& \nabla_p(g_{ij}\,\extd x^i\tens\extd x^j)=(g_{ij,p}-g_{rj}\,\Gamma^r_{pi}-g_{ir}\,\Gamma^r_{pj})\,\extd x^i\tens\extd x^j\ ,
\end{eqnarray*}
and using this, if the lowered $T_{ijk}$ is totally antisymmetric
\begin{eqnarray} \label{uyvuytdc}
\extd(g_{ij}\,H^{ij}) 
&=& (g_{rj}\,\Gamma^r_{pi}+g_{ir}\,\Gamma^r_{pj} - g_{rj}\,\Gamma^r_{ip}    -  g_{ir}\,  \Gamma^r_{jp})\,\extd x^p\wedge H^{ij}\cr
&=& (g_{rj}\,T^r_{pi}+g_{ir}\,T^r_{pj} )\,\extd x^p\wedge H^{ij}\cr
&=& (T_{jpi}+T_{ipj} )\,\extd x^p\wedge H^{ij} =0\ .
\end{eqnarray}
(2) Clearly $\wedge_1(g_Q)=\lambda\CR$ so $\wedge_1(g_1)=0$. Likewise $q^2\nabla_Qg_1=q^2\nabla_Qg_Q-\lambda\nabla\CR=-\lambda\nabla\CR$ by Lemma~\ref{lemmagQ}, where the last term here is viewed as an element of $\Omega^1(M)^{\tens_0 3}$ by an antisymmetric lift. The antisymmetric lift commutes with $\nabla$ so $\nabla_1g_1=0$ if and only if  $\nabla\CR=0$ on $\CR$ as a 2-form.  To give the formulae here more explicitly, we remember our 2-form conventions so that
\begin{equation}\label{CRtensor} \CR={1\over 2}\CR_{nm}\extd x^m\wedge\extd x^n,\quad \CR_{mn}=\tfrac{1}{2}\,g_{ij} \omega^{is}(T^j_{nm;s}- R^j{}_{nms}+R^j{}_{mns}).\end{equation}
in which case, 
\[ g_1=g_Q+{\lambda\over 2}\CR_{mn} \extd x^m\tens_1\extd x^n.\]
 (3) Over $\C$, we also have the condition $\Upsilon^{-1}(\star\tens_1\star)g_1=\overline{g_1}$, as the correction is both imaginary and antisymmetric. $\nabla_Q$ is still star-preserving because that statement is not dependent on the metric (which means that it is also Hermitian-metric compatible with the corresponding Hermitian metric $(\star\tens\id)g_1$).   \endproof

In general we may not have either of these properties of $\CR$ but we do have $\nabla_Q g_1$ being  order $\lambda$  and that is enough to make $g_1$ commute with elements of $A_1$ to  order $\lambda$ which is what we wanted to retain at this point. The terminology for $\CR$ comes from the K\"ahler case which is a subcase of the following special case.

\begin{corollary} \label{leviequal} If the quantising connection $\nabla$ is the Levi-Civita one, 
\newline\noindent (1)\quad Poisson-compatibility reduces to $\omega$ covariantly constant. 
\newline\noindent (2)\quad $\nabla_Q$ is quantum torsion free and $\CR={1\over 2}\omega^{ji}R_{inmj}\extd x^m\wedge\extd x^n$ is closed.
\newline\noindent (3)\quad  $\nabla_Qg_1=0$, i.e. $\nabla_Q$  is a quantum-Levi-Civita connection for $g_1$, if and only if $\nabla\CR=0$.
\end{corollary}
\proof This is a special caae of Proposition~\ref{nablaQCR}.  For the quantum torsion we use Proposition~\ref{torsionQ} where the torsion $T$ of $\nabla$ is currently being assumed to be zero. In this case $\extd \CR=0$ as $T=0$ is antisymmetric. Note that if $\nabla\CR\ne 0$ we still have $\nabla_Qg_1$ is order $\lambda$ by Lemma~\ref{lemmagQ}. \endproof

\subsection{Relating general  $\nabla$ and the Levi-Civita $\widehat{\nabla}$}

In general the quantising connection $\nabla$ may not be the same as the classical Levi-Civita connection $\widehat\nabla$ for our chosen metric on $M$. In this section we write the latter in the general form $\nabla_S=\nabla+S$ for some $S:\Omega^1 (M)\to \Omega^1 (M)\tens_0 \Omega^1 (M)$ and we assume that the quantising connection $\nabla$ obeys $\nabla g =0$.  The quantising connection has torsion $T$  and we lower its indices by the Riemannian metric $T_{abc}=g_{ad}\,T^d_{bc}$. It is well-known (see \cite{hehl1}) that given an arbitrary torsion $T$, there is a unique metric compatible covariant derivative $\nabla$ with that torsion, given by 
\begin{eqnarray}  \label{bcilskjc}
\Gamma^a_{bc}=\widehat \Gamma^a_{bc}+\tfrac12 g^{ad}(T_{dbc}-T_{bcd}-T_{cbd})\ .
\end{eqnarray}
Here $\Gamma^a_{bc}$ in our case is the Christoffel symbols for the quantising connection and $\widehat \Gamma^a_{bc}$ is the Christoffel symbols for the Levi-Civita connection so that $\nabla_S(\extd x^a)=-\widehat \Gamma^a_{bc}\,\extd x^b\tens \extd x^c$. Hence 
\begin{equation}\label{SfromT}S^a_{bc}=\tfrac12 g^{ad}(T_{dbc}-T_{bcd}-T_{cbd}).\end{equation}
As a quick check of conventions, note that this formula is consistent with (\ref{bhcuiiycxyr}).  Throughout this section $T$ is arbitrary which fixes $\nabla$ such that this is metric compatible, and $S$ is the above function of $T$ so that
$\nabla_S=\widehat\nabla$, the Levi-Civita connection. 

\begin{lemma}\label{RLC}
The curvatures are related by
\begin{eqnarray*}
\widehat R^{l}_{\phantom{l}ijk}
&=& R^{l}_{\phantom{l}ijk} - S^l_{ki;j}+ S^l_{ji;k} - T^m_{jk}\, S^l_{mi}
+  S^m_{ki}\, S^l_{jm}-S^m_{ji}\, S^l_{km}  \ ,
\end{eqnarray*}
where semicolon is derivative with respect to $\nabla$. \end{lemma}
\proof This is elementary: $\widehat \Gamma^m_{ji}=\Gamma^m_{ji}-S^m_{ji}$ so that 
\begin{eqnarray*}
\widehat R^{l}_{\phantom{l}ijk} &=& \widehat  \Gamma^l_{ki,j}- \widehat \Gamma^l_{ji,k}+
\widehat \Gamma^m_{ki}\,\widehat \Gamma^l_{jm}-\widehat \Gamma^m_{ji}\,\widehat \Gamma^l_{km}\cr
&=& R^{l}_{\phantom{l}ijk} - S^l_{ki,j}+ S^l_{ji,k}- \Gamma^m_{ki}\,S^l_{jm}+\Gamma^m_{ji}\,S^l_{km}- S^m_{ki}\,\Gamma^l_{jm}+S^m_{ji}\,\Gamma^l_{km}\cr
&& +\,  S^m_{ki}\,S^l_{jm}-S^m_{ji}\,S^l_{km}  \cr
&=& R^{l}_{\phantom{l}ijk} - S^l_{ki;j}+ S^l_{ji;k} - T^m_{jk}\,S^l_{mi}
+  S^m_{ki}\,S^l_{jm}-S^m_{ji}\,S^l_{km}  \ .\quad\square\end{eqnarray*}

\medskip
This gives a different point of view on some of the formulae below, if we wish to rewrite expressions in terms of the Levi-Civita connection. In the same vein:

\begin{proposition}  \label{leviomega}
Suppose that a connection $\nabla$ is metric-compatible. Then $(\nabla,\omega)$ are Poisson-compatible if and only if 
\[ (\widehat\nabla_k\omega)^{ij}+ \omega^{ir}\, S^j_{rk} -\, \omega^{jr} S^i_{rk} =0
\] 
or equivalently
\begin{eqnarray*}
\omega^{jm}S^i_{mk}
 &=&\tfrac 12\left( (\widehat\nabla_k\omega)^{ij}\, - (\widehat\nabla_r\omega)^{mj}\, g^{ri}\,g_{mk}
+(\widehat\nabla_r\omega)^{im}\,g^{rj}\,g_{mk}\  \right) .
\end{eqnarray*}
\end{proposition}
\proof 
The compatibility condition gives
\begin{eqnarray*}
0 &=& (\widehat\nabla_m\omega)^{ij}+\omega^{ik}\,(T^j_{km}+\tfrac12 g^{jd}(T_{dmk}-T_{mkd}-T_{kmd})) \cr
&& + \, \omega^{kj}\,(T^i_{km}+\tfrac12 g^{id}(T_{dmk}-T_{mkd}-T_{kmd})) \cr
&=& (\widehat\nabla_m\omega)^{ij}+\omega^{ik}\,\tfrac12 g^{jd}(T_{dkm}-T_{mkd}-T_{kmd}) \cr
&& + \, \omega^{kj}\,\tfrac12 g^{id}(T_{dkm}-T_{mkd}-T_{kmd}) \cr
&=& (\widehat\nabla_m\omega)^{ij}+\omega^{ik}\,\tfrac12 g^{jd}(T_{dkm}+T_{mdk}-T_{kmd}) \cr
&& + \, \omega^{kj}\,\tfrac12 g^{id}(T_{dkm}+T_{mdk}-T_{kmd}) \cr
&=& (\widehat\nabla_m\omega)^{ij}+\tfrac12\,(\omega^{ik}\, g^{jd}- \omega^{jk}\,g^{id})(T_{dkm}+T_{mdk}-T_{kmd}) \ 
\end{eqnarray*}
which is the first condition stated in terms of $S$.  From this,
\begin{eqnarray*}
(\widehat\nabla_m\omega)^{ij}\,g_{ir}\,g_{js}&=& - \omega^{ik}\,g_{ir}\, S_{skm}
 +\, \omega^{jk}\,g_{js}S_{rkm} \ .
\end{eqnarray*}
Now define
\begin{eqnarray*}
-\,\Theta_{mrs} &:=& (\widehat\nabla_m\omega)^{ij}\,g_{ir}\,g_{js} -  \omega^{jk}\,g_{js} 2  S_{rkm} = -\, \omega^{ik}\,g_{ir}\, S_{skm} - \omega^{jk}\,g_{js} S_{rkm}
\end{eqnarray*}
and note that $\Theta_{mrs}$ is symmetric on swapping $r,s$. Rearranging this gives
\begin{eqnarray}\label{bcuovuv}
\omega^{jk}\,g_{js}2\, S_{rkm}  &=& (\widehat\nabla_m\omega)^{ij}\,g_{ir}\,g_{js} +\Theta_{mrs}\ ,\cr
\omega^{jk}2\, S_{rkm}  &=& (\widehat\nabla_m\omega)^{ij}\,g_{ir}+\Theta_{mrs}\,g^{sj}\ .
\end{eqnarray}
We can also write
\begin{eqnarray*}
\Theta_{mrs} &=& -\,(\widehat\nabla_m\omega)^{ij}\,g_{ir}\,g_{js} +  \omega^{jk}\,g_{js}2\, S_{rkm}\ ,
\end{eqnarray*}
and from this we get the following condition, which we repeat with permuted indices
\begin{eqnarray} \label{ncjdvfb}
\Theta_{mrs}+\Theta_{rms} &=& -\,(\widehat\nabla_m\omega)^{ij}\,g_{ir}\,g_{js}
-(\widehat\nabla_r\omega)^{ij}\,g_{im}\,g_{js}\ ,\cr
\Theta_{rsm}+\Theta_{srm} &=& -\,(\widehat\nabla_r\omega)^{ij}\,g_{is}\,g_{jm}
-(\widehat\nabla_s\omega)^{ij}\,g_{ir}\,g_{jm}\ ,\cr
\Theta_{smr}+\Theta_{msr} &=& -\,(\widehat\nabla_s\omega)^{ij}\,g_{im}\,g_{jr}
-(\widehat\nabla_m\omega)^{ij}\,g_{is}\,g_{jr}\ .
\end{eqnarray}
Taking the first line of (\ref{ncjdvfb}), subtracting the second and adding the third gives
\begin{eqnarray}
 \Theta_{mrs} &=&
 (\widehat\nabla_r\omega)^{ij}\,g_{is}\,g_{jm}
+(\widehat\nabla_s\omega)^{ij}\,g_{ir}\,g_{jm}\ .
\end{eqnarray}
Now we rewrite (\ref{bcuovuv}) as
\begin{eqnarray*}
\omega^{jk}2\, S_{rkm}  &=& (\widehat\nabla_m\omega)^{ij}\,g_{ir}+( (\widehat\nabla_r\omega)^{it}\,g_{is}\,g_{tm}
+(\widehat\nabla_s\omega)^{it}\,g_{ir}\,g_{tm})\,g^{sj}\cr
 &=& (\widehat\nabla_m\omega)^{ij}\,g_{ir}- (\widehat\nabla_r\omega)^{ij}\,g_{im}
+(\widehat\nabla_s\omega)^{it}\,g_{ir}\,g_{tm}\,g^{sj}\ 
\end{eqnarray*}
which we write as stated. \endproof

\subsection{ Metric compatibility in the general case} \label{secmetriccompat}
Now we look for a quantum Levi Civita connection in the general case where the quantising connection $\nabla$ is not the Levi-Civita one. 
As  in Section~\ref{bchijadflkjhov} we assume a metric $g\in\Omega^{1^{\tens 2}}(M)$ and $\nabla g=0$ and as in Section~5.2 we let $S$ be a function of $T$ so that $\nabla_S=\nabla+S=\widehat\nabla$, the classical Levi-Civita connection for $g$.  We do the straight metric compatibility in this section (which makes sense over any field) and the Hermitian version in the next section (recall that the two versions of the metric-compatibility coincide if the quantum connection is star-preserving).

\begin{lemma}\label{QSmetrictor} For $\nabla_S$ the Levi-Civita connection, the quantum metric compatibility tensor and quantum torsion $T_{\nabla_{QS}}=\tfrac\lambda2\, \xi_p A^p_{nm} \extd x^m\wedge \extd x^n$ are given respectively by
\[
q^2\,\nabla_{QS\tens_1 QS}(g_Q) =-\,\lambda\,\omega^{ij}\,g_{rs}\,S^s_{jn}(R^r{}_{mki}+S^r_{km;i})(\extd x^k\tens\extd x^m \tens \extd x^n)\]
\[ A^p_{nm}=- \tfrac14\omega^{ij}\left( g^{pd}\,(T_{isd}+T_{sid})\,( T^s_{nm;j}-R^s{}_{nmj} +R^s{}_{mnj}) +T^s_{nm}\,R^p{}_{sij}-T^p_{sm}\,R^s{}_{nij}+T^p_{sn}\,R^s{}_{mij}\right).\]
\end{lemma}
\proof We look at Proposition~\ref{jdolobvlo}, and set
\begin{eqnarray*}
H=S\tens\id_F+(\tau\tens\id)(\id\tens S):\Omega^{\tens 2} (M)\to \Omega^1 (M)\tens_0 \Omega^{\tens 2} (M)\ .
\end{eqnarray*}
As classically both $\nabla$ and $\nabla_S$ preserve $g$, we get $H(g)=0$. By Lemma~\ref{hilskljhvcj} again, we get $Q(H)(g)=0$ and $\nabla_{QH}(g)=0$. Now applying Proposition~\ref{jdolobvlo} gives
\begin{eqnarray}\label{bciosouiyvviu}
q^2\,\nabla_{QS\tens_1 QS}(q^{-1}g) &=& (q\,\nabla_{QH}\,q+\lambda\,\mathrm{rem})(q^{-1}g)=\lambda\,\mathrm{rem}(q^{-1}g)\ ,
\end{eqnarray}
where, using $S(f)=\extd x^k\tens S_k(f)$
\begin{eqnarray*}
\mathrm{rem}(e\tens_1 f) = \omega^{ij}\,( \extd x^k \tens [\nabla_k,\nabla_i]e
- \nabla_i(S)(e)  )  \tens S_j(f)\ .
\end{eqnarray*}
 We now have, by (\ref{bciosouiyvviu}),
\begin{eqnarray*}
q^2\,\nabla_{QS\tens_1 QS}(q^{-1}g) &=& \lambda\,\mathrm{rem}(q^{-1}g)=
 \lambda\,\mathrm{rem}(g_{rs}\,\extd x^r\tens_1\extd x^s) \cr
 &=&\lambda\, \omega^{ij}\,g_{rs}\,( \extd x^k \tens [\nabla_k,\nabla_i](\extd x^r)
- \nabla_i(S)(\extd x^r)  )  \tens S_j(\extd x^s) \cr
 &=&\lambda\, \omega^{ij}\,g_{rs}\,S^s_{jn}\, ( \extd x^k \tens [\nabla_k,\nabla_i](\extd x^r)
- \nabla_i(S)(\extd x^r)  )  \tens \extd x^n \ 
\end{eqnarray*}
which we write as stated. For the torsion we used  $S^p_{ij}+S^p_{ji}=-\,g^{pd}(T_{ijd}+T_{jid})$ in Corollary~\ref{kjgccghgc} and relabelled. \endproof

We see that the quantisation  $\nabla_{QS}$ given by the procedure outlined in Section~\ref{jytcxfufgx} is only quantum metric compatible to an error of order $\lambda$. However we have the freedom to add an  order $\lambda$ correction to $g_Q$ as above and an  order $\lambda$ correction to the proposed quantum connection:

\begin{theorem} \label{qlevi} Let $\nabla_S$ be the Levi-Civita connection. There is a unique  quantum connection of the form $\nabla_1=\nabla_{QS}+\lambda K$ 
such  that the quantum torsion and merely the symmetric part of $\nabla_1g_1$ vanish. The antisymmetric part,
\begin{eqnarray*}  
(\id\tens\wedge)q^2\,\nabla_1g_1& =& -\lambda\widehat\nabla\CR -\,\lambda\,\omega^{ij}\,g_{rs}\,S^s_{jn}(R^r{}_{mki}+S^r_{km;i})\,\extd x^k\tens\extd x^m \wedge \extd x^n           \ ,
\end{eqnarray*}
is independent of $K$. A fully metric compatible torsion free $\nabla_1$ exists if and only if the above expression vanishes, in which case it is given by the unique $\nabla_1$ discussed. 
\end{theorem}
\proof We  write $K(\xi)=\xi_p\,K^p_{nm}\,\extd x^n\tens \extd x^m$, then (where semicolon is given by the quantising connection) the results in the preceding lemma are
clearly adjusted to 
\begin{eqnarray*}  \label{kjhvcxjhgc}
q^2\,\nabla_{1}(g_1)& =& -\,\lambda\,\omega^{ij}\,g_{rs}\,S^s_{jn}(R^r{}_{mki}+S^r_{km;i})\,\extd x^k\tens\extd x^m \tens \extd x^n           \cr
&& - \, \tfrac{\lambda}{4}\,\nabla_{S\tens S}(g_{ij} \, \omega^{is}(T^j_{nm;s}- R^j{}_{nms}+R^j{}_{mns})\, \extd x^m\tens\extd x^n) \cr
&&+\, \lambda\,(  g_{pn}\,K^p_{km}+g_{mp}\,K^p_{kn}   )\, \extd x^k\tens\extd x^m \tens \extd x^n\  \cr
T_{\nabla_1}(\xi) &=&
\tfrac\lambda2\, \xi_p\,( K^p_{nm} -K^p_{mn} - A^p_{nm})\,\extd x^n\wedge \extd x^m  \ .
\end{eqnarray*}
Looking at the first expression reveals that the second term is purely antisymmetric in $nm$, whereas the third term (the only one to contain the  order $\lambda$ correction $K^a_{bc}$) is purely symmetric in $nm$. Hence there is nothing we can do by adding $K^a_{bc}$ to make the part of the metric compatibility tensor which is antisymmetric in $nm$ vanish, it will have the value stated, but we show that we can choose $K^a_{bc}$ to make the part which is symmetric in $nm$ vanish, namely by setting
\begin{eqnarray*} 
  g_{np}\,K^p_{km}+g_{mp}\,K^p_{kn}    & =& B_{knm}
\end{eqnarray*}
where
\begin{eqnarray*} 
B_{knm}   & =& \tfrac12\,\omega^{ij}\,g_{rs}\,(S^s_{jn}(R^r{}_{mki}+S^r_{km;i}) 
  +S^s_{jm}(R^r{}_{nki}+S^r_{kn;i}) )     \ 
\end{eqnarray*}
while for vanishing torsion, clearly we need $K^p_{nm} -K^p_{mn} =A^p_{nm}$. If we set $K_{nkm}=g_{np}\,K^p_{km}$ then these
conditions become
\begin{eqnarray*}
K_{nkm}+K_{mkn}=B_{knm}\ ,\quad K_{knm}-K_{kmn}=g_{kp}\,A^p_{nm}\ .
\end{eqnarray*}
Now
\begin{eqnarray*}
K_{nkm} &=& B_{knm} - K_{mkn}=B_{knm} +g_{mp}\,A^p_{nk} - K_{mnk}\ ,
\end{eqnarray*}
and continuing in this manner six times gives a unique value of $K$, 
\begin{eqnarray} \label{koszul}
K_{nkm} &=& {1\over 2}\left(B_{knm} - B_{nkm}+ B_{mnk} +g_{mp}\,A^p_{nk} +g_{kp}\,A^p_{nm} 
 +g_{np}\,A^p_{km}\right) \ .
\end{eqnarray}
where
\begin{eqnarray*}
&& B_{knm}  - B_{nkm} 
+ B_{mnk} \cr
&=&  \tfrac12\,\omega^{ij}\,g_{rs}\,(S^s_{jn}(R^r{}_{mki}+S^r_{km;i}) 
  +S^s_{jm}(R^r{}_{nki}+S^r_{kn;i}) )   \cr
&& -\, \tfrac12\,\omega^{ij}\,g_{rs}\,(S^s_{jk}(R^r{}_{mni}+S^r_{nm;i}) 
  +S^s_{jm}(R^r{}_{kni}+S^r_{nk;i}) ) \cr
  && + \,  \tfrac12\,\omega^{ij}\,g_{rs}\,(S^s_{jn}(R^r{}_{kmi}+S^r_{mk;i}) 
  +S^s_{jk}(R^r{}_{nmi}+S^r_{mn;i}) )   \cr
&=&  \tfrac12\,\omega^{ij}\,g_{rs}\,(S^s_{jn}(R^r{}_{mki}+R^r{}_{kmi}+S^r_{mk;i}+S^r_{km;i}) 
  +S^s_{jm}(R^r{}_{nki}-R^r{}_{kni}-S^r_{nk;i}+S^r_{kn;i})    \cr
 && + \, S^s_{jk}(R^r{}_{nmi}-R^r{}_{mni}-S^r_{nm;i}+S^r_{mn;i}) )   \cr
&=&  \tfrac12\,\omega^{ij}\,g_{rs}\,(S^s_{jn}(R^r{}_{mki}+R^r{}_{kmi} - g^{rd}\, (T_{mkd;i}+T_{kmd;i})) 
  +S^s_{jm}(R^r{}_{nki}-R^r{}_{kni}+T^r_{kn;i})    \cr
 && + \, S^s_{jk}(R^r{}_{nmi}-R^r{}_{mni}+T^r_{mn;i}) )  \end{eqnarray*}
 using 
$S^r_{mk}+S^r_{km}=-\,g^{rd}\, (T_{mkd}+T_{kmd})$. \endproof

\medskip
This clearly reduces to Corollary~\ref{leviequal} in the case where $T=0$ but  more generally we have a free parameter, the value of $T$ for the quantizing connection provided only that $(\nabla,\omega)$ are Poisson-compatible. We might hope to use this freedom to set $R=0$ so that our differential calculus remains associative at the next order in $\lambda$, and/or we might hope to choose $T$ so that the the antisymmetric part of the quantum metric compatibility tensor also vanishes. In the nontrivial example black-hole below this expression will not even depend on $T$ within the class discussed, i.e. can have more of a topological character. Hence we  can't always obtain full metric compatibility but rather can have an unavoidable quantum correction. In that case we still have a `best possible' choice of $\nabla_1$ given by the formula (\ref{koszul}).

\subsection{Hermitian-metric compatibility}

Here we again assume that  our Poisson-compatible connection $\nabla$ obeys $\nabla g=0$ and that $\nabla_S=\nabla+S$ is the  Levi-Civita connection for $g$. We set $\nabla_1=\nabla_{QS}+\lambda K$, for some real $K$, and ask this time that $\nabla_1$ is Hermitian-metric compatible with the Hermitian metric $(\star \tens\id)g_1$ corresponding to $g_1$. This is a potentially different condition from straight metric compatibility unless  $\nabla_1$ is star-preserving, in which case it is equivalent.

\begin{proposition} \label{hermetriccompat}
Over $\C$ and with $\nabla_S$ the Levi-Civita connection,  the condition for $\nabla_{QS}+\lambda K$ to be Hermitian-metric compatible with $g_1$ is, where $\hat;$ denotes the Levi-Civita derivative,
\begin{eqnarray*}
K_{npm}-K_{mpn} = \CR_{nm\hat;p}+ \tfrac12\, \omega^{ij}\, ( g_{rm}\,\nabla_i(\nabla_j(S))^r_{pn}-g_{nr}\,\nabla_i(\nabla_j(S))^r_{pm})\ .
\end{eqnarray*}
This can always be solved simultaneously with vanishing of the quantum torsion.
\end{proposition}
\proof  
(1) If we write the quantum correction to the metric in Proposition~\ref{nablaQCR} as $g_1=g_Q-\lambda\,g_c$, then Hermitian-metric compatibility tensor for $\nabla_{QS}$ becomes
\begin{eqnarray}\label{hermmetric1}
&& ((\id\tens\star^{-1})\Upsilon \, \overline{q^{-1}Q(S)}\tens\id+\id\tens q^{-1}Q(S))(\star\tens\id)g_Q \cr
&& -\lambda(\bar\nabla_{S}\tens\id+\id\tens\nabla_{S})(\star\tens\id)g_c
\end{eqnarray}
From Proposition~\ref{kihjvfcyiouy} we can write this as
\begin{eqnarray*}
&& (q^{-1}(\id\tens\star^{-1})\Upsilon \, \overline{Q(S)}\tens_1\id+\id\tens_1 q^{-1}Q(S))(\star\tens\id)g_Q \cr
&& -\lambda(\bar\nabla_{S}\tens\id+\id\tens\nabla_{S})(\star\tens\id)g_c\ .
\end{eqnarray*}
The definition of $Q(S)$ gives
\begin{eqnarray*}
&& (q^{-1}(\id\tens\star^{-1})\Upsilon \, \overline{S}\tens_1\id+\id\tens_1 q^{-1}S)(\star\tens\id)g_Q \cr
&& +\tfrac\lambda2\, \omega^{ij}\, (-\,(\id\tens\star^{-1})\Upsilon \, \overline{\nabla_i\circ\nabla_j(S)}\tens\id+\id\tens \nabla_i\circ\nabla_j(S))(\star\tens\id)g \cr
&& -\,\lambda((\id\tens\star^{-1})\Upsilon \, \overline{\nabla_{S}}\tens\id+\id\tens\nabla_{S})(\star\tens\id)g_c\ .
\end{eqnarray*}
Now we use Proposition~\ref{q} and apply $q^2$, noting that
\begin{eqnarray*} \label{oiuyccxru}
&& q^2(q^{-1}(\id\tens\star^{-1})\Upsilon \, \overline{S}\tens_1\id+\id\tens_1 q^{-1}S)(\star\tens\id)g_Q \cr
&=& q((\id\tens\star^{-1})\Upsilon \, \overline{S}\tens_1\id+\id\tens_1 S)(\star\tens\id)\,q^{-1}g \cr
&=& ((\id\tens\star^{-1})\Upsilon \, \overline{S}\tens_0\id+\id\tens_0 S)(\star\tens\id)\,g \cr
&&+\, \tfrac\lambda2\,\omega^{ij}\, 
(\nabla_i((\id\tens\star^{-1})\Upsilon \, \overline{S})\tens_0\nabla_j+\nabla_i\tens_0 \nabla_j(S))(\star\tens\id)\,g \cr
&=& ((\id\tens\star^{-1})\Upsilon \, \overline{S}\tens_0\id+\id\tens_0 S)(\star\tens\id)\,g \cr
&&+\, \tfrac\lambda2\,\omega^{ij}\, 
((\id\tens\star^{-1})\Upsilon \, \overline{\nabla_i(S)}\tens_0\nabla_j+\nabla_i\tens_0 \nabla_j(S))(\star\tens\id)\,g \cr
&=& ((\id\tens\star^{-1})\Upsilon \, \overline{S}\tens_0\id+\id\tens_0 S)(\star\tens\id)\,g \cr
&&-\, \tfrac\lambda2\,\omega^{ij}\, 
((\id\tens\star^{-1})\Upsilon \, \overline{\nabla_i(S)\circ \nabla_j}\tens_0\id+\id\tens_0 \nabla_j(S)\circ\nabla_i)(\star\tens\id)\,g \ 
\end{eqnarray*}
as $g$ is preserved by $\nabla$ so that $(\nabla_i\tens\id)g=-(\id\tens\nabla_i)g$.

Then $q^2$ applied to (\ref{hermmetric1}) gives
\begin{eqnarray}  \label{kuyccy}
&& ((\id\tens\star^{-1})\Upsilon \, \overline{S}\tens_0\id+\id\tens_0 S)(\star\tens\id)\,g \cr
&&-\, \tfrac\lambda2\,\omega^{ij}\, 
(-(\id\tens\star^{-1})\Upsilon \, \overline{\nabla_j(S)\circ \nabla_i}\tens_0\id+\id\tens_0 \nabla_j(S)\circ\nabla_i)(\star\tens\id)\,g \cr
&& +\tfrac\lambda2\, \omega^{ij}\, (-\,(\id\tens\star^{-1})\Upsilon \, \overline{\nabla_i\circ\nabla_j(S)}\tens\id+\id\tens \nabla_i\circ\nabla_j(S))(\star\tens\id)g \cr
&& -\,\lambda((\id\tens\star^{-1})\Upsilon \, \overline{\nabla_{S}}\tens\id+\id\tens\nabla_{S})(\star\tens\id)g_c\cr
&=& ((\id\tens\star^{-1})\Upsilon \, \overline{S}\tens_0\id+\id\tens_0 S)(\star\tens\id)\,g \cr
&& +\tfrac\lambda2\, \omega^{ij}\, (-\,(\id\tens\star^{-1})\Upsilon \, \overline{\nabla_i(\nabla_j(S))}\tens\id+\id\tens \nabla_i(\nabla_j(S)))(\star\tens\id)g \cr
&& -\,\lambda((\id\tens\star^{-1})\Upsilon \, \overline{\nabla_{S}}\tens\id+\id\tens\nabla_{S})(\star\tens\id)g_c\ .
\end{eqnarray}
Now set $g=g_{nm}\,\extd x^n\tens\extd x^m$ and $\nabla_i(\nabla_j(S))(\extd x^a)=\nabla_i(\nabla_j(S))^a_{nm}\,\extd x^n\tens\extd x^m$, and using the reality of $S$ the first two lines of the result of (\ref{kuyccy}) become
\begin{eqnarray*}
&& (g_{rm}\,S^r_{pn}+g_{nr}\,S^r_{pm})\,\overline{\extd x^n} \tens \extd x^p \tens \extd x^m \cr
&&  +\tfrac\lambda2\, \omega^{ij}\, ( -\,g_{rm}\,\nabla_i(\nabla_j(S))^r_{pn}+g_{nr}\,\nabla_i(\nabla_j(S))^r_{pm})\,\overline{\extd x^n} \tens \extd x^p \tens \extd x^m\ ,
\end{eqnarray*}
and the first line of this vanishes as $\nabla_S$ preserves $g$. Now we write $g_c=-{1\over 2}\CR_{nm}\,\extd x^n\tens\extd x^m$ where $\CR_{nm}$ is antisymmetric giving
\[  q^2(\bar\nabla_{QS}\tens\id+\id\tens\nabla_{QS})(\star\tens\id)g_1=-\ \lambda\, C_{npm}\,\overline{\extd x^n} \tens \extd x^p \tens \extd x^m; \]
\[C_{npm}=-\tfrac12\,\CR_{nm\hat;p}+ \tfrac12\, \omega^{ij}\, ( g_{rm}\,\nabla_i(\nabla_j(S))^r_{pn}-g_{nr}\,\nabla_i(\nabla_j(S))^r_{pm}) .\ \]

(2) 
Now we look at $\nabla_1=\nabla_{QS}+\lambda K$, then clearly
\begin{eqnarray*} \label{klhjvciyufew}
&&(\bar\nabla_{1}\tens\id+\id\tens\nabla_{1})(\star\tens\id)g_1=  \lambda\,(g_{na}\, K^a_{pm}-
g_{ma}\,K^a_{pn}- C_{npm})\,\overline{\extd x^n} \tens \extd x^p \tens \extd x^m
\end{eqnarray*}
so we need to solve $K_{npm}-K_{mpn}=C_{npm}$ to preserve the Hermitian metric, and also
$K_{knm}-K_{kmn}=g_{ks}\,A^s_{nm}$ if we want to have zero torsion as in the previous section. These equations have a required compatibility condition
\begin{eqnarray*}
C_{npm}+C_{mnp}+C_{pmn}+g_{ms}\,A^s_{pn}+g_{ps}\,A^s_{nm}+g_{ns}\,A^s_{mp}=0\ .
\end{eqnarray*}
We use the formula (\ref{SfromT}) for $S^a_{bc}$ in terms of the torsion to write
\begin{eqnarray*}
C_{npm} &=& -\tfrac12\, \CR_{nm\hat;p}+ \tfrac12\, \omega^{ij}\, ( g_{rm}\,S^r_{pn;\hat j i}-g_{nr}\,S^r_{pm;\hat j i}) \cr
&=& -\tfrac12\, \CR_{nm\hat;p}+ \tfrac14\, \omega^{ij}\, ( (T_{mpn}-T_{pnm}-T_{npm})-(T_{npm}-T_{pmn}-T_{mpn}))_{;\hat j i} \cr
&=& -\tfrac12\, \CR_{nm\hat;p}+ \tfrac12\, \omega^{ij}\,  (T_{mpn}-T_{pnm}-T_{npm})_{;\hat j i} \ ,
\end{eqnarray*}
and taking the cyclic sum gives
\begin{eqnarray*}
C_{npm}+C_{mnp}+C_{pmn} &=& -\tfrac12\, \CR_{nm\hat;p}-\tfrac12\, \CR_{pn\hat;m} -\tfrac12\, \CR_{mp\hat;n}+ \tfrac12\, \omega^{ij}\,  (T_{mpn}+T_{pnm}-T_{npm})_{;\hat j i} \ .
\end{eqnarray*}
We have
\begin{eqnarray*}
-\,g_{pa}\,A^a_{nm} 
&=& \tfrac14\omega^{ij}\,(T_{isp}+T_{sip})\,( T^s_{nm;j}-R^s{}_{nmj}+R^s{}_{mnj})  \cr
&&+\, \tfrac14\,\omega^{ij}\,g_{pa}\,\,(T^s_{nm}\,R^a{}_{sij}-T^a_{sm}\,R^s{}_{nij}+T^a_{sn}\,R^s{}_{mij})  \cr
&=&  \tfrac14\omega^{ij}\,(T_{isp}+T_{sip})\,( T^s_{nm;j}-R^s{}_{nmj}+R^s{}_{mnj})  + \tfrac12\,\omega^{ij}\,T_{pnm;\hat ji}\ ,
\end{eqnarray*}
so now the cyclic sum becomes
\begin{eqnarray} \label{oiuyvcfyu}
C_{npm}+C_{mnp}+C_{pmn} &=&
 -\tfrac12\, \CR_{nm\hat;p} -  \tfrac14\omega^{ij}\,(T_{isp}+T_{sip})\,( T^s_{nm;j}-R^s{}_{nmj}+R^s{}_{mnj}) \cr
&&  -\tfrac12\, \CR_{pn\hat;m} -  \tfrac14\omega^{ij}\,(T_{ism}+T_{sim})\,( T^s_{pn;j}-R^s{}_{pnj}+R^s{}_{npj}) \cr
&&  -\tfrac12\, \CR_{mp\hat;n} -  \tfrac14\omega^{ij}\,(T_{isn}+T_{sin})\,( T^s_{mp;j}-R^s{}_{mpj}+R^s{}_{pmj}) \ .
\end{eqnarray}
This is totally antisymmetric in $npm$, so we may equivalently consider the 3-form
\begin{eqnarray*}
\alpha &=& \big( -\tfrac12\, \CR_{nm\hat;p} -  \tfrac14\omega^{ij}\,(T_{isp}+T_{sip})\,( T^s_{nm;j}-R^s{}_{nmj}+R^s{}_{mnj}) \big)\,\extd x^p\wedge\extd x^n\wedge\extd x^m \cr
&=& \extd x^p\wedge(-\tfrac12\, \CR_{nm\hat;p}\,\extd x^n\wedge\extd x^m)+ (T_{isp}+T_{sip})\,\extd x^p\wedge H^{is} \cr
&=& \extd x^p\wedge\widehat\nabla_p(g_{ij}\,H^{ij})+ 2\, T_{isp}\,\extd x^p\wedge H^{is}
\end{eqnarray*}
where we use $H^{ij}=\tfrac14\omega^{is}( T^j_{nm;s}-2 R^j{}_{nms})\extd x^m\wedge\extd x^n$
and the symmetry of $H^{ij}$. Now we have as in (\ref{uyvuytdc}) (but not requiring this to be zero)
\begin{eqnarray*}
\extd \CR&=& (T_{jpi}+T_{ipj} )\,\extd x^p\wedge H^{ij} =-\,2\,T_{ijp}\,\extd x^p\wedge H^{ij}\ ,
\end{eqnarray*}
so vanishing of $\alpha =\extd x^p\wedge\widehat\nabla_p(\CR)-\extd\CR$ is the condition for a joint solution. But this is zero as the Levi-Civita connection is torsion free. \quad$\square$

\medskip Note that Proposition~\ref{hermetriccompat} does not say that such a torsion free quantum connection preserving the Hermitian metric is unique.  If we take the collection of $K_{ijk}$ for permutations of the $ijk$, then the equations fix what the relative value of, for example $K_{ijk}-K_{kij}$ will be, but we can add an overall factor to each of these classes under the permutation group $S_3$.

\begin{corollary}\label{jointcompat}   Over $\C$ and with $\nabla_S$ the Levi-Civita connection,  if a torsion free metric compatible quantum connection of the form $\nabla_1=\nabla_{QS}+\lambda K$ exists, it is star-preserving and coincides with the unique star-preserving quantum connection in Theorem~\ref{starpresK}.
\end{corollary}
\proof
From Lemma~\ref{bgyiuouy1} and Theorem~\ref{starpresK} the star preserving connection is given by $K^a_{nm} = \tfrac14\,\omega^{ij} D^a_{ijnm}$, or
\begin{eqnarray*}
K^a_{nm} &=& \tfrac14\,\omega^{ij}\,(2\, S^a_{ip}\,S^p_{nm;j}-(S^b_{nm}\,R^a{}_{bij}
-S^a_{rm}\,R^r{}_{nij}   - S^a_{nr}\,R^r{}_{mij}   ) -2\,S^a_{jr}\,R^r{}_{mni} )  \cr
&=&  \tfrac14\,\omega^{ij}\,(2\, S^a_{ip}\,S^p_{nm;j} +2\,S^a_{ip}\,R^p{}_{mnj}  -(S^b_{nm}\,R^a{}_{bij}
-S^a_{rm}\,R^r{}_{nij}   - S^a_{nr}\,R^r{}_{mij}   ) ) \cr
&=&  \tfrac12\,\omega^{ij}\,S^a_{ip}\,(S^p_{nm;j} +R^p{}_{mnj} )-\tfrac14\,\omega^{ij}\,
([\nabla_i,\nabla_j] S)^a_{nm}   \cr
&=&  \tfrac12\,\omega^{ij}\,S^a_{ip}\,(S^p_{nm;j} +R^p{}_{mnj} )-\tfrac12\,\omega^{ij}\,
\nabla_i(\nabla_j( S))^a_{nm}\ .
\end{eqnarray*}
From this we get
\begin{eqnarray*}
K_{npm} &=&  \tfrac12\,g_{an}\,\omega^{ij}\,S^a_{is}\,(S^s_{pm;j} +R^s{}_{mpj} )-\tfrac12\,g_{nr}\,\omega^{ij}\, \nabla_i(\nabla_j( S))^r_{pm}\ .
\end{eqnarray*}
From Proposition~\ref{hermetriccompat} the condition for $\nabla_1=\nabla_{QS}+\lambda K$ to be Hermitian-metric compatible is the following, where $\hat;$ denotes Levi-Civita derivative
\begin{eqnarray*}
K_{npm}-K_{mpn} = -\tfrac12\, \CR_{nm\hat;p}+ \tfrac12\, \omega^{ij}\, ( g_{rm}\,\nabla_i(\nabla_j(S))^r_{pn}-g_{nr}\,\nabla_i(\nabla_j(S))^r_{pm})\ ,
\end{eqnarray*}
so on substituting for $K_{npm}$ we find the single condition 
\[  \widehat\nabla \CR =-\,\omega^{ij}\,g_{rs}\,S^s_{jn}(R^r{}_{mki}+S^r_{km;i})\,\extd x^k\tens\extd x^m \wedge \extd x^n  \ .\] 
This is the same as the condition for existence of a fully metric compatible torsion free connection of our assumed form in Theorem~\ref{qlevi}. So, if such a connection exists, our star-preserving one gives it. 
The converse direction is also proved, but obvious (if our star-preserving connection is Hermitian-metric compatible then it is also straight metric compatible and hence the stated condition must hold by Theorem~\ref{qlevi}.)  \endproof

\section{Quantized Surfaces and K\"ahler-Einstein manifolds}

We have seen that our theory applies in particular to any Riemannian manifold equipped with a covariantly constant Poisson-bivector, with the choice $\nabla=\nabla_{LC}$. We then always have a quantum differential algebra by Theorem~\ref{dga} and Corollary~\ref{leviequal} says that the nicest case is when the $\omega$-contracted Ricci tensor  is covariantly constant. In this case we have a quantum symmetric $g_1$ and a quantum-Levi-Civita connection for it. 

\begin{proposition} In the case of a K\"ahler manifold, $\CR$  in Corollary~\ref{leviequal} is  the Ricci 2-form. A sufficient condition for this to be covariantly constant is for the metric to be K\"ahler-Einstein. \end{proposition}
\proof Here $\omega^{ij}=-g^{ik}J_k{}^j=J_k{}^i g^{kj}$ where $J^2=-\id$ and $\CR={1\over 2}\CR_{nm}\extd x^m\wedge\extd x^n$ in our conventions so in Corollary~\ref{leviequal} we have $\CR_{nm}=\omega^{ji}R_{inmj}=g_{kj}\omega^{ji}R_{inm}{}^k=-J_j{}^iR_{inm}{}^j$. Now we use standard complexified local coordinates $z^a,\bar z^{a}$ in which $J_a{}^b=\imath\delta_a{}^b$ and $J_{\bar a}{}^{\bar b}=-\imath\delta_{\bar a}{}^{\bar b}$. The only nonzero elements of Riemann are then of the form
\[ R_{\bar a bc}{}^d=-R_{b\bar a c}{}^d,\quad R_{a\bar b\bar c}{}^{\bar d}=-R_{\bar b a\bar c}{}^{\bar d}.\]
Hence $\CR_{\bar n m}=-\imath R_{a\bar n m}{}^a=\imath R_{\bar n am}{}^a=\imath {\rm Ricci}_{\bar n m}$ and similarly $\CR_{n\bar m}=-\imath {\rm Ricci}_{n\bar m}=-\CR_{\bar m n}$ by symmetry of {\rm Ricci}. Then $\CR_{ij}=-J_i{}^k{\rm Ricci}_{kj}$ in our conventions for 2-form components. Equivalently, $\CR={1\over 2}\CR_{a\bar b}\extd z^{\bar b}\wedge\extd z^a+{1\over 2}\CR_{\bar b a}\extd z^a\wedge\extd\bar z^b=\imath{\rm Ricci}_{a\bar b}\extd z^a\wedge\extd \bar z^b$ as usual.  Clearly in the K\"ahler-Einstein case we have also that ${\rm Ricci}=\alpha g$ for some constant $\alpha$. Then $\CR_{ij}=-J_i{}^k\alpha g_{kj}=-\alpha\omega_{ij}$ in terms of the inverse $\omega_{ij}$ of the Poisson tensor, or $\CR=\alpha\omega/2$ in terms of the symplectic 2-form $\omega=\omega_{ij}\extd x^i\wedge\extd x^j$. This is covariantly constant by our assumption of Poisson-compatibility by Lemma~\ref{compatT}. \endproof

Note that the Ricci 2-form here is closed and represents the 1st Chern class. It is known that every K\"ahler manifold with ${\rm c}_1\le 0$ admits a K\"ahler-Einstein metric and that this is also true under certain stability conditions for positive values. This includes Calabi-Yau manifolds (admitting a Ricci flat metric) and  $\C P^n$ with its Fubini-study metric. Also note that on a K\"ahler manifold the $J$ is also covariantly constant and we may hope to have a noncommutative complex structure in the sense of   \cite{ebsmNCcomplex} to order $\lambda$. This will be considered elsewhere.

Any orientable surface can be given the structure of a  K\"ahler manifold so that the above applies. In fact we do not make use above of the full K\"ahler structure and in the case of an orientable surface we can consider any metric and Poisson tensor $\omega=-{\rm Vol}^{-1}$ as obtained from the volume form, which will be covariantly constant. The generalised Ricci 2-form is then a constant multiple of $S {\rm Vol}$ where $S$ is the Ricci scalar (this follows from the Ricci tensor being $g S/2$ for any surface). So $\CR$ will be covariantly constant if and only if $S$ is constant, i.e. the  case of constant curvature.  

Some general formulae for any surface are as follows, in local coordinates $(x,y)$. Here ${\rm Vol}=\sqrt{\det(g)}\extd x\extd y$ where $g=(g_{ij})$ is the metric. The Poisson tensor $\omega=-{\rm Vol}^{-1}$ is then
\[ \omega=w \left({\del\over\del x}\tens {\del\over\del y}-{\del\over\del y}\tens {\del\over\del x}\right);\quad \omega^{12}= w(x,y):={1\over \sqrt{\det(g)}}\]
which of course gives our product as $x\bullet x=x^2$, $y\bullet y=y^2$, $x\bullet y=xy+{\lambda\over 2} w$, $y\bullet x=xy-{\lambda\over 2}w$, or commutation relations
$[x,y]_\bullet=\lambda w$ on the generators. Similarly, the bimodule commutation relations from the form of $\omega$ are
\[ [f,\xi]_\bullet=\lambda w\left({\del f\over \del x}\nabla_y-{\del f\over\del y}\nabla_x\right)\xi\]
where $\nabla_x,\nabla_y$ are the covariant derivatives along $\del\over\del x$ and $\del\over\del y$ respectively. In terms of $\Gamma$ we have
\[ [f,\extd x^j]_\bullet=\lambda w (f_{,2}\Gamma^j{}_{1m}\extd x^m- f_{,1}\Gamma^j{}_{2m}\extd x^m)\]
or on generators and with $\eps^{12}=1$ antisymmetric,
\[ [x^i,\extd x^j]_\bullet=-\lambda w \eps^{in}\Gamma^j_{nm}\extd x^m.\]
There are similar expressions for $\bullet$ itself in terms of the the classical product plus half of the relevant commutator.

Next, ${\rm Ricci}={S\over 2}g$ implies by symmetries of the Riemann tensor that
\[ R_{1212}={S\over 2}\det(g)=:\rho(x,y)\]
say, with other components determined by its symmetries. In this case 
\[ \CR_{12}=-\CR_{21}=-\omega^{is}R_{i12s}=w\rho,\quad \CR=-w\rho \extd x\extd y=- {S\over 2} {\rm Vol}\]
and 
\[ H^{ij}=-{1\over 2} \omega^{is}R^i{}_{nms}\extd x^m \extd x^n\]
which we compute first with $j$ lowered by the metric as
\[ H^1{}_1=H^2{}_2=-{w\rho\over 2} \extd x\extd y=-{S\over 4}{\rm Vol},\quad H^1{}_2=H^2{}_1=0\]
so we conclude in terms of the inverse metric that 
\[ H^{ij}=-{S\over 4}g^{ij}{\rm Vol}.\]
By Theorem~\ref{dga} we necessarily have a differential graded algebra to order $\lambda$. Here Proposition~\ref{wedge1leib}  in our case becomes
\[ \extd x^i\bullet \extd x^j=\extd x^i\wedge\extd x^j+\tfrac{\lambda}{2} w\left(\Gamma^i_{11}\Gamma^j_{22}-2\Gamma^i_{12}\Gamma^j_{12}+\Gamma^i_{22}\Gamma^j_{11}\right)\extd x\wedge\extd y-\tfrac{\lambda S}{4}g^{ij}{\rm Vol}\]
so that the anticommutation relations for the quantum wedge product have the form
\[ \{\extd x^i, \extd x^j\}_\bullet={\lambda} \left(w^2\left(\Gamma^i_{11}\Gamma^j_{22}-2\Gamma^i_{12}\Gamma^j_{12}+\Gamma^i_{22}\Gamma^j_{11}\right)-\tfrac{S}{2}g^{ij}\right){\rm Vol}.\]

Finally the quantized metric, from (\ref{gQ}) and since $\nabla g=0$,
\begin{eqnarray}\label{gQsurf}g_1&=&g_Q+{\lambda}\CR_{12}\widetilde{\rm Vol}=\tilde g+{\lambda w\over 2}\eps^{ij}g_{ma}\Gamma^a_{ib}\Gamma^b_{jn}\extd x^m\tens_1\extd x^n+{\lambda S\over 2}\widetilde{\rm Vol}\end{eqnarray}
where the first two terms are $g_Q$ and  
\[ \tilde g:=g_{ij}\extd x^i\tens_1\extd x^j,\quad \widetilde{\rm Vol}:={1\over 2 w}(\extd x\tens_1\extd y-\extd y\tens_1\extd x)\]
are shorthand notations. Similarly, the connection $\nabla_Q$ is computed from the local formula (\ref{nablaQ}).  As explained, we will have $\nabla_Q g_1=0$ at order $\lambda$ if and only if  $S$ is constant. We compute  further details for the two basic examples.

\subsection{Quantised hyperbolic space}

As the basic example we look at the Poincar\'e upper half plane with its hyperbolic metric
\[ M=\{ (x,y)\in \R^2\ |\ y>0\},\quad g={1\over y^2}(\extd x\tens\extd x+ \extd y\tens\extd y)\]
which is readily found to have nonzero Christoffel symbols
\[  \Gamma^1_{12}=\Gamma^1_{21}=\Gamma^2_{22}=-y^{-1}\ ,\quad \Gamma^2_{11}=y^{-1}\]
or $\Gamma^i_{1j}=-\eps^{ij}y^{-1}$ and $\Gamma^i_{2j}=-\delta^i_jy^{-1}$.  The bivector $\omega^{12}=-\omega^{21}= y^2$
is easily seen to be the unique solution to (\ref{poissoncomp}) up to normalisation. This is the inverse of the volume form ${\rm Vol}=y^{-2}\extd x\extd y$. 

Clearly from the Poisson tensor
\[ \omega=y^2\left({\del\over\del x}\tens {\del\over\del y}-{\del\over\del y}\tens {\del\over\del x}\right)\]
we have $ [x, y]_\bullet={\lambda}y^2$,
which relations also occur for the standard bicrossproduct model spacetime in 2-dimensions in terms of  inverted coordinates in \cite{BegMa4}. Also note that $[x,y^{-1}]_\bullet=\lambda$. Note that although the relations do extend to  an obvious associative algebra $A_\lambda$, this is not unique and not immediately  relevant.

Next, from $\Gamma$ we see that
\[ [f,\extd x]_\bullet=\lambda y\left({\del f\over\del x}\extd x- {\del f\over\del y}\extd y\right),\quad [f,\extd y]_\bullet=\lambda y\left({\del f\over\del y}\extd x+ {\del f\over\del x}\extd y\right)\]
or on generators we have
\[ [x,\extd x]_\bullet=[y,\extd y]_\bullet=\lambda y\extd x,\quad [x,\extd y]_\bullet=-[y,\extd x]_\bullet=\lambda y\extd y.\]
There are similar expressions for $\bullet$ itself in terms of the the classical product.

The Ricci scalar here is $S=-2$ so 
\[ \CR={\rm Vol},\quad H^{ij}={1\over 2}g^{ij}{\rm Vol}\]
and from the latter  we obtain
\[ \extd x^i\bullet \extd x^j=\extd x^i\wedge\extd x^j+\tfrac{\lambda}{2} y^2\left(\Gamma^i_{11}\Gamma^j_{22}-2\Gamma^i_{12}\Gamma^j_{12}+\Gamma^i_{22}\Gamma^j_{11}\right)\extd x\wedge\extd y+\tfrac{\lambda}{2}\delta_{ij}\extd x\wedge\extd y\]
which from the form of $\Gamma$ simplifies further to
\[ \extd x^i\bullet \extd x^j=\extd x^i\wedge\extd x^j-\tfrac{\lambda}{2}\delta_{ij}\extd x\wedge\extd y,\quad \{\extd x^i, \extd x^j\}_\bullet=-{\lambda}\delta_{ij}\extd x\wedge\extd y\]
which has a `Clifford algebra-like' form. The result here is the same as obtained by applying $\extd$ to the bimodule relations, i.e. is consistent with the
maximal prolongation of the first order calculus.

Finally, we have  our constructions of noncommutative Riemannian geometry. In our case
\[ \eps^{ij}g_{ma}\Gamma^a_{ib}\Gamma^b_{jn}=0\]
so that $g_Q$ has he same form as classically but with $\tens_1$ and 
\[ g_1={\extd x^i\over y^2}\tens_1\extd x^i-{\lambda}\widetilde{\rm Vol}\]
(sum over $i$). Similarly, one may compute using the form of $\Gamma$ in (\ref{nablaQ}) that
\[ \nabla_Q\extd x^i={\extd y\over y}\tens_1\extd x^i+{\extd x\over y}\tens_1\eps_{ij}\extd x^j \]
which again has the same form as classically.  There is an associated generalised brading $\sigma_Q$ making this a bimodule connection. As per our general theory, $\nabla_Q$ is quantum torsion free and metric compatible with  $g_1$.

All constructions above are invariant under $SL_2(\R)$ and hence under the modular group and other discrete subgroups. Indeed, the metric is well known to be invariant. The volume form can also easily be seen to be and correspondingly $\omega$ is invariant. As these are the
only inputs into the theory it follows that the deformed structures are likewise compatible with this action. The quotient of the constructions corresponds to
replacing the Poincar\'e upper half plane by a Riemann surface of constant  negative curvature, constructed as quotient.  Therefore this is achieved in principle. We might reasonably then expect a role for modular forms in the deeper aspects of the noncommutative geometry.

\subsection{Quantised  sphere}

The case of a surface of constant positive curvature, the sphere, is the $n=1$ case of $\C P^n$ which will be covered elsewhere  in holomorphic coordinates.  Here we give it is as an example of the analysis for surfaces above. 

We work in the upper hemisphere in standard cartesian coordinates, with similar formulae for the lower hemisphere, so 
\[ M=\{(x,y)\ |\ x^2+y^2<1\},\quad z=\sqrt{1-x^2-y^2},\]
\[g={1\over z^2}\left((1-y^2)\extd x\tens\extd x+xy(\extd x\tens\extd y+\extd y\tens\extd x)+(1-x^2)\extd y\tens\extd y\right)\]
which is readily found to have symmetric Christoffel symbols
\[ \Gamma^1_{11}={x\over z^2}(1-y^2),\quad \Gamma^1_{22}={x\over z^2}(1-x^2),\quad \Gamma^1_{12}={x^2 y\over z^2}\]
\[ \Gamma^2_{11}={y\over z^2}(1-y^2),\quad \Gamma^2_{22}={y\over z^2}(1-x^2),\quad \Gamma^2_{12}={x y^2\over z^2}.\]
or compactly $\Gamma^i_{jk}=x^i g_{jk}$. 

The inverse of the volume form ${\rm Vol}=z^{-1}\extd x\extd y$ gives the Poisson bivector
\[ \omega=z({\del\over\del x}\tens{\del\over\del y}-{\del\over\del y}\tens{\del\over\del x})\]
so we have relations 
\[ [x,y]_\bullet=\lambda z,\quad [z,x]_\bullet=\lambda y,\quad [y,z]_\bullet=\lambda x,\]
the standard relations of the fuzzy sphere. In this case there is an associative quantisation to all orders as the enveloping algebra $U(su_2)$
modulo a constant value of the quadratic Casimir. It is known that this algebra does not admit an associative 3D rotationally invariant calculus\cite{BegMa:coch}
so there won't be a zero-curvature Poisson-compatible preconnection. At present we use the Levi-Civita connection according to Corollary~\ref{leviequal}. Then from $\Gamma$ we have 
\[ [f,\extd x^j]_\bullet= -\lambda z\, x^j f_{,i}\eps^{ik}g_{km}\extd x^m \]
for the bimodule relations of the quantum differential calculus, where $\eps^{12}=1$ is antisymmetric. Explicitly,
\[ [x,\extd x^i]_\bullet=-\lambda{x^i\over z}(xy\extd x+(1-x^2)\extd y),\quad  [y,\extd x^i]_\bullet=\lambda{x^i\over z}((1-y^2)\extd x+ xy\extd y).\]
Next, the Ricci scalar of the unit sphere is $S=2$ so 
\[  \CR=-{\rm Vol},\quad  H^{ij}=-{1\over 2}g^{ij}{\rm Vol}\ ,\quad \mathrm{where}\
g^{ij}=\begin{pmatrix}1-x^2 & -xy \cr -xy & 1-y^2\end{pmatrix}\ .
\]
From this and $\Gamma$ we obtain
\[ \extd x^i \wedge_\bullet \extd x^j=\extd x^i\wedge\extd x^j+ \lambda\left(x^ix^j-{1\over 2}g^{ij}\right){\rm Vol},\quad \{\extd x^i, \extd x^j\}_\bullet= \lambda\left(2 x^ix^j-g^{ij}\right){\rm Vol}\]
for the exterior algebra relations. One can verify that this is the maximal prolongation of the bidmodule relations.

Finally, we note that  from the form of $g$ that  $x^ag_{ai}=x^i z^{-2}$ and $\omega^{ab}g_{ai}g_{bj}=\eps_{ij} z^{-1}$. The first of these and the form of $\Gamma$ gives us
\[ g_{ij,k}=\Gamma^a_{ki}g_{aj}+\Gamma^a_{kj}g_{ia}=x^a(g_{ki}g_{aj}+g_{kj}g_{ia})=z^{-2}(x^jg_{ki}+x^ig_{kj})\]
using metric compatibility. One may then compute the quantum metric and connection from (\ref{gQsurf}) and (\ref{nablaQ}) respectively as
\[ g_1=\tilde g+{\lambda\over 2 z^3}x^m\extd x^m\tens_1  x^a\eps_{an}\extd x^n+{\lambda}\widetilde{\rm Vol}\]
\begin{eqnarray*} \nabla_Q\extd x^i&=&-x^i\tilde g-{\lambda}x^i\widetilde{\rm Vol}-{\lambda\over 2 z}x^m\extd x^m\tens_1( \eps^{ib}g_{bn}+ {x^i x^b\over z^2}\eps_{bn})\extd x^n\\
&=&-x^i g_1-{\lambda\over 2 z}x^m\extd x^m\tens_1 \eps^{ib}g_{bn}\extd x^n=-x^i\bullet g_1\end{eqnarray*}
(sum over $m$). Here on the left  $x^i\tilde g$ is a shorthand notation for the previously defined element of $\Omega^1A_1\tens_1\Omega^1A_1$ but now with an extra classical $x^i$ in the definition. One can think if it as made with the classical product when the classical and quantum vector spaces are identified, and ditto for $x^i g_1$. The expression $x^i\bullet g_1$ is computed with the quantised product but only on the first tensor factor of $g_1\in\Omega^1A_1\tens_1\Omega^1A_1$ (since this is the relevant bimodule structure). 

In fact, this example  at the level of first order differentials was proposed in \cite{BegMa:coch} where we showed that the Levi-Civita connection arises as a cochain twist of the classical exterior algebra by a certain action of the Lorenz group. This could potentially be used to construct the full noncommutative nonassociative Riemannian geometry by twisting\cite{BegMa:twi}.

\section{$\nabla$ far from Levi-Civita: bicrossproduct and black-hole models} \label{bcwkkuhvcj}

In this section we give two contrasting examples where  where we cannot take  $\nabla=\widehat\nabla$. The black-hole for a natural rotationally invariant Poisson bracket will provide an example where some of the obstructions in the general theory hold, i.e. we cannot find a quantising connection $\nabla$ (expressed by choice of its torsion $T)$ such that the quantum `Levi-Civita' connection of the form $\nabla_1=\nabla_{QS}+\lambda K$ is star preserving, torsion free and metric compatible at the same time; one or more of these features necessarily gets an  order $\lambda$ correction. Here $\nabla_S=\widehat{\nabla}$ is the classical Levi-Civita. On the other hand, we will find that the generalised Ricci 2-form $\CR=0$  so that the quantum metric $g_1=g_Q$, the functorial one. We'll find in the black-hole case that $\nabla$ necessarily has curvature and hence the quantum differential calculus will be nonassociative. 

Before doing that we give an easier warm-up example which also illustrates all our semiclassical theory and where the algebraic version is already exactly solved by computer algebra\cite{BegMa4}. In this 2D example all the obstructions vanish and there is a unique quantum $\nabla_1$ that is star-preserving, torsion free and metric compatible. Here the existing differential calculus, derived from the theory of quantum groups, gives $\omega,\nabla$ while  $\nabla g=0$ then forces the metric. Moreover,  $\nabla$ has torsion but no curvature and yet the 2-form $\CR\ne 0$,  in contrast to the Schwarzschild black-hole case. On the other hand this 2D model still has a physical interpretation as a toy model with strong gravitational source, so strong that even light can't escape (so something like the inside of a black hole but with decaying rather than zero Ricci tensor).  We refer to \cite{BegMa4} for details and for a different, cosmological, interpretation as well.

\subsection{The 2D bicrossproduct model}

Setting $x^0=t$ and $x^1=r$, we have $\omega^{10}=-\omega^{01}=r$ as the semiclassical data behind the bicrossproduct model commutation relations $[t,r]_\bullet=\lambda r$. It is known that this model has a standard 2D differential calculus with nonzero relations
\[ [r,\extd t]_\bullet=\lambda\extd r,\quad [t,\extd t]_\bullet=\lambda\extd t,\]
which has as its underlying semiclassical data a connection with Christoffel symbols  $\Gamma^0_{01}=-r^{-1}$ and $\Gamma^0_{10}=r^{-1}$ and all other Christoffel symbols zero. This has torsion $T^0_{10}=-T^0_{01} =2\,r^{-1}$ and $T^1_{ij}=0$ and one
 can check that it Poisson-compatible. One can calculate
\begin{eqnarray*}
T^i_{01;p} &=& T^i_{01,p} +\Gamma^i_{pn}\,T^n_{01}-\Gamma^n_{p0}\,T^i_{n1}-\Gamma^n_{p1}\,T^i_{0n}  \cr
 &=& \delta_{0i}\,(T^0_{01,p} +\Gamma^0_{pn}\,T^n_{01}-\Gamma^n_{p0}\,T^0_{n1}-\Gamma^n_{p1}\,T^0_{0n})  \cr
&=& \delta_{0i}\,(T^0_{01,p} +\Gamma^0_{p0}\,T^0_{01}-\Gamma^0_{p0}\,T^0_{01})  = \delta_{0i}\,\delta_{1p}\,T^0_{01,1}  = 2\,r^{-2}\, \delta_{0i}\,\delta_{1p}\ 
\end{eqnarray*}
and that the curvature is zero, as it should since the standard calculus is associative to all orders. To see this, without
loss of generality, we look at $j=0$, $k=1$:
\begin{eqnarray*}
R^{l}_{\phantom{l}i01} &=& \frac{\partial \Gamma^l_{1i}}{\partial
x^0}\,-\, \frac{\partial \Gamma^l_{0i}}{\partial x^1}\,+\,
\Gamma^m_{1i}\,\Gamma^l_{0m}\,-\,\Gamma^m_{0i}\,\Gamma^l_{1m}\cr
&=& \delta_{0l}\,( \frac{\partial \Gamma^0_{1i}}{\partial
t}\,-\, \frac{\partial \Gamma^0_{0i}}{\partial r}\,+\,
\Gamma^0_{1i}\,\Gamma^0_{00}\,-\,\Gamma^0_{0i}\,\Gamma^0_{10})   \cr
&=& \delta_{0l}\,( -\, \frac{\partial \Gamma^0_{0i}}{\partial t}\,-\,\Gamma^0_{0i}\,\Gamma^0_{10})   
= \delta_{0l}\,\delta_{1i}\,( -\, \frac{\partial \Gamma^0_{01}}{\partial r}\,-\,\Gamma^0_{01}\,\Gamma^0_{10})   =0\ .
\end{eqnarray*}
Next we compute,
\begin{eqnarray*}
H^{ij}&:=&
\tfrac14\omega^{is}\big(T^j_{nm;s}
-2 R^j_{nms}  \big)\,\extd x^m\wedge\extd x^n  \cr
&=&  \tfrac14\,\delta_{0j}\,\omega^{is}T^0_{nm;1}\,\delta_{1s}
\,\extd x^m\wedge\extd x^n  
=  \tfrac14\,\delta_{0j}\,\omega^{i1}T^0_{nm;1}
\,\extd x^m\wedge\extd x^n  \cr
&=&  \tfrac14\,\delta_{0j}\,\delta_{0i}\,\omega^{01}T^0_{nm;1}
\,\extd x^m\wedge\extd x^n  \cr
&=&  \tfrac14\,\delta_{0j}\,\delta_{0i}\,\omega^{01} (T^0_{01;1}
\,\extd x^1\wedge\extd x^0 + T^0_{10;1}\,\extd x^0\wedge\extd x^1) \cr
&=&  \tfrac14\,\delta_{0j}\,\delta_{0i}\,\omega^{01}\big(   T^0_{01;1}
- T^0_{10;1}\big)\,\extd x^1\wedge\extd x^0  \cr
&=&  \tfrac14\,\delta_{0j}\,\delta_{0i}\,\omega^{01} 2\, T^0_{01;1}
\,\extd r\wedge\extd t  \cr
&=&  \tfrac12\,\delta_{0j}\,\delta_{0i}\,(-r)   2\,r^{-2}
\,\extd r\wedge\extd t  
=\delta_{0j}\,\delta_{0i}\,r^{-1}\,\extd t\wedge\extd r  \ .
\end{eqnarray*}
The wedge product obeying the Leibniz rule in Theorem~\ref{dga} is then;
\begin{eqnarray} \label{sbhdfvckjc}
\xi\wedge_1\eta \,&=&\,\xi\wedge\eta+\tfrac{\lambda}{2} \omega^{ij}\,\nabla_i\xi\wedge\nabla_j\eta \nonumber\\ &&
+ (-1)^{|\xi|+1}\, \lambda\,r^{-1} \,\extd t\wedge\extd r\wedge (\partial_0 \, \righthalfcup\, \xi)\wedge(\partial_0 \, \righthalfcup\, \eta) \ .
\end{eqnarray}
For $\xi$ and $\eta$ being either $\extd r$ or $\extd t$, the only potentially deformed case is
\begin{eqnarray*}
\extd t\wedge_1\extd t &=& \tfrac\lambda2\,\omega^{ij}\,\nabla_i(\extd t)\wedge\nabla_j(\extd t)
+ \lambda\,r^{-1} \,\extd t\wedge\extd r\wedge (\partial_0 \, \righthalfcup\, \extd t)\wedge(\partial_0 \, \righthalfcup\, \extd t) \cr
&=& \tfrac\lambda2\,\big(\omega^{01}\,\nabla_0(\extd t)\wedge\nabla_1(\extd t)+\omega^{10}\,\nabla_1(\extd t)\wedge\nabla_0(\extd t)   \big)
+ \lambda\,r^{-1} \,\extd t\wedge\extd r\cr
&=& \tfrac\lambda2\,\omega^{01}\,\big(\nabla_0(\extd t)\wedge\nabla_1(\extd t)-\nabla_1(\extd t)\wedge\nabla_0(\extd t)   \big)
+ \lambda\,r^{-1} \,\extd t\wedge\extd r=0\ .
\end{eqnarray*}
The exterior algebra among these basis elements is therefore undeformed, in agreement with the noncommutative algebraic picture where this is known
(and holds to all orders). 

Our goal is to study the semiclassical geometry of this model using our functorial methods. First of all, the above connection is not compatible with the flat metric, but {\em is} compatible with the metric
 \begin{eqnarray*}
g=g_{ij}\,\extd x^i\tens\extd x^j=b\,r^2\,\extd t\tens\extd t -b\,r\,t\,(\extd t\tens\extd r +\extd r\tens\extd t) +(1+b\,t^2)\,\extd r\tens\extd r\ .
\end{eqnarray*}
where $b$ is a non-zero real parameter. This is our semiclassical analogue of the obstruction discovered in \cite{BegMa4}. For our purposes it is better to write the metric as the following, where $v=r\,\extd t-t\,\extd r$
 \begin{eqnarray*}
g=\extd r\tens\extd r+b\,v\tens v\ .
\end{eqnarray*}
Note that $\nabla$ applied to both $\extd r$ gives zero. We 
quantise the classical bicrossproduct spacetime with this metric. First
\begin{eqnarray*}
q^{-1}(g)
&=& \extd r\tens_1\extd  r+b \, v\tens_1 v\ .
\end{eqnarray*}
From the expression for $H^{ij}$, we have $\CR=g_{ij}H^{ij}=b\,r\,\extd t\wedge\extd r=b\,v\wedge\extd r=\pm \sqrt{|b|}{\rm Vol}$ and 
 according to our  general scheme, we take
\[ g_1=\extd r\tens_1\extd  r+b\, v\tens_1 v+\tfrac{b\lambda}{2}(\extd r\tens_1 v-v\tens_1\extd r).\]
To compare with \cite{BegMa4}, if we let
\begin{equation}\label{2Dcompare} \nu:=r\bullet_\lambda\extd t-t\bullet_\lambda \extd r=v+\tfrac{\lambda}{2}\extd r,\quad \nu^*:=(\extd t)\bullet_\lambda r-( \extd r)\bullet_\lambda t=v-\tfrac{\lambda}{2}\extd r\end{equation}
and identify these with $v,v^*$ in \cite{BegMa4} (apologies for the clash of notation) then the quantum metric there gives the same answer as $g_1$ above, i.e.\ this is the leading order part of the noncommutative geometry. From Theorem~\ref{functor} we get $\nabla_Q$ vanishing on both $\extd r$ and $v$, and for all 1-forms $\xi$, $\sigma_Q(\extd r\tens_1\xi)=\xi\tens_1\extd r$ and $\sigma_Q(v\tens_1\xi)=\xi\tens_1 v$.

Next we express the classical Levi-Civita connection for the above metric in the form $\nabla_S$. We use (\ref{SfromT}) together with the only nonvanishing downstairs torsions being $T_{010}=-T_{001}=2\,b\,r$
and $T_{110}=-T_{101}=-\,2\,b\,t$ and
\begin{eqnarray*}
S^a_{bc} &=& \tfrac12 g^{a0}(T_{0bc}-T_{bc0}-T_{cb0})+\tfrac12 g^{a1}(T_{1bc}-T_{bc1}-T_{cb1})  \ ,
\end{eqnarray*}
to give 
\begin{eqnarray*}
S^a_{11} &=& -\, g^{a0}\,T_{110}\ ,\ S^a_{00}=-\, g^{a1}\,T_{001}\ ,\ S^a_{10} =  g^{a1}\,T_{110} \ ,\ 
S^a_{01} =  g^{a0}\, T_{001}\  .
\end{eqnarray*}
The upstairs metric is given by
\begin{eqnarray*}
g^{00}=(1+b\,t^2)/(b\,r^2)\ ,\ g^{01}=g^{10}=t\,r^{-1}\ ,\ g^{11}=1\ .
\end{eqnarray*}
\begin{eqnarray*}
S^a_{11} &=& 2\,b\,t\, g^{a0}\ ,\ S^a_{00}=2\,b\,r\, g^{a1}\ ,\ S^a_{10} =  -\,2\,b\,t\,g^{a1} \ ,\ 
S^a_{01} =  -\,2\,b\,r\,g^{a0},\end{eqnarray*}
which we write compactly, along with its covariant derivative, as
\[  S^a_{ij}=2b \eps_{im} x^m\eps_{jn}g^{an},\quad S^a_{ij;k}=2b\eps_{ik}\eps_{jm}g^{am}\]
where $\eps_{01}=1$ is antisymmetric and $g$ is preserved by $\nabla$ (as expressed by $;i$). We also have $\widehat{\nabla}\CR=0$ since $\CR$ was a multiple of the volume form,  and $R=0$ for the curvature of $\nabla$, so the obstruction in Theorem~\ref{qlevi} for a torsion free metric compatible quantum connection is
\begin{eqnarray*}  \widehat\nabla \CR &+&\omega^{ij}\,g_{rs}\,S^s_{jn}(R^r{}_{mki}+S^r_{km;i})\,\extd x^k\tens\extd x^m \wedge \extd x^n  \\
&& =\omega^{ij}\,g_{rs}\,S^s_{jn}S^r_{km;i}\,\extd x^k\tens\extd x^m \wedge \extd x^n=0
\end{eqnarray*}
when we put in the compact form of $S$ and its covariant derivative. Hence Theorem~\ref{qlevi} tells us that there is a unique such quantum connection of the form $\nabla_1=\nabla_{QS}+\lambda K$. Corollary~\ref{jointcompat} tells us that this is also the unique star-preserving connection of this form. In short, all obstructions vanish and we have a unique quantum Levi-Civita connection with all our desired properties.

It only remains to compute $\nabla_1$.  We take the liberty of changing the basis to write for $K$ real
\begin{eqnarray*}
K(v)=K^v_{vv}\,v\tens v+K^v_{rv}\,\extd r\tens v+K^v_{vr}\,v\tens \extd r+K^v_{rr}\,\extd r\tens \extd r\ ,\cr
K(\extd r)=K^r_{vv}\,v\tens v+K^r_{rv}\,\extd r\tens v+K^r_{vr}\,v\tens \extd r+K^r_{rr}\,\extd r\tens \extd r\ .
\end{eqnarray*}

\begin{proposition} \label{loivcyu}
The unique star-preserving quantum connection  of the form $\nabla_1=\nabla_{QS}+\lambda K$ is also torsion free and metric compatible (`quantum Levi-Civita') and given by non-zero components
\begin{eqnarray*}
 K^r_{vr}=K^v_{vv}=-2\,b\,r^{-1}\ 
\end{eqnarray*}
in our basis, leading to
\[ \nabla_{1} \extd r={2 b v\over r}\tens_1 v-{2 b\lambda\over r}v\tens_1 \extd r,\quad \nabla_{1} v= -{2 v\over r}\tens_1\extd r-{2 b \lambda\over r} v\tens_1 v\ .\]
\end{proposition}
\proof  
Note that $v^*=v$ and $\extd r^*=\extd r$ and also that Theorem~\ref{starpresK} tells us the value of $K$ which can be computed out as the value stated. But we still need to compute $\nabla_1$ and, moreover, since this is an illustrative example we will also verify its properties directly as a nontrivial check of all our main theorems.   

First we compute $S$ as an operator from the components stated above (or one can readily compute the classical Levi-Civita connection and find $S$ as the difference between this and $\nabla$). Either way,
\[ 
S(\extd r) = 2\,b\, r^{-1}\, v\tens v,\quad  S(\extd t) = 2\,b\, t\,r^{-2}\,v\tens v -  2\,r^{-2}\,v\tens\extd r,\quad 
S(v) = -  2\,r^{-1}\,v\tens\extd r \]

Next we compute $\nabla_{QS}$ and its associated generalised braiding. In the following calculation, $\nabla_0,\nabla_1$ denote the components $\nabla_i$ of the classical connection $\nabla$ (apologies for the clash of notation). We have $\nabla_0(S)=0$ and 
\begin{eqnarray*}
\nabla_1(S)(v)&=&\nabla_1(S(v))=\nabla_1(-  2\,r^{-1}\,v\tens\extd r)=2\,r^{-2}\,v\tens\extd r\ ,\cr
\nabla_1(S)(\extd r)&=&\nabla_1(S(\extd r))=\nabla_1(2\,b\, r^{-1}\, v\tens v)=-\,2\,b\, r^{-2}\, v\tens v\ .
\end{eqnarray*}
From Proposition~\ref{quantQS},
\begin{eqnarray*}
\sigma_{QS}(v\tens_1\xi) &=& \sigma_{Q}(v\tens_1\xi)+\lambda\,\omega^{01}\,\xi_{0}\,\nabla_1(S)(v)  \cr
&=& \xi \tens_1 v-\lambda\,r\,\xi_{0}\,\nabla_1(S)(v)  \cr
&=& \xi \tens_1 v-2\,\lambda\,\xi_{0}\,r^{-1}\,v\tens\extd r\ ,  \cr
\sigma_{QS}( \extd r \tens_1\xi) &=& \xi \tens_1 \extd r - \lambda\,r\,\xi_{0}\,\nabla_1(S)(\extd r)  \cr
&=& \xi \tens_1 \extd r +2\, \lambda\,\xi_{0}\,b\, r^{-1}\, v\tens v  \ .
\end{eqnarray*}
and 
\begin{eqnarray*}
Q(S)(v) &=& S(v)+\tfrac\lambda2\,\omega^{ij}\,\nabla_i(\nabla_j(S)(v)) 
= S(v)+\tfrac\lambda2\,\omega^{01}\,\nabla_0(\nabla_1(S)(v)) \cr
&=& S(v) = -  2\,r^{-1}\,v\tens\extd r\ ,\cr
Q(S)(\extd r) &=& S(\extd r)+\tfrac\lambda2\,\omega^{ij}\,\nabla_i(\nabla_j(S)(\extd r)) 
= S(\extd r)+\tfrac\lambda2\,\omega^{01}\,\nabla_0(\nabla_1(S)(\extd r)) \cr
&=& S(\extd r)=2\,b\, r^{-1}\, v\tens v\ .
\end{eqnarray*}
Then
\begin{eqnarray*}
\nabla_{QS}(v) &=& \nabla_Q(v)+q^{-1}Q(S)(v)=-  2\,q^{-1}(r^{-1}\,v\tens\extd r)= -  2\,r^{-1}\,v\tens_1\extd r\ ,\cr
\nabla_{QS}(\extd r) &=& \nabla_Q(\extd r)+q^{-1}Q(S)(\extd r)=2\,b\, q^{-1}(r^{-1}\, v\tens v) =
2\,b\, r^{-1}\, v\tens_1 v\ .
\end{eqnarray*}
We can add this to the $K$ obtained from Theorem~\ref{starpresK} to obtain the result stated for the quantum Levi-Civita connection. 

For illustrative purposes  let us also see directly why this adjustment is necessary and that it succeeds. Firstly,  
to be star preserving we need $(\id\tens\star)\nabla_{QS}(\xi)=(\star^{-1}\tens\id)\Upsilon\,\overline{\sigma_{QS}^{-1}\nabla_{QS}(\xi^*)}$ for our two cases, $\xi=v$ and $\xi=\extd r$. It is more convenient to rearrange this as 
\begin{eqnarray*}
\overline{\nabla_{QS}(\xi^*)}=\overline{\sigma_{QS}}\,\Upsilon^{-1}(\star\tens\star)\nabla_{QS}(\xi)\ .
\end{eqnarray*}
We do this for the two cases, where $b$ is real
\begin{eqnarray*}
\overline{\sigma_{QS}}\,\Upsilon^{-1}(\star\tens\star)\nabla_{QS}(v) &=& 
\overline{\sigma_{QS}}\,\Upsilon^{-1}(\star\tens\star)(-  2\,r^{-1}\,v\tens_1\extd r)\cr
&=& - \,  2\,\overline{\sigma_{QS}}\,\Upsilon^{-1}(\overline{r^{-1}\,v}\tens_1\overline{\extd r})  \cr
&=& - \,  2\,\overline{\sigma_{QS}  (  \extd r  \tens_1  r^{-1}\,v)   }  \cr
&=& - \,  2\,(\overline{  r^{-1}\,v\tens_1 \extd r +2\, \lambda\,b\, r^{-1}\, v\tens v  } )\ ,\cr
\overline{\sigma_{QS}}\,\Upsilon^{-1}(\star\tens\star)\nabla_{QS}(\extd r) &=& 
\overline{\sigma_{QS}}\,\Upsilon^{-1}(\star\tens\star)(2\,b\, r^{-1}\, v\tens_1 v) \cr
&=& 2\,b\, \overline{\sigma_{QS}}\,\Upsilon^{-1}(\overline{r^{-1}\, v}\tens_1 \overline{v}) \cr
&=& 2\,b\, \overline{\sigma_{QS}  ( v \tens_1  r^{-1}\, v  ) } \cr
&=& 2\,b\, (\overline{r^{-1}\, v \tens_1 v-2\,\lambda\,r^{-1}\,v\tens\extd r }) \ .
\end{eqnarray*}
The difference in going clockwise minus anticlockwise round the diagram in 
Lemma~\ref{oiuyfvcftyytrs} is now
\begin{eqnarray*}
\overline{\sigma_{QS}}\,\Upsilon^{-1}(\star\tens\star)\nabla_{QS}(v) - \overline{\nabla_{QS}(v)} 
&=& - \,  2\,(\overline{ 2\, \lambda\,b\, r^{-1}\, v\tens v  } )\ ,\cr
\overline{\sigma_{QS}}\,\Upsilon^{-1}(\star\tens\star)\nabla_{QS}(\extd r) - \overline{\nabla_{QS}(\extd r)} 
&=& 2\,b\, (\overline{-2\,\lambda\,r^{-1}\,v\tens\extd r }) \ .
\end{eqnarray*}
Thus $\nabla_{QS}$ is not star preserving. However we follow 
 Theorem~\ref{starpresK} to see that $\nabla_{QS}+\lambda K$ is star preserving  if and only if the only nonzero $K^a_{bc}$ in this basis are
\begin{eqnarray*}
K^v_{vv}=K^r_{vr}=-\,2\,b\,r^{-1}\ .
\end{eqnarray*}
This completes the direct derivation of the values in this example.

Although implied by our theory, let us also see how the quantum torsion and metric compatibility conditions get to hold. We calculate the torsions from (\ref{sbhdfvckjc}) by 
\begin{eqnarray*}
T_{\nabla_{QS}}(v) &=& -  2\,r^{-1}\,v\wedge_1\extd r-\extd v=-  2\,r^{-1}\,v\wedge\extd r-\extd v =0\ ,  \cr
T_{\nabla_{QS}}(\extd r) &=&  2\,b\, r^{-1}\, v\wedge_1 v
=    2\,b\,   \lambda\,\extd t\wedge\extd r\ .
\end{eqnarray*}
 The condition for $\nabla_{QS}+\lambda K$ to be torsion free is that
\begin{eqnarray*}
0=\lambda\,(\wedge K(v))\ ,\ 0=\lambda\,(\wedge K(\extd r))+2\,b\,\lambda\,\extd t\wedge\extd r\ ,
\end{eqnarray*}
which becomes 
\begin{eqnarray*}
K^v_{rv}=K^v_{vr}\ ,\ K^r_{vr}-K^r_{rv}+2\,b\,r^{-1}=0\ 
\end{eqnarray*}
which indeed holds for our found values. 

Next,  right connections on the conjugate of the 1-forms are
\begin{eqnarray*}
\bar\nabla_{QS}(\overline{v}) &=& (\id\tens\star^{-1})\Upsilon\, \overline{\nabla_{QS}(v)}
=  -  2\,(\id\tens\star^{-1})\Upsilon\, \overline{r^{-1}\,v\tens_1\extd r} 
= -  2\,\overline{\extd r} \tens_1 r^{-1}\,v\ ,\cr
\bar\nabla_{QS}(\overline{\extd r}) &=& (\id\tens\star^{-1})\Upsilon\, \overline{\nabla_{QS}(\extd r)}
=2\,b\, (\id\tens\star^{-1})\Upsilon\, \overline{r^{-1}\, v\tens_1 v} = 2\,b\, 
\overline{ v} \tens_1 r^{-1}\, v\ .
\end{eqnarray*}
Now we apply $\nabla_{QS}$ to the Hermitian metric $(\star\tens\id)g_1$ as follows:
\begin{eqnarray*}
(\bar\nabla_{QS}\tens\id+\id\tens \nabla_{QS})(\star\tens\id)g_1=0\ .
\end{eqnarray*}
Hence the condition for  $\nabla_1=\nabla_{QS}+\lambda K$ to preserve the Hermitian metric is
\begin{eqnarray*}
0 &=& (\id\tens\star^{-1})\Upsilon\,\overline{\lambda\,K(\extd r)}\tens \extd r
+ b\, (\id\tens\star^{-1})\Upsilon\,\overline{\lambda\,K(v)}\tens v \cr 
&& +\, \overline{\extd r} \tens \lambda\,K(\extd r) + b\, \overline{v} \tens \lambda\,K(v)\ .
\end{eqnarray*}
We split this into two parts, depending on whether we end in $\extd r$ or $v$, to give, as $\lambda$ is imaginary,
\begin{eqnarray*}
0 &=& K^r_{vr} \overline{\extd r}\tens v+K^r_{rr} \overline{\extd r}\tens \extd r+b\, K^v_{vr} \overline{v}\tens v+b\, K^v_{rr} \overline{v}\tens \extd r -(\id\tens\star^{-1})\Upsilon\,\overline{K(\extd r)} \cr
&=& (b\, K^v_{vr}-K^r_{vv}) \overline{v}\tens v+(b\, K^v_{rr}-K^r_{rv}) \overline{v}\tens \extd r\ ,\cr
0 &=& K^r_{vv} \overline{\extd r}\tens v+K^r_{rv} \overline{\extd r}\tens \extd r+b\, K^v_{vv} \overline{v}\tens v+b\, K^v_{rv} \overline{v}\tens \extd r -b\,(\id\tens\star^{-1})\Upsilon\,\overline{K(v)} \cr
&=& (K^r_{vv}-b\,K^v_{vr}) \overline{\extd r}\tens v+(K^r_{rv}-b\,K^v_{rr}) \overline{\extd r}\tens \extd r\ .
\end{eqnarray*}
Thus the conditions for $\nabla_1$ to preserve  the Hermitian quantum metric reduce to $K^r_{vv}=b\,K^v_{vr}$ and $K^r_{rv}=b\,K^v_{rr}$, which again holds for our found values. 

Finally we consider straight quantum metric compatibility $\nabla_1 g_1$. Again, this has to follow from Hermitian-metric compatibiity since $\nabla_1$ is star-preserving, but we check directly what is needed  as per Theorem~\ref{qlevi}. First we apply $\nabla_{QS}$ to $g_1$,
\begin{eqnarray*}
\nabla_{QS}(\extd r \tens_1 \extd r) &=& \nabla_{QS}(\extd r) \tens_1 \extd r +(\sigma_{QS}\tens\id)
(\extd r \tens_1 \nabla_{QS}(\extd r) ) \cr
&=& 2\,b\, r^{-1}\, v\tens_1 v\tens_1 \extd r + 2\,b\,  \sigma_{QS}(\extd r \tens_1 r^{-1}\, v)\tens_1 v \cr
&=& 2\,b\, r^{-1}\, v\tens_1 (v\tens_1 \extd r +\extd r\tens_1 v)
+ 4\,b^2\,  \lambda\, r^{-1}\, v\tens_1 v\tens_1 v\ ,\cr
\nabla_{QS}(\extd r \tens_1 v) &=& \nabla_{QS}(\extd r) \tens_1 v +(\sigma_{QS}\tens\id)
(\extd r \tens_1 \nabla_{QS}(v) ) \cr
&=& 2\,b\, r^{-1}\, v\tens_1 v \tens_1 v -2\, \sigma_{QS}(\extd r \tens_1   r^{-1}\,v)\tens_1\extd r \cr
&=& 2\, r^{-1}\, v\tens_1 (b\,v \tens_1 v -\extd r \tens_1\extd r  ) -4\,  \lambda\,b\, r^{-1}\, v\tens v \tens_1\extd r \ ,\cr
\nabla_{QS}(v \tens_1 \extd r) &=& \nabla_{QS}(v) \tens_1 \extd r +(\sigma_{QS}\tens\id)
(v \tens_1 \nabla_{QS}(\extd r) ) \cr
&=& -  2\,r^{-1}\,v\tens_1\extd r \tens_1 \extd r + 2\,b\,  \sigma_{QS}(v \tens_1 r^{-1}\, v)\tens_1 v \cr
&=& 2\, r^{-1}\, v\tens_1 (b\,v \tens_1 v -\extd r \tens_1\extd r  ) - 4\,b\,  \lambda\,r^{-1}\,v\tens\extd r\tens_1 v\ ,\cr
\nabla_{QS}(v \tens_1 v) &=& \nabla_{QS}(v) \tens_1 v +(\sigma_{QS}\tens\id)
(v \tens_1 \nabla_{QS}(v) ) \cr
&=& -  2\,r^{-1}\,v\tens_1\extd r \tens_1 v - 2\,\sigma_{QS}(v \tens_1 r^{-1}\,v)\tens_1\extd r \cr
&=& -  2\,r^{-1}\,v\tens_1(\extd r \tens_1 v + v \tens_1\extd r )  + 4\,\lambda\,r^{-1}\,v\tens\extd r\tens_1\extd r \ .
\end{eqnarray*}
Adding these with the appropriate weights gives
\begin{eqnarray*}
\nabla_{QS}(g_1) &=& 4\,b^2\,  \lambda\, r^{-1}\, v\tens_1 v\tens_1 v+ 4\,b\,\lambda\,r^{-1}\,v\tens\extd r\tens_1\extd r \ .
\end{eqnarray*}
Hence the condition for $\nabla_1=\nabla_{QS}+\lambda K$ to preserve $g_1$ is
\begin{eqnarray*}
0 &=& 4\,b^2\,  \lambda\, r^{-1}\, v\tens v\tens v+ 4\,b\,\lambda\,r^{-1}\,v\tens\extd r\tens\extd r \cr
&& +\, \lambda\,(\tau\tens\id)(\extd r \tens K(\extd r))
+ b\, \lambda\,(\tau\tens\id)(v \tens K(v))  \cr
&& +\, \lambda\,K(\extd r)\tens\extd r+ b\,\lambda\,K(v)\tens v\ .
\end{eqnarray*}
Again splitting into the endings, the derivative of $g_1$ is the following $\tens\,\extd r$
\begin{eqnarray*}
 && 4\,b\,r^{-1}\,v\tens\extd r + K(\extd r) \cr
&& +\, (\tau\tens\id)(\extd r \tens (K^r_{rr}\,\extd r+K^r_{vr}\,v))
+ b\, (\tau\tens\id)(v \tens (K^v_{rr}\,\extd r+K^v_{vr}\,v))  \cr
&=& 4\,b\,r^{-1}\,v\tens\extd r + K^r_{vv}\,v\tens v+K^r_{rv}\,\extd r\tens v+K^r_{vr}\,v\tens \extd r+K^r_{rr}\,\extd r\tens \extd r\cr
&& +\, (\tau\tens\id)(\extd r \tens (K^r_{rr}\,\extd r+K^r_{vr}\,v))
+ b\, (\tau\tens\id)(v \tens (K^v_{rr}\,\extd r+K^v_{vr}\,v))   \cr
&=&  (K^r_{vv}+b\,K^v_{vr})\,v\tens v+(K^r_{rv}+b\,K^v_{rr})\,\extd r\tens v+(4\,b\,r^{-1} +2\,K^r_{vr})\,v\tens \extd r+2\,K^r_{rr}\,\extd r\tens \extd r\ ,
\end{eqnarray*}
plus the following $\tens v$
\begin{eqnarray*}
 && 4\,b^2\,  r^{-1}\, v\tens v +  b\,K(v)   \cr
&& +\, (\tau\tens\id)(\extd r \tens (K^r_{rv}\,\extd r+K^r_{vv}\,v))
+ b\, (\tau\tens\id)(v \tens (K^v_{rv}\,\extd r+K^v_{vv}\,v))  \cr
&=& 4\,b^2\,  r^{-1}\, v\tens v +  b\,(K^v_{vv}\,v\tens v+K^v_{rv}\,\extd r\tens v+K^v_{vr}\,v\tens \extd r+K^v_{rr}\,\extd r\tens \extd r)   \cr
&& +\, (\tau\tens\id)(\extd r \tens (K^r_{rv}\,\extd r+K^r_{vv}\,v))
+ b\, (\tau\tens\id)(v \tens (K^v_{rv}\,\extd r+K^v_{vv}\,v))  \cr
&=& (4\,b^2\,  r^{-1} +  2\,b\,K^v_{vv})\,v\tens v+2\,b\,K^v_{rv}\,\extd r\tens v+(b\,K^v_{vr}+K^r_{vv})\,v\tens \extd r+(b\,K^v_{rr}+K^r_{rv})\,\extd r\tens \extd r   \ .
\end{eqnarray*}
Thus the condition for the derivative of $g_1$ to be zero is
\begin{eqnarray*}
K^r_{rr}=K^v_{rv}=0\ ,\  K^r_{vr}=K^v_{vv}=-2\,b\,r^{-1}\ ,\ K^r_{rv}=-b\,K^v_{rr}\ ,\ K^r_{vv}=-b\,K^v_{vr}\ .
\end{eqnarray*}
Combined with the condition for vanishing quantum torsion again gives us our stated values.  \endproof 

One can also check that this quantum connection is indeed the part to order $\lambda$ of the full connection found in \cite{BegMa4} by algebraic methods, provided we make the identification (\ref{2Dcompare}).  In summary, all steps can be made to work in the 2D bicrossproduct model quantum spacetime including a quantum metric $g_1$ and quantisation of the Levi-Civita connection so as to be $*$-preserving and at the same time compatible with $g_1$ and torsion free. That this was possible was not in doubt but we see in detail how it arises at the semiclassical level.

\subsection{Semiquantisation of the Schwarzschild black hole}

We take polar coordinates plus $t$ for 4-dimensional space, where $\phi$ is the angle of rotation about the $z$-axis and $\theta$ is the angle to the $z$-axis.
We take any static isotropic form of metric (including the Schwarzschild case)
\begin{eqnarray}
 g= - e^{N(r)}\extd t\tens\extd t+e^{P(r)}\extd r\tens\extd r+ r^2(\extd\theta\tens\extd\theta+ \sin^2(\theta)\extd\phi\tens\extd\phi)
\end{eqnarray}
The Levi-Civita Christoffel symbols are zero except for
\begin{eqnarray} \label{jyfxxu}
&&\widehat\Gamma^0_{01}=\widehat\Gamma^0_{10}=\tfrac12\,{N'},\quad \widehat\Gamma^1_{11}=\tfrac12\,{P'},\quad \widehat\Gamma^1_{00}=\tfrac12\,{N'}\,e^{N-P}\cr
&& \widehat\Gamma^1_{22}=-r\,e^{-P},\quad \widehat\Gamma^1_{33}=-r\,e^{-P}\sin^2(\theta),\quad \widehat\Gamma^2_{12}=\widehat\Gamma^2_{21}=\widehat\Gamma^3_{13}=\widehat\Gamma^3_{31}=r^{-1}\cr
&&\widehat\Gamma^2_{33}=-\sin(\theta)\cos(\theta),\quad\widehat\Gamma^3_{23}=\widehat\Gamma^3_{32}=\cot(\theta)\ .
\end{eqnarray}

\medskip We shall only consider rotationally invariant Poisson tensors $\omega$. Consider a bivector and rotation invariance in the spherical polar coordinate system. To generate the Lie algebra of the rotation group, we only need two infinitesimal rotations, about the $z$ axis and about the $y$ axis.
For the first, denoting change under the infinitesimal rotation by $\delta$, we get $\delta(\theta)=0$, $\delta(\phi)=1$, and $\delta(\extd\theta)=\delta_A(\extd\phi)=0$.
The infinitesimal rotation about the $y$ axis  is rather more complicated in polar coordinates:
\begin{eqnarray*}
&&\delta(\theta)=\cos\phi\ ,\ \delta(\phi)=-\cot\theta\sin\phi\ ,\ \delta(\extd\theta)=-\,\sin\phi\,\extd\phi\ ,\cr
&& \delta(\extd\phi)=-\,\cot\theta\cos\phi\,\extd \phi+\csc^2\theta\sin\phi\,\extd\theta\ .
\end{eqnarray*}
It is now easily checked that a rotation invariant 2-form on the sphere is, up to a multiple, $\sin\theta\,\extd\theta\wedge\extd\phi$. It follows that a rotation invariant bivector on the sphere is, up to a multiple, given in polars by $\omega^{23}=\csc\theta$.

\begin{proposition}  \label{vbcdhsujhgcx}
If  $\omega$ is rotationally invariant and independent of $x^0$, then only $\omega^{01}=-\omega^{10}=k(r)$ and $\omega^{23}=-\omega^{32}=f(r)/\sin\theta$ are non-zero.
The condition to be a poisson tensor is that $\omega^{01}\,\omega^{23}_{\phantom{23},1}=0$, i.e.\ $k(r)\,f'(r)=0$.
\end{proposition}
\proof  We now suppose that $\omega$ is rotationally invariant as a bivector field. To analyse this is it useful to use our Minkowski-polar coordinates to view $E^i=\omega^{0i}$ as a spatial vector in polar coordinates and to view $\omega^{ij}$ where $i,j\ne 0$ as a spatial 2-form which we view as another vector, $B$. Now consider their values at the north pole of a sphere of radius $r$. Under rotation about the $z$-axis the north pole does not move so there is no orbital angular momentum. There is, however, rotation of the vector indices unless both $E, B$ point along the $z$-axis. This applies equally at any point of the sphere, i.e. $E,B$ must point radially. Equation (\ref{cyclic}) gives the Poisson result.
  \endproof

We now write the Christoffel symbols $\Gamma^a_{bc}$ for the quantising connection $\nabla$ in terms of its torsion $T$ and use Mathematica to get the following result:

\begin{proposition} \label{vhukkhgc} Assume time independence and axial symmetry (i.e.\ that the torsions $T_{ijk}$ are independent of the coordinates $t$ and $\phi$). Then the general solution for the Poisson-compatibility and metric-compatiblity conditions for $(\nabla,\omega)$ is  given by $\omega^{23}=1/\sin\theta$
(up to a constant multiple set to one), $\omega^{01}=0$, and the following restrictions on $T_{ijk}$,
apart from the obvious $T_{ijk}=-T_{ikj}$:

\begin{tabular}{ccc}$T_{012}=T_{201}+T_{102}$ & $T_{013}=T_{301}+T_{103}$ & $T_{023}=0$ \\$T_{123}=0$ & $T_{202}=0$ & $T_{203}=-T_{302}$ \\$T_{212}=r$ & $T_{213}=-T_{312}$ & $T_{223}=0$ \\$T_{303}=0$ & $T_{313}=r\,\sin^2(\theta)$ & $T_{323}=0$\end{tabular}
\end{proposition}

As $T_{313}$ and $T_{212}$ are non-zero, we cannot  take for $\nabla$ the Levi-Civita connection. 
We get the following value of $H^{ij}$, independently of any choice in the torsions: 
\begin{eqnarray*}
H^{ij}\,=\, \left\{\begin{array}{cc}-\,\tfrac12\,\sin\theta\,\extd\theta\wedge\extd\phi & i=j=2  \\\frac12\,\csc\theta\,\extd\theta\wedge\extd\phi & i=j=3  \\0 & \mathrm{otherwise} \end{array}\right.\ .
\end{eqnarray*}
From this  $\CR=g_{ij}\,H^{ij}=0$, so the correction to the metric is zero, $g_1=g_Q$. 

Moreover, we find in Theorem~\ref{qlevi} that (with semicolons WRT the quantising connection) that antisymmetric part of $\nabla_1g_1$ is proportional to 
\begin{eqnarray*}
 \omega^{ij}\,g_{rs}\,S^s_{jn}(R^r{}_{mki}+S^r_{km;i})
&-& \omega^{ij}\,g_{rs}\,S^s_{jm}(R^r{}_{nki}+S^r_{kn;i}) \cr
&=&
\left\{\begin{array}{cc}    -r\,\sin\theta \phantom{bvh} &  (k,m,n)=(2,3,1) \ \&\  (k,m,n)=(3,2,1) \\
  r\,\sin\theta \phantom{bvh} &  (k,m,n)=(2,1,3) \ \&\  (k,m,n)=(3,1,2) 
\\0 & \mathrm{otherwise}\end{array}\right.
\end{eqnarray*}
independently of $\nabla$. Thus there is an obstruction  and no adjustment $\nabla_1$ exactly preserves the metric.

We now specialise to the case where the 
$T_{ijk}$ are rotationally symmetric, which gives the following as the only non-zero torsions,
apart from the obvious $T_{ijk}=-T_{ikj}$:

\begin{tabular}{ccc}$T_{001}=f_1(r)$ & $T_{101}=f_2(r)$ & $T_{203}=-\,T_{302}=-\,f_3(r)\,\sin\theta$ 
\\$T_{212}=r$ & $T_{313}=r\,\sin^2(\theta)$ &  $T_{213}=-\,T_{312}=-\,f_4(r)\,\sin\theta$  
\end{tabular}

\noindent where $f_1(r),f_2(r),f_3(r),f_4(r)$ are arbitrary functions of $r$ only. 

Finallly, we specialise further to the Schwarzschild case, where $e^N=\mathrm{c}^2\,(1-r_s/r)$ and $e^P=(1-r_s/r)^{-1}$, where $r_s$ is the Schwarzschild radius. A short calculation with Mathematica then gives

\begin{lemma} For the Schwarzschild metric the non-zero $R^i{}_{jkl}$, up to the obvious $R^i{}_{jkl}=-R^i{}_{jlk}$ are
\newline
$R^1{}_{010}=R^0{}_{110}=-\,\dfrac{f_1'(r)+\mathrm{c}^2\, r_s\,r^{-3}}{\mathrm{c}^2\,(1-r_s/r)}$
\quad
$R^2{}_{310}=\sin\theta\,(2\,f_3(r)-r\,f_3'(r))\,r^{-3}$
\newline
$R^3{}_{210}=-\csc\theta\,(2\,f_3(r)-r\,f_3'(r))\,r^{-3}$
\quad
$R^3{}_{223}=-\,1$
\quad
$R^2{}_{323}=\sin^2\theta$.

In particular, the curvature cannot vanish entirely.
\end{lemma}

We also have (using row $i$ column $j$ notation)
\begin{eqnarray*}
S^0_{ij} &=& \left(
\begin{array}{cccc}
 0 & -e^{-N} f_1(r) & 0 & 0 \\
 0 & -e^{-N} f_2(r) & 0 & 0 \\
 0 & 0 & 0 & 0 \\
 0 & 0 & 0 & 0 \\
\end{array}
\right)    \ ,\quad
S^1_{ij} = \left(
\begin{array}{cccc}
 -e^{-P} f_1(r) & 0 & 0 & 0 \\
 -e^{-P} f_2(r) & 0 & 0 & 0 \\
 0 & 0 & e^{-P} r & 0 \\
 0 & 0 & 0 & e^{-P} r \sin ^2(\theta )  \\
\end{array}
\right) \ ,  \cr
S^2_{ij} &=& \left(
\begin{array}{cccc}
 0 & 0 & 0 & -\frac{f_3(r) \sin (\theta )}{r^2} \\
 0 & 0 & 0 & -\frac{f_4(r) \sin (\theta )}{r^2} \\
 0 & -\frac{1}{r} & 0 & 0 \\
 0 & 0 & 0 & 0 \\
\end{array}
\right) \ ,\quad 
S^3_{ij} =\left(
\begin{array}{cccc}
 0 & 0 & \frac{f_3(r) \csc (\theta )}{r^2} & 0 \\
 0 & 0 & \frac{\csc (\theta ) f_4(r)}{r^2} & 0 \\
 0 & 0 & 0 & 0 \\
 0 & -\frac{1}{r} & 0 & 0 \\
\end{array}
\right)\ ,
\end{eqnarray*}
and the Christoffel symbols for the quantising connection are are
\begin{eqnarray*}
\Gamma^0_{ij} &=& \left(
\begin{array}{cccc}
 0 & \frac{N'(r)}{2}-e^{-N} f_1(r) & 0 & 0
   \\
 \frac{N'(r)}{2} & -e^{-N} f_2(r) & 0 & 0 \\
 0 & 0 & 0 & 0 \\
 0 & 0 & 0 & 0 \\
\end{array}
\right) \ ,\quad
\Gamma^2_{ij} = \left(
\begin{array}{cccc}
 0 & 0 & 0 & -\frac{f_3(r) \sin (\theta )}{r^2} \\
 0 & 0 & \frac{1}{r} & -\frac{f_4(r) \sin (\theta )}{r^2} \\
 0 & 0 & 0 & 0 \\
 0 & 0 & 0 & -\cos (\theta ) \sin (\theta ) \\
\end{array}
\right) \ ,
\cr
\Gamma^1_{ij} &=& \left(
\begin{array}{cccc}
 \frac{1}{2} e^{-P} \left(e^{N}
   N'(r)-2 f_1(r)\right) & 0 & 0 & 0 \\
 -e^{-P} f_2(r) & -\frac{1}{2} N'(r) & 0 & 0
   \\
 0 & 0 & 0 & 0 \\
 0 & 0 & 0 & 0 \\
\end{array}
\right) \ ,\quad
\Gamma^3_{ij} =\left(
\begin{array}{cccc}
 0 & 0 & \frac{f_3(r) \csc (\theta )}{r^2} & 0 \\
 0 & 0 & \frac{\csc (\theta ) f_4(r)}{r^2} & \frac{1}{r} \\
 0 & 0 & 0 & \cot (\theta ) \\
 0 & 0 & \cot (\theta ) & 0 \\
\end{array}
\right)\ . 
\end{eqnarray*}

Clearly we can chose the functions here to minimise but not eliminate either the torsion or the curvature. 

Now we consider a quantum connection of the form $\nabla_{QS}+\lambda K$ and ask to what extent it preserves the metric, given that we know that it cannot be fully metric-compatible with  $g_1$. 

\begin{proposition}  \label{bioudyvty} For the Schwarzschild metric, the unique $K$ that gives a torsion free quantum connection with the symmetric part of the metric-compatibility tensor vanishing (in Theorem~\ref{qlevi}) coincides with the unique real $K$ that gives a  star-preserving quantum conneciton (in Theorem~\ref{starpresK}). Its 
non-zero components are
\begin{eqnarray*}
K^1_{23}=K^1_{32}=-\,e^{-P}r\,\sin\theta/2\ .
\end{eqnarray*}
\end{proposition}
\proof The unique star-preserving connection (which is also torsion free) is computed by Mathematica from the theorem as stated. For the torsion free connection with vanishing symmetric part of 
the metric compatibility tensor a Mathematica computation gives  $A^p_{nm}=0$ in Lemma~\ref{QSmetrictor} while the only nonzero $B_{knm}$ for the quantity in the proof of Theorem~\ref{qlevi} are
\begin{eqnarray*}
B_{213}=B_{231}=B_{312}=B_{321}=-\,\tfrac12\,r\,\sin\theta\ .
\end{eqnarray*}
Following the proof of Theorem~\ref{qlevi}, the conditions to preserve the symmetric part of the metric and be torsion free are, respectively
\begin{eqnarray*}
K_{nkm}+K_{mkn}=B_{knm}\ ,\quad K_{knm}-K_{kmn}=g_{kp}\,A^p_{nm}\ .
\end{eqnarray*}
leading us in our case to $K_{ijk}=0$ except for
\begin{eqnarray}
K_{123}=K_{132}=-\,\tfrac12\,r\,\sin\theta\ .
\end{eqnarray}
This is the same solution as stated after raising an index. \endproof

This would therefore be a good candidate for `best possible' quantum Levi-Civita' connection as it is suitably `real' (unitary in a suitable context) as expressed by being $*$-preserving and also comes as close as possible to metric-compatible given the order $\lambda$ obstruction. 

Meanwhile, from Proposition~\ref{hermetriccompat} the condition for preserving the corresponding Hermitian metric is $K_{npm}=K_{mpn}$. If we combine this with quantum torsion freeness (as in the proof above), the condition for torsion free Hermitian-metric compatibiity is that $K_{ijk}$ is totally symmetric in $ijk$. This means that there is a 4-functional moduli of Hermitian-metric compatible connections (but  necessarily disjoint from the point in Proposition \ref{bioudyvty}, since that point is star-preserving but not fully metric compatible).  Also, from (\ref{kjhvcxjhgc}), in the case where $\mathcal{R}=0$, the cotorsion (as in the approach used in \cite{Ma:sph})  is
\begin{eqnarray*}  
(\wedge\tens\id)\,q^2\,\nabla_{1}(g_1)& =& -\,\lambda\,\omega^{ij}\,g_{rs}\,S^s_{jn}(R^r{}_{mki}+S^r_{km;i})\,\extd x^k\wedge\extd x^m \tens \extd x^n           \cr
&&+\, \lambda\,(  g_{pn}\,K^p_{km}+g_{mp}\,K^p_{kn}   )\, \extd x^k\wedge\extd x^m \tens \extd x^n\  \ ,
\end{eqnarray*}
and from this the condition for vanishing cotorsion comes down to $ K_{nkm}+K_{mkn} - K_{nmk}-K_{kmn} =0 $. Adding to this the condition for vanishing torsion, $ K_{nkm}= K_{nmk} $, we see that $ K_{mkn} =K_{kmn}  $, so the resulting $K_{kmn}$ are totally symmetric. Thus the space of torsion free Hermitian-metric compatible $\nabla_1$  in the case of the black-hole is the same as the space of torsion and cotorsion free ones.

In summary we have shown that we inevitably have curvature of $\nabla$ and hence a nonassociative calculus  at order $\lambda^2$ if we try to quantise the black-hole and keep rotational invariance, an anomaly in line with experience in quantum group models\cite{BegMa1}. As with those models, the alternative is to quantise associatively but have an extra cotangent dimension as in the wave-operator quantisation of the black hole recently achieved to all orders in \cite{Ma:alm}. And by perturbing $\nabla_{QS}$ we can either have a unique `unitary' connection in the sense of star-preserving (and best possible but not fully metric compatible) or a moduli of Hermitian-metric compatible (but not star-preserving) connections perhaps fixed by further requirements. Which type of $K$ to take would depend on the application or experimental assumption as to which feature we want to fix in the quantisation process.

\end{document}